\pdfoutput=1
\RequirePackage{ifpdf}
\ifpdf 
\documentclass[pdftex]{sigma}
\else
\documentclass{sigma}
\fi

\numberwithin{equation}{section}

\newtheorem{Theorem}{Theorem}[section]
\newtheorem*{Theorem*}{Theorem}
\newtheorem{Corollary}[Theorem]{Corollary}
\newtheorem{Lemma}[Theorem]{Lemma}
\newtheorem{Proposition}[Theorem]{Proposition}

\theoremstyle{definition}
\newtheorem{Definition}[Theorem]{Definition}

\newtheorem{Example}[Theorem]{Example}
\newtheorem{Remark}[Theorem]{Remark}

\usepackage{tikz-cd}
\usepackage{paralist}
\usepackage{xstring}

\tikzset{
curvarr/.style={
 to path={ -- ([xshift=2ex]\tikztostart.east)
 |- (#1) [near end]\tikztonodes
 -| ([xshift=-2ex]\tikztotarget.west)
 -- (\tikztotarget)}
 }
}

\newcommand{\lie}[1]{\operatorname{\mathfrak{#1}}}

\newcommand{\sln}{\lie{sl}}

\newcommand\C{{\mathbb C}}

\newcommand{\R}{{\mathbb R}}

\renewcommand{\i}{\mathsf{i}}

\newcommand{\pder}[2]{
		\IfEqCase{#1}{
			{}{\frac{\partial}{\partial {#2}}}
			{t}{\frac{\partial}{\partial {#2}}}
			}}

\begin{document}

\allowdisplaybreaks

\newcommand{\arXivNumber}{2102.10853}

\renewcommand{\PaperNumber}{008}

\FirstPageHeading

\ShortArticleName{Geometry of the Space of Sections of Twistor Spaces with Circle Action}

\ArticleName{Geometry of the Space of Sections of Twistor Spaces \\ with Circle Action}

\Author{Florian BECK~$^{\rm a}$, Indranil BISWAS~$^{\rm b}$, Sebastian HELLER~$^{\rm c}$ and Markus R\"OSER~$^{\rm a}$}

\AuthorNameForHeading{F.~Beck, I.~Biswas, S.~Heller and M.~R\"oser}

\Address{$^{\rm a)}$~Fachbereich Mathematik, Universit\"at Hamburg, 20146 Hamburg, Germany}
\EmailD{\mail{flrn.beck@gmail.com}, \mail{markus.roeser@uni-hamburg.de}}

\Address{$^{\rm b)}$~School of Mathematics, Tata Institute of Fundamental Research,\\
\hphantom{$^{\rm b)}$}~Homi Bhabha Road, Mumbai 400005, India}
\EmailD{\mail{indranil@math.tifr.res.in}}

\Address{$^{\rm c)}$~Beijing Institute of Mathematical Sciences and Applications,\\
\hphantom{$^{\rm c)}$}~No.~544, Hefangkou Village Huaibei Town, Beijing 101408, P.R.~China}
\EmailD{\mail{sheller@bimsa.cn}}

\ArticleDates{Received March 25, 2025, in final form January 13, 2026; Published online February 02, 2026}

\Abstract{We study the holomorphic symplectic geometry of (the smooth locus of) the space of holomorphic sections of a twistor space with rotating circle action. The twistor space carries a line bundle with meromorphic connection constructed by Hitchin. We give an interpretation of Hitchin's meromorphic connection in the context of the Atiyah--Ward transform of the corresponding hyperholomorphic line bundle. It is shown that the residue of the meromorphic connection serves as a moment map for the induced circle action, and furthermore the critical points of this moment map are studied. Particular emphasis is given to the example of Deligne--Hitchin moduli spaces.}


\Classification{53C26; 53C28; 53C43; 14H60; 14H70}

\section{Introduction}

Given a compact Riemann surface $\Sigma$ of genus at least two and a complex semisimple Lie group $G$, the moduli
space of $G$-Higgs bundles on $\Sigma$ has a hyperK\"ahler structure \cite{Hitchin87}. The corresponding twistor
space $\mathcal{M}_{\mathrm{DH}}(\Sigma, G)$ is known as the Deligne--Hitchin moduli space \cite{Si-Hodge}. This space $\mathcal{M}_{\mathrm{DH}}(\Sigma, G)$ has been the topic of many papers in recent years; see, for example, \cite{CollierWentworth,DFKMMN,DM,HuHuang,Huang} and \cite{BeHeRo,BH,BHR,HelHel,Hel} and references therein.

In \cite{BHR}, the last three
authors initiated a detailed investigation of the space $\mathcal S$ of holomorphic
sections of the natural holomorphic projection $\mathcal{M}_{\mathrm{DH}}(\Sigma, G) \longrightarrow \C P^1$. Our aim here is
to go deeper into the systematic study of $\mathcal S$. There are two distinct sources
of motivation for doing this: One is the integrable systems approach to the theory of harmonic maps
from Riemann surfaces into symmetric spaces, the second one coming from hyperK\"ahler geometry.

Let $K$ be a compact real form of a complex semisimple Lie group
$G$ with Lie algebra $\mathfrak{g}$. Hitchin's self-duality equations on $\Sigma$ are given by~\cite{Hitchin87}
\begin{equation}\label{eq:SDintro}
F^\nabla + [\Phi\wedge\Phi^*] = 0 = {\overline{\partial}}^\nabla\Phi.
\end{equation}
Here $\nabla$ is a $K$-connection and $\Phi \in \Omega^{1,0}(\Sigma, \mathfrak{g})$, which is called a Higgs
field, while $\Phi \longmapsto \Phi^*$ is given by the negative of the Cartan involution
on $\mathfrak{g}$ associated with the given real form~${K \subset G}$. Every solution of \eqref{eq:SDintro} gives rise to a family of flat $G$-connections
on $\Sigma$ parametrized by $\lambda\in \C^* = \C \setminus\{0\}$, $
\nabla^\lambda = \nabla + \lambda^{-1}\Phi + \lambda\Phi^*
$
that satisfies a certain reality condition which is determined by the symmetric pair $(G, K)$.
The nonabelian Hodge correspondence \cite{Corlette,Don,Hitchin87, SimpsonIHES1992} allows us to invert
this process, meaning any (reductive) flat $G$-connection on $\Sigma$ actually
occurs in a~suitable $\C^*$-family of flat $G$-connections
given by a solution of the self-duality equations \eqref{eq:SDintro}.

Recall from \cite{Si-Hodge} that the Deligne--Hitchin moduli space $\mathcal{M}_{\mathrm{DH}}(\Sigma, G)$ is obtained by gluing the
two Hodge moduli spaces $\mathcal{M}_{\mathrm{Hod}}(\Sigma, G)$ and $\mathcal{M}_{\mathrm{Hod}}\bigl(\overline{\Sigma}, G\bigr)$ of $\lambda$-connections via the Riemann--Hilbert correspondence, where $\overline{\Sigma}$ is the conjugate of
$\Sigma$ (see Section~\ref{s:DeligneHitchin} for details). It enters the picture via the interpretation
of the $\C^*$-family of flat $G$-connections associated with a solution of \eqref{eq:SDintro} as a family of
$\lambda$-connections, which extends to all of $\C P^1$. From this point of view, we may think of the above
$\C^*$-family as a holomorphic section of the fibration $\mathcal{M}_{\mathrm{DH}}(\Sigma, G) \longrightarrow \C P^1$ satisfying an
appropriate reality condition. This observation fits naturally into the twistor theory of hyperK\"ahler
manifolds, as we shall explain next.

The moduli space $\mathcal{M}_{\mathrm{SD}}(\Sigma, G)$ of solutions of \eqref{eq:SDintro} is a (typically singular) hyperK\"ahler mani\-fold which
comes equipped with a rotating circle action given by $(\nabla, \Phi) \longmapsto
\bigl(\nabla, {\rm e}^{{\rm i}\theta}\Phi\bigr)$ for all~${\theta \in \mathbb R}$. Here the word
\emph{rotating} means that the circle action is isometric, it preserves one of the three K\"ahler forms $\omega_I$, $\omega_J$, $\omega_K$, let us say $\omega_I$,
while rotating
the other two in the sense that it acts with weight one on the holomorphic symplectic form $\omega_J+\i \omega_K$ (see Sections~\ref{ss:rotating} and~\ref{ss: MSD} for details). The Deligne--Hitchin moduli
space $\mathcal{M}_{\mathrm{DH}}(\Sigma, G)$ is the twistor space associated with~$\mathcal{M}_{\mathrm{SD}}(\Sigma, G)$. The $\C^*$-family of flat
connections on $\Sigma$ associated with a solution to the self-duality equations
in \eqref{eq:SDintro} can then be interpreted as the twistor line corresponding to the point in~$\mathcal{M}_{\mathrm{SD}}(\Sigma, G)$
represented by $(\nabla, \Phi)$ in \eqref{eq:SDintro}. In this way, we can view $\mathcal{M}_{\mathrm{SD}}(\Sigma, G)$ in a natural
way as a subset of the space $\mathcal S$ of holomorphic sections of the twistor family~$\mathcal{M}_{\mathrm{DH}}(\Sigma, G) \longrightarrow
\C P^1$.

More generally, if $Z$ is the twistor space of a hyperK\"ahler manifold $M$, then the twistor construction produces a hyperK\"ahler metric on the space $M'$ of all real holomorphic sections (with
appropriate normal bundle) of the fibration $Z \longrightarrow \C P^1$. From this point of view, it is a natural question
whether $M = M'$ \cite{Si-Hodge}. In the context of the self-duality equations, this motivates the question whether every real holomorphic section of $\mathcal{M}_{\mathrm{DH}}(\Sigma,G) \longrightarrow \C P^1$ is actually obtained from a solution to
the self-duality
equations. In \cite{BHR}, an answer to this question was given in the case of $G = \operatorname{SL}(2,\C)$ by showing that in
general $\mathcal{M}_{\mathrm{SD}}(\Sigma, G)$ is strictly contained in the space of real sections. Furthermore, in \cite{BHR} the real
holomorphic sections belonging to $\mathcal{M}_{\mathrm{SD}}(\Sigma, G)$ were characterized in terms of a certain $\mathbb{Z}/2\mathbb{Z}$-valued
invariant that can be associated with any real section.

The starting point for the results of \cite{BHR} is the well-known fact that Hitchin's self-duality equations \eqref{eq:SDintro} can be
interpreted as the gauge-theoretic equations for a (twisted) harmonic map $\Sigma \longrightarrow G/K$
(a so-called \emph{harmonic metric}).
One observes that the twisted harmonic maps into other (pseudo-Riemannian) symmetric spaces of the form $G/G_{\R}$
and their duals (where $G_\R$ is a real form of $G$) have analogous gauge-theoretic interpretations and give rise to $\C^*$-families of flat
$G$-connections satisfying different reality conditions depending on the target of the harmonic map.
For example, if we have $G = \operatorname{SL}(2,\C)$ then $K = {\rm SU}(2)$, and the harmonic maps, from~$\Sigma$, into the symmetric spaces
$\operatorname{SL}(2,\C)/{\rm SU}(2)$, $\operatorname{SL}(2,\C)/{\rm SU}(1,1)$, ${\rm SU}(2)$, ${\rm SU}(1,1)$ as well as into the hyperbolic disc
${\rm SU}(1,1)/{\rm U}(1)$ can all be encoded in $\C^*$-families of flat $\operatorname{SL}(2,\C)$ connections on~$\Sigma$ satisfying appropriate
reality conditions; see for instance \cite{DPW,Hitchin90,Pohl,Uh}.
Considering these $\C^*$-families as families of $\lambda$-connections on $\Sigma$ with $\lambda \in \C P^1$ we can --
at least formally -- interpret twisted harmonic maps from $\Sigma$ into various symmetric spaces associated
with real forms of $G$ as holomorphic sections of $\mathcal{M}_{\mathrm{DH}}(\Sigma, G)$; see~\cite{BHR}.

This interaction with the theory of harmonic maps has provided useful guidance and intuition for the
predecessors~\cite{BHR} and \cite{BeHeRo} of this project. For example, the new real holomorphic sections for the
$\operatorname{SL}(2,\C)$-case constructed in~\cite{BHR} are closely related to twisted harmonic maps~${\Sigma \longrightarrow
\operatorname{SL}(2,\C)/{\rm SU}(1,1)}$.
We may view the
space $\mathcal S$ of holomorphic sections of the fibration $\mathcal{M}_{\mathrm{DH}}(\Sigma,G) \longrightarrow \C P^1$ as a \emph{master space}
for the moduli spaces of twisted harmonic maps from $\Sigma$ into various (pseudo-Riemannian) symmetric spaces
associated with the group $G$. These moduli spaces are contained in $\mathcal S$ as the fixed point loci of certain antiholomorphic involutions. It is thus natural to study the global geometry of $\mathcal S$.

Many of our results and constructions are best understood from the point of view of the general twistor theory of hyperK\"ahler manifolds with a rotating circle action. These are an active field of research in hyperK\"ahler geometry,
especially their connection with the well-known hyperK\"ahler/quaternion K\"ahler correspondence is extensively investigated \cite{ACM,Haydys,Hitchin-HKQK,Hitchin-HKQK2}. Such manifolds are never compact, a main source of examples is the construction of Feix \cite{Fei01} (and independently Kaledin \cite{Kaledin}) of a hyperK\"ahler metric on (a neighbourhood of the zero section in) the cotangent bundle of a K\"ahler manifold, where the circle action is just the usual action by scalar multiplication in the fibers. Metrics of this type also arise on many moduli spaces of solutions to gauge theoretic equations (e.g., magnetic monopoles, Higgs bundles). See the papers~\cite{Hitchin-HKQK,Hitchin-HKQK2} for a discussion of examples of interest to physicists as well as of a purely geometric or Lie theoretic origin.

HyperK\"ahler manifolds with rotating circle action feature naturally in the hyperK\"ah\-ler/quater\-nion K\"ahler (HK/QK) correspondence. Haydys \cite{Haydys} has observed the following: If~$(M, \omega_I, \omega_J,
\omega_K)$ is a hyperK\"ahler manifold with rotating circle action such that $\omega_I$ is integral, then there
exists a complex line bundle $L \longrightarrow M$ with a unitary connection whose curvature is~${\omega_I + dd^c_I\mu}$, where $\mu$ is
the moment map for the $S^1$-action with respect to the K\"ahler form~$\omega_I$. This unitary connection on $L$ is
hyperholomorphic, in the sense that its curvature is of type $(1, 1)$ with respect to every complex structure in the
family of K\"ahler structures on $M$ parametrized by the sphere $S^2$. The moment map
facilitates a lift of the $S^1$-action on $M$ to an $S^1$-action on
the total space of the principal $S^1$-bundle $P \longrightarrow M$ associated with $L$. Haydys
has shown that the corresponding
quotient $Q = P/S^1$ for this action carries a quaternionic K\"ahler metric. The natural $S^1$-action on the
principal bundle $S^1$-bundle $P$ descends to an
isometric circle action on $Q$, and $M$ can be recovered as a hyperK\"ahler quotient of the Swann bundle of $Q$ by
the lift of this isometric circle action.
Hitchin \cite{Hitchin-HKQK,Hitchin-HKQK2} has described the HK/QK correspondence from a
purely twistorial point of view. A natural starting point for him is the observation that, by the Atiyah--Ward
correspondence, the hyperholomorphic line bundle on $M$ produces a~holomorphic line bundle $L_Z$ on the twistor
space $Z$ associated to $M$. Hitchin has shown how to construct $L_Z$ directly on the twistor space $Z$ and explained how the
circle action determines a distinguished meromorphic connection on $L_Z$. This meromorphic connection plays a key
role in the construction, from $Z$, of the twistor space $Z_Q$ of the associated quaternionic K\"ahler manifold
$Q$ in the sense that it determines the contact distribution on $Z_Q$. Therefore, in order to make progress towards the
 understanding of the quaternionic K\"ahler manifold associated with $\mathcal{M}_{\mathrm{SD}}(\Sigma,G)$, it is clearly important to obtain
information on this meromorphic connection.

In \cite{BeHeRo}, it is proved that on the space $\mathcal S$ of holomorphic sections of the twistor family
$\mathcal{M}_{\mathrm{DH}}(\Sigma, G) \longrightarrow \C P^1$ there exists an interesting
and useful holomorphic functional, called the energy functional, whose
evaluation on a real section coming from a harmonic map is (a~constant multiple of) the Dirichlet energy of the
associated harmonic map. This functional is also intimately related to the Willmore energy of certain immersions of
$\Sigma$ into the 3-sphere. It was then shown that this functional has a natural interpretation from the
hyperK\"ahler point of view. In fact, it can be interpreted as a holomorphic extension of the moment map,
associated with the rotating circle action, from the space~$M$ of twistor lines to the whole of $\mathcal S$.
Essentially, it is given by associating to any~$s \in \mathcal S$ the residue of the meromorphic connections along
$s$.

In this article, we continue our study of this setup in the context of the twistor space $Z$ of a~general
hyperK\"ahler manifold $M$ with a rotating circle action. The natural geometric structure on the space $\mathcal
S$ of holomorphic sections of $Z \longrightarrow \C P^1$, viewed as a complexification of the hyperK\"ahler
manifold $M$, has been elucidated by Jardim and Verbitsky \cite{JV2,JV1}. They show that $\mathcal S$ comes
equipped with a certain family of closed two-forms called a trisymplectic structure. In particular, the
complexifications of the K\"ahler forms on $M$ are contained in this family. We describe in this paper how the
picture gets enriched by the presence of a~rotating circle action. It turns out that the energy functional
$\mathcal E$ is a~moment map for the circle action on $\mathcal S$ induced from the circle action on $M$ with
respect to a natural holomorphic symplectic form $\Omega_0$ which is a~part of the canonical trisymplectic
structure on $\mathcal S$ (see Theorem~\ref{EnergyMomentMap}). In~particular, the critical points of $\mathcal E$
turn out to be exactly the fixed points of the circle action.\looseness=1

We also explain how to use the Atiyah--Ward transform to obtain a natural holomorphic extension of the line
bundle $L \longrightarrow M$ to $\mathcal L \longrightarrow \mathcal S$, and how the meromorphic $1$-form given by the difference between the meromorphic connection and the Atiyah--Ward connection
on $L_Z$ can be described naturally in terms of $\mathcal L$ and the data of the circle action (see Theorem
\ref{AWH} and Section~\ref{subsec:alternativeconn}). This is implicit in Hitchin's work \cite{Hitchin-HKQK}, but we believe our point of view
in terms of the geometry of the space $\mathcal S$ of holomorphic sections might shed a new light onto his
constructions and will be useful in explicit computations. Moreover, our results allow us to extend the knowledge about the global structure of $\mathcal S$ (see
Theorem~\ref{Thm:S0neqZ0Zinfty}).

We then apply this general framework to the space $\mathcal S$ of holomorphic sections of $\mathcal{M}_{\mathrm{DH}}(\Sigma,\allowbreak \operatorname{SL}(n,\C))$ and are able to strengthen some of our general results in this particular setup. Even
though $\mathcal S$ is expected to be singular in this case, the general theory tells us that the critical points
of the energy functional are actually closely related to the fixed points of the $\C^*$-action on $\mathcal S$. These sections arise, for instance, in the context of variations of Hodge structures. There are interesting examples of such sections obtained from $2\pi$-grafting constructions (\cite{Hel}, see also Example \ref{ex:grafting}). We analyze
these $\C^*$-invariant sections in detail (see Proposition~\ref{Prop:C*_sectionLiftC}), building on~\cite{CollierWentworth}. We obtain explicit formulas for the energy of a $\C^*$-invariant section
(see Proposition~\ref{Prop: energyC*-fixed_sections}) and we describe the second variation of the energy at such a
section (see Proposition~\ref{prop1}). Moreover, we establish an explicit formula for the degree of the hyperholomorphic
line bundle restricted to a $\C^*$-invariant section (see Proposition~\ref{Prop: deghyperhol_C*fixedsection}). We hope such explicit formulas will be useful in future investigations, for example of the grafting sections mentioned above. We~apply these results at the end of the article to study the global topology of $\mathcal S$. By showing that there exist sections along which the hyperholomorphic line bundle has non-zero
degree, we are able to deduce that $\mathcal S$ is in general not connected (see Theorem~\ref{Thm:
Hyperholnontrivial}).

The paper is organized as follows. In Section~\ref{sectionaboutsections}, we discuss the twistor fibration $\varpi\colon Z
 \longrightarrow \C P^1$ associated with
a hyperK\"ahler manifold $M$ and describe the geometric structure induced on the space $\mathcal S$ of
holomorphic sections of $\varpi$ in a way that is suitable for our purpose. We then go on in Section~\ref{sectionhyperholo} to explain how
the geometry of~$\mathcal S$ is enriched by the presence of a~rotating circle action on $M$ and hence on $Z$. In
particular, we explain how the energy functional~$\mathcal E$ ties up naturally with the holomorphic symplectic
geometry on $\mathcal S$ discussed in Section~\ref{sectionaboutsections} and with the meromorphic connection on the
hyperholomorphic line bundle. We have added a~large amount of background material in Sections~\ref{sectionaboutsections} and~\ref{sectionhyperholo} to keep the
paper somewhat self-contained. The particular example of the Deligne--Hitchin moduli space is described in Section~\ref{s:DeligneHitchin}. In
Section~\ref{ss:GeometrySectionDH}, the abstract geometric framework developed in the first two sections is illustrated in the context of
the Deligne--Hitchin moduli space $\mathcal{M}_{\mathrm{DH}}(\Sigma, G)$ and the results mentioned above are proved.\looseness=1

\section[Geometry of the space of holomorphic sections of a twistor space]{Geometry of the space of holomorphic sections\\ of a twistor space}\label{sectionaboutsections}

In this section, we collect some aspects of the geometry of holomorphic sections of a hyperK\"ahler twistor space
that will be used later. Useful references are \cite{HKLR,JV2,JV1}. See also \cite{LeBrun} for
the similar, but different, quaternionic K\"ahler case.

\subsection{Twistor space}\label{ss:twistor}

Let $(M, g_0, I, J, K)$ be a hyperK\"ahler manifold of complex dimension $2d$,
where $I$, $J$, $K$ are almost complex structures and $g_0$ a Riemannian metric on $M$.
The associated K\"ahler forms are~${\omega_L = g_0(L\cdot, \cdot)}$, $ L \in \{I, J, K\}$.
For convenience, the complex manifolds $(M, I)$ and
$(M, -I)$ will sometimes be denoted by simply $M$ and $\overline{M}$, respectively; it will
be ensured that this abuse of notation does not create any confusion.

There is a family of K\"ahler structures on $M$ with complex structures
\begin{equation}\label{e2}
\bigl\{I_x := x_1I + x_2J+x_3K \mid x := (x_1, x_2, x_3) \in S^2 \bigr\}
\end{equation}
parametrized by the sphere
$S^2 := \bigl\{(x_1, x_2, x_3) \in {\mathbb R}^3 \mid x^2_1+x^2_2+x^2_3 = 1\bigr\}$.
The twistor space~${Z = Z(M)}$ of $(M, g_0, I, J, K)$ is a complex manifold whose underlying smooth
manifold is~${S^2\times M}$ \cite[Section~3\,(F)]{HKLR}.
The almost complex structure $I_Z$ of $Z$ at any point~$(x, m) \in S^2\times M$~is
\[I_Z|_{(x,m)} = (I_{\mathbb CP^1}\vert_{T_x S^2}) \oplus (I_x\vert_{T_mM}) ,\]
 where
$I_{\mathbb CP^1}$ is the standard almost complex structure on $S^2 = \C P^1$ and $I_x$
is the almost complex structure in \eqref{e2}. Here
we identify $S^2$ with $\C P^1$ using the stereographic projection from~${(-1, 0, 0)}$ to the plane in
${\mathbb R}^3$ spanned by the $x_2$ axis and the $x_3$ axis.
In particular, $(1, 0, 0) \in S^2$ corresponds to $0 \in \C P^1$ and $(-1, 0, 0)$ corresponds to
$\infty \in \C P^1$.
Throughout we shall use an affine coordinate $\lambda$ on $\C P^1$, so that $\C P^1 \cong \C\cup\{\infty\}$
having two coordinate functions~$\lambda$ and $\lambda^{-1}$.

The hyperK\"ahler structure on $M$ is encoded in the following complex-geometric data on $Z$.
The natural projection $S^2 \times M \longrightarrow S^2$ corresponds to a holomorphic submersion
\begin{equation}\label{evp}
\varpi \colon\ Z \longrightarrow \C P^1
\end{equation}
with fibers
\begin{equation}\label{zl}
\varpi^{-1}(\lambda) =: Z_\lambda = (M, I_\lambda) .
\end{equation}
Here $I_\lambda$ is the complex structure on $M$ in \eqref{e2} corresponding to $\lambda \in \C P^1 \cong S^2$.
For any complex vector bundle $V$ on $Z$ and any integer $m$, we use the notation $V(m) := V \otimes \varpi^*\mathcal{O}_{\C P^1}(m)$.
Let
\[T_\varpi = T_\varpi Z := (\ker {\rm d}\varpi) \subset TZ\]
be the relative holomorphic tangent bundle (also called
the vertical tangent bundle) for the projection $\varpi$ in \eqref{evp}.
Then $Z$ carries the twisted relative holomorphic symplectic form
$\omega \in H^0\bigl(Z, \bigl(\Lambda^2 T_{\varpi}^*\bigr)(2)\bigr)$ given by
\begin{equation}\label{t6}
\omega = \left(\omega_J+\i \omega_K + 2\i \lambda\omega_I +
\lambda^2(\omega_J-\i\omega_K)\right)\otimes\frac{\partial}{\partial \lambda} \in
H^0\bigl(Z, \bigl(\Lambda^2T_\varpi^*\bigr)(2)\bigr) .
\end{equation}
Moreover, $Z$ carries an anti-holomorphic involution (or real structure)
\begin{equation}\label{etz}
\tau_Z \colon\ Z \longrightarrow Z
\end{equation}
given by the map $S^2 \times M \longrightarrow S^2 \times M$, $ (x, m) \longmapsto (-x, m)$,
using the diffeomorphism $Z \cong S^2\times M$. Note that $\tau_Z$ covers the antipodal map
\begin{equation}\label{tcp}
\tau_{\C P^1}\colon\ \C P^1 \longrightarrow \C P^1, \qquad \lambda \longmapsto -(\overline{\lambda})^{-1},
\end{equation}
which is an anti-holomorphic involution.
The relative twisted symplectic form $\omega$ is real with respect to $\tau$ in the sense that $\tau^*\overline\omega
 = \omega$.

Let $\mathcal{S}$ denote the space of all holomorphic sections of
the projection $\varpi$ in \eqref{evp}; this space is discussed in detail in Section \ref{ss:HolSects}.
We have the embedding
\begin{equation}\label{te}
\iota \colon \ M \hookrightarrow \mathcal{S}
\end{equation}
that sends any $m \in M$ to the constant section $x \longmapsto (x, m) \in S^2\times M$,
which is called the twistor line $s_m$ associated to $m$.
By \cite[Section~3, (F)]{HKLR} it is known that the normal bundle of a~twistor line is isomorphic to $\mathcal{O}_{\C P^1}(1)^{\oplus 2d}$.
The space $\mathcal{S}$ has a real structure defined by
\begin{equation}\label{tau}
\tau \colon \ \mathcal{S} \longrightarrow \mathcal{S}, \qquad s \longmapsto \tau_Z\circ s\circ\tau_{\C P^1} ,
\end{equation}
where $\tau_Z$ and $\tau_{\C P^1}$ are the involutions in \eqref{etz} and \eqref{tcp},
respectively; we note that while the section $\tau(s)$ is holomorphic for fixed $s$, the map
$\tau$ itself is anti-holomorphic. The manifold $M$, considered as space of twistor lines, is a component of the
fixed point locus $\mathcal{S}^\tau \subset \mathcal{S}$, also known as the space of real sections.

We have described the complex-geometric data on $Z$ induced by the hyperK\"ahler structure on $M$. Conversely,
suppose we have a complex manifold $Z$ with a holomorphic submersion~${\varpi\colon Z
\longrightarrow \C P^1}$, a twisted relative
symplectic form~${\omega\in H^0\bigl(Z, \bigl(\Lambda^2T^*_\varpi\bigr)(2)\bigr)}$ and an antiholomorphic involution $\tau\colon Z
 \longrightarrow Z$
covering the antipodal map on $\C P^1$ satisfying $\tau^*\overline{\omega} = \omega$. Then the parameter space $M$ of real
sections of $\varpi$ with normal bundle isomorphic to $\mathcal{O}_{\C P^1}(1)^{\oplus 2d}$ is a~(pseudo-)hyperK\"ahler
manifold of real dimension $4d$, if it is non-empty; see \cite[Theorem 3.3]{HKLR}.

\subsection{Space of holomorphic sections as a complexified hyperK\"ahler manifold}\label{ss:HolSects}

To examine the local structure of $\mathcal{S}$, for any $s \in {\mathcal S}$, let
$N_s = (s^*TZ) / {\rm d}s\bigl(T\C P^1\bigr)$
be the normal bundle
of $s\bigl(\C P^1\bigr) \subset Z$. Since $s$ is a section of $\varpi$, we have a canonical isomorphism~${
s^*T_\varpi Z \cong N_s }$,
where $T_\varpi Z \subset TZ$ as before is the kernel of the differential ${\rm d}\varpi$ of the map~$\varpi$. The
following proposition is well known, however a proof of it is given because parts of the proof will be
used later.

\begin{Proposition}\label{prop:SpaceOfSections}
Let $Z$ be the twistor space of a hyperK\"ahler manifold $M$ of complex dimension~$2d$. Then the set $\mathcal{S}$ of
holomorphic sections of
$\varpi\colon Z \longrightarrow \C P^1$
is a~complex space in a~natural way. The
tangent space of $s \in
\mathcal{S}$ is $H^0\bigl({\mathbb C}P^1, N_s\bigr)$, and $\mathcal{S}$
is smooth at a point $s \in \mathcal{S}$ if~${H^1\bigl({\mathbb C}P^1, N_s\bigr) = 0}$. If
$H^1\bigl({\mathbb C}P^1, N_s\bigr) = 0$, then $\dim T_s\mathcal S = 4d$.
\end{Proposition}

\begin{proof}
For any $s \in \mathcal{S}$, the sufficiently small deformations of the complex submanifold
${s\bigl({\mathbb C}P^1\bigr) \subset Z}$ continue to be the image of a section of $\varpi$. Consequently, $\mathcal{S}$
is an open subset of the corresponding Douady space of rational curves in $Z$; see also \cite[Theorem 2]{Namba}. In particular,
$T_s\mathcal{S} = H^0\bigl({\mathbb C}P^1, N_s\bigr)$ for all $s \in \mathcal{S}$.

The complex structure of $\mathcal{S}$ around a point $s \in \mathcal{S}$ is constructed in the
following way. There are open neighbourhoods
$V_s \subset H^0\bigl({\mathbb C}P^1, N_s\bigr)$ and $
U_s \subset \mathcal{S}$
 of $0 \in H^0\bigl({\mathbb C}P^1, N_s\bigr)$ and $s \in \mathcal{S}$ respectively,
as well as a holomorphic map
\[
\mathsf{k}_s \colon \ H^0\bigl({\mathbb C}P^1, N_s\bigr) \longrightarrow H^1\bigl({\mathbb C}P^1, N_s\bigr) ,
\]
which is sometimes called the Kuranishi map.
Then there is a natural isomorphism
\begin{equation}\label{eq:kuranishi}
U_s \cong V_s\cap \mathsf{k}_s^{-1}(0)
\end{equation}
taking $s$ to $0$. The complex structure on $U_s$ is given by the complex structure of $V_s\cap \mathsf{k}_s^{-1}(0)$
using the isomorphism in \eqref{eq:kuranishi}.

The statement on the smoothness of $\mathcal S$ follows from \eqref{eq:kuranishi}, because $\mathsf{k}_s$ is the
zero map if~${H^1\bigl({\mathbb C}P^1, N_s\bigr) = 0}$. The dimension of $T_s{\mathcal S} \cong
H^0\bigl({\mathbb C}P^1, N_s\bigr)$ is computed from the Riemann--Roch theorem applied to $N_s$
\begin{equation}\label{b1}
h^0\bigl(\C P^1, N_s\bigr) = \deg(N_s) + \mathrm{rank}(N_s) = 2d + 2d = 4d.
\end{equation}
Here we have used the fact that we have $\deg N_s = 2d$ for any holomorphic section $s\colon \C P^1
\longrightarrow Z$ (cf.\ \cite{Fei01,Mayrand}).
To see the equality in \eqref{b1}, note that we have $N_s = s^*T_\varpi Z$ and hence the twisted relative symplectic form
$\omega \in H^0\bigl(Z, \bigl(\Lambda^2T_\varpi^*\bigr)(2)\bigr)$ induces an isomorphism $s^*\omega\colon N_s
\longrightarrow N_s^*\otimes \mathcal{O}_{\C P^1}(2)$; consequently, we have
$(\det N_s)^{\otimes 2} \cong \mathcal{O}_{\C P^1}(4d)$.
\end{proof}

As mentioned earlier, the normal bundle $N_{s_m}$ of any twistor line $s_m$, $m \in M$, is isomorphic to
${\mathcal O}_{\C P^1}(1)^{\oplus 2d}$ so that we have $H^1\bigl(\C P^1, N_{s_m}\bigr) = 0$.
Consequently, by Proposition \ref{prop:SpaceOfSections}, the image~${M \subset \mathcal{S}}$ of the embedding in \eqref{te}
is contained in the smooth locus of $\mathcal S$. Define the complex manifold
\begin{equation}\label{t3}
{\mathcal S}' := \bigl\{s \in \mathcal S \mid N_s \cong
{\mathcal O}_{\C P^1}(1)^{\oplus 2d}\bigr\} \subset \mathcal S,
\end{equation}
which is of complex dimension $\dim{\mathcal S}' = 4d = 2\dim M$ by Proposition \ref{prop:SpaceOfSections}.
Since the vector bundle ${\mathcal O}_{\C P^1}(1)^{\oplus 2d}$ is semistable, from the openness
of the semistability condition, \cite[p.~635, Theorem 2.8\,(B)]{Maru}, it follows immediately
that ${\mathcal S}'$ in \eqref{t3} is an open subset of $\mathcal S$.

We want to transfer geometric objects from $Z$ to $\mathcal S$. Towards this end, it is useful to introduce
the correspondence space
\begin{equation}\label{ncf}
\mathcal F = \C P^1 \times \mathcal S
\end{equation}
as well as
\begin{equation}\label{n2}
\mathcal F' = \C P^1 \times \mathcal{S}'.
\end{equation}

\begin{Remark}
The correspondence space $\mathcal F$ is usually described as the space of pairs $(z, \ell)$, where $\ell$ is a complex
line in $Z$ passing through $z \in Z$.
In our work, the lines $\ell$ in $Z$ are treated as sections, meaning $\ell \in \mathcal{S}$.
 This yields an isomorphism
from this space of pairs $(z, \ell)$ to $\mathcal{F}$ in~\eqref{ncf} by mapping
any $(z, \ell)$ to $(\varpi(z), \ell)$.
\end{Remark}

Consider the evaluation map
$
\operatorname{ev}\colon\ \mathcal F \longrightarrow Z $, $ \operatorname{ev}(\lambda , s) = s(\lambda) $.
For a fixed $\lambda \in \C P^1 = \C\cup\{\infty\}$ define the map
\begin{equation}\label{t5}
\operatorname{ev}_\lambda\colon \ {\mathcal S} \longrightarrow Z_\lambda := \varpi^{-1}(\lambda)
 ,\qquad \operatorname{ev}_\lambda(s) = \operatorname{ev}(\lambda, s) .
\end{equation}
We will denote the restrictions of $\operatorname{ev}$ to $\mathcal F'$ and of $\operatorname{ev}_\lambda$ to $\mathcal S'$
(see \eqref{n2} and \eqref{t3} for ${\mathcal F}'$ and ${\mathcal S}'$ respectively) by $\operatorname{ev}$ and
$\operatorname{ev}_\lambda$ respectively.

\begin{Lemma}\label{Lem:evsubmersion}
The evaluation map $\operatorname{ev}\colon \mathcal F' \longrightarrow Z$ is a surjective holomorphic submersion.
\end{Lemma}

\begin{proof}
The surjectivity of $\operatorname{ev}\colon \mathcal F' \longrightarrow Z$ is clear, since there is a twistor line through each point
of $Z = \C P^1 \times M$.
If
\[(l, V) \in T_{(\lambda,s)}\mathcal F' = T_\lambda\C P^1 \oplus T_s\mathcal S'
 = T_\lambda\C P^1 \oplus H^0(s^*T_\varpi Z)\] (see \eqref{t3} for ${\mathcal S}'$), then we have
\[
({\rm d}\operatorname{ev})_{(\lambda,s)}(l, V) = ({\rm d}s)_\lambda(l)+V(\lambda).
\]
Now, since the evaluation map $H^0\bigl(\C P^1, {\mathcal O}_{\C P^1}(1)\bigr) \longrightarrow \mathcal{O}_{\C P^1}(1)_\lambda$ evaluated
at $\lambda$ is surjective,
the map ${\rm d}\operatorname{ev}_{(\lambda,s)}$ is surjective as well.
\end{proof}

We thus obtain the following commutative diagram of maps in which every arrow is a holomorphic submersion:
\begin{equation}\label{eq: doublefibration}
\begin{tikzcd}
Z \arrow[r, leftarrow, "\operatorname{ev}"] \arrow[d, "\varpi"'] & \mathcal F' =
\C P^1 \times \mathcal{S}' \arrow[dl, "\pi_1"'] \arrow[d, "\pi_2"] \\
\C P^1 & \mathcal{S}'.
\end{tikzcd}
\end{equation}

If $V \longrightarrow \mathcal F'$ is a holomorphic vector bundle such that $H^{q'}\bigl(\pi_2^{-1}(s), V\bigr) = 0$ for all
$s \in {\mathcal S}'$ (defined in \eqref{t3}) and $q' \not= q$, then the $q$-th direct image $R^q (\pi_2)_*V$
is a holomorphic vector bundle\footnote{We identify a holomorphic vector bundle with the locally free coherent analytic
sheaf given by its local sections.} on $\mathcal S'$ (see \cite[Chapter~10]{GrauertRemmert}) whose fiber over any $s \in \mathcal{S}'$
is $((\pi_2)_*V)_s = H^q\bigl(\pi_2^{-1}(s), V\bigr)$.
Hence the sheaves
\begin{equation}\label{eq: defEH}
\mathcal V := (\pi_2)_*\operatorname{ev}^*(T_\varpi Z(-1)), \qquad \mathcal H := (\pi_2)_*\pi_1^*\mathcal{O}_{\C P^1}(1)
\end{equation}
are holomorphic vector bundles over $\mathcal{S}'$ (defined in \eqref{t3}) because $s^*T_\varpi Z(-1) \cong \mathcal{O}_{\C P^1}^{2d}$
for any~${s \in \mathcal S'}$. Note that $(\pi_2)_*\pi_1^*W$ is in fact a trivial vector bundle with fiber $H^0\bigl(\C P^1, W\bigr)$.

\begin{Lemma}\label{lem:TSEH}
Consider ${\mathcal S}'$ in \eqref{t3}.
There is a canonical isomorphism $T\mathcal S' = \mathcal V\otimes\mathcal H$
$($defined in \eqref{eq: defEH}$)$, and the bundles
$\mathcal V$ and $\mathcal H$ carry natural holomorphic symplectic forms $\omega_{\mathcal V}$ and
$\omega_{\mathcal H}$. Thus, $\mathcal S'$ comes naturally equipped with the holomorphic Riemannian
metric $g = \omega_{\mathcal V}\otimes\omega_{\mathcal H}$.
\end{Lemma}

\begin{proof}
To prove the first statement, observe that
\[
T_{s}\mathcal S' = H^0\bigl(\C P^1, s^*T_\varpi Z\bigr) = H^0\bigl(\C P^1, s^*T_\varpi Z(-1)\bigr)\otimes H^0\bigl(\C P^1, \mathcal{O}_{\C P^1}(1)\bigr) = \mathcal V_s\otimes \mathcal H_s
\]
for every $s\in \mathcal{S}'$. This produces the identification $T\mathcal S' = \mathcal V\otimes\mathcal H$ in the lemma.

To obtain the symplectic forms, note that the twisted relative symplectic structure
$\omega \in H^0\bigl(Z, \bigl(\Lambda^2 T^*_\varpi Z\bigr)(2)\bigr)$ induces a natural symplectic form on the vector
bundle $T_\varpi Z(-1)$. Since~${\operatorname{ev}^*T_\varpi Z(-1)}$ is trivial on each fiber $\pi_2^{-1}(s)
 = \{s\}\times\C P^1$, $s \in \mathcal S'$, the pullback $\operatorname{ev}^*\omega$ induces a~symplectic form $\omega_{\mathcal V}\big\vert_s$ on
\[
\mathcal V_s = H^0\bigl(\pi_2^{-1}(s), \operatorname{ev}^*T_\varpi Z(-1)\bigr) = H^0\bigl(\C P^1, s^*T_\varpi Z(-1)\bigr).
\]
Finally, the symplectic form $\omega_{\mathcal{H}}\big\vert_s$ on $\mathcal H_s = H^0\bigl(\C P^1, \mathcal{O}_{\C P^1}(1)\bigr)$ is
induced from the Wronskian on $\pi_1^*\mathcal{O}_{\C P^1}(1)$. More
explicitly, take $\psi_1, \psi_2 \in \mathcal{H}_s$, and denote by ${\rm d}\psi_i$ the derivative
of $\psi_i$ defined in terms of local trivialisations. Then we set
\begin{equation}\label{eq:omegaH}
\omega_{\mathcal H}|_s(\psi_1,\psi_2) := \frac{1}{2}\left(
\psi_1\otimes({\rm d}\psi_2) - ({\rm d}\psi_1)\otimes\psi_2\right),
\end{equation}
which is a well-defined element of
\[H^0\bigl(\C P^1, K_{\C P^1}\otimes \mathcal{O}_{\C P^1}(1)^{\otimes 2}\bigr)
 = H^0\bigl(\C P^1, \mathcal{O}_{\C P^1}\bigr) = \C.\]
Since the bilinear forms $\omega_{\mathcal V}$ and $\omega_{\mathcal H}$ are non-degenerate on $\mathcal V$ and $\mathcal H$
respectively, it follows that~${g = \omega_{\mathcal V}\otimes\omega_{\mathcal H}}$ is non-degenerate on
$\mathcal V\otimes\mathcal H = T\mathcal S'$. Since both $\omega_{\mathcal V}$ and $\omega_{\mathcal H}$ are
anti-symmetric, we know that $g$ is symmetric. So $g$ is a holomorphic Riemannian metric on ${\mathcal S}'$.
This completes the proof.
\end{proof}

The restriction to $M \subset \mathcal S'$ of $g$ in Lemma \ref{lem:TSEH} coincides
with the Riemannian metric $g_0$ on the hyperK\"ahler manifold $M$.

The above considerations yield natural integrable distributions
\begin{equation}\label{b2}
T_{\operatorname{ev}}\mathcal F' := \ker({\rm d}\operatorname{ev})\qquad
\text{and}\qquad T_{\operatorname{ev}_x} \mathcal{S}' = \ker({\rm d}\operatorname{ev}_x),\qquad x \in \C P^1,
\end{equation}
on $\mathcal{F}'$ and $\mathcal{S}'$, respectively. The associated leaves of the foliations are the fibers
\[
\mathcal F_z := \operatorname{ev}^{-1}(z) \cong \mathcal{S}_z := \operatorname{ev}_x^{-1}(z), \qquad z \in Z,
\]
where $x = \varpi(z) \in \C P^1$.

\begin{Lemma}\label{lem:sl2web}
For any $x \in \C P^1$, the above integrable distribution $T_{\operatorname{ev}_{x}}\mathcal S'$ in
\eqref{b2} is maximally isotropic
with respect to $g$ in Lemma {\rm\ref{lem:TSEH}}. For two distinct points $x \neq y \in \C P^1$,
\[T_{\operatorname{ev}_{x}}\mathcal S'\cap T_{\operatorname{ev}_{y}}\mathcal S' = \{0\} .\]
\end{Lemma}

\begin{proof}
For any $s \in {\mathcal S}'$, we have
\[T_{\operatorname{ev}_{x}}\mathcal S'|_{s} = \bigl\{V \in
H^0\bigl(\C P^1, s^*T_\varpi Z\bigr) \mid V(x) = 0\bigr\},\]
and this subbundle is of rank $2d = \frac{1}{2}\dim \mathcal S'$. Now note that $\operatorname{ev}^*T_\varpi Z(-1)$ is trivial
on $\pi_2^{-1}(s) = \{s\}\times \C P^1$. We may choose an affine coordinate $\lambda$ on $\C P^1$ such that $\lambda(x) =
0$, and then view $\lambda$ as an element of $H^0\bigl(\C P^1, \mathcal{O}_{\C P^1}(1)\bigr)$. We may therefore write any $V \in
T_{\operatorname{ev}_{x}}\mathcal S'|_{s}$ in the form~${V = v\otimes \lambda}$. Thus, if we evaluate
$g(V_1, V_2)$ at $x \in \C P^1$, it follows from the formula for $\omega_\mathcal H$ (see~\eqref{eq:omegaH}) that we
actually get zero.

Recall that \[T_s{\mathcal S}' = H^0\bigl(\C P^1, s^*T_\varpi Z\bigr) = H^0\bigl(\C P^1, \mathcal{O}_{\C P^1}(1)^{\oplus 2d}\bigr).\]
Any holomorphic section of $\mathcal{O}_{\C P^1}(1)$ vanishing at two distinct points of $\C P^1$ must be identically
zero. This implies the second part of the lemma.
\end{proof}

Given any $x \in \C P^1$, we can thus define an associated non-degenerate holomorphic two-form~$\Omega_x$ on $\mathcal S'$
by taking the natural skew-form on
\begin{equation}\label{eq:splitlambda}
T\mathcal S' = T_{\operatorname{ev}_{x}}\mathcal S' \oplus T_{\operatorname{ev}_{\tau_{\C P^1}(x)}}\mathcal S'
\end{equation}
(recall that $\tau_{\C P^1}(x) \not= x$) induced from the holomorphic Riemannian metric
$g$: Write $V, W \in T\mathcal S' $ as $V = V_x + V_{\tau_{\C P^1}(x)}$, $W =
W_x + W_{\tau_{\C P^1}(x)}$ with respect to the splitting in \eqref{eq:splitlambda}, and put
\begin{equation}\label{eq:defOmegalambda}
\Omega_x(V, W) := -\i(g(V_x,W_{\tau_{\C P^1}(x)}) - g(V_{\tau_{\C P^1}(x)},W_x)).
\end{equation}
We may now define an endomorphism $I_x \in H^0(\mathcal{S'}, \mathrm{End}(T\mathcal S'))$ via
\begin{equation}\label{eq:Ilambda}
I_x = \begin{pmatrix}
-\i & 0 \\ 0& \i
\end{pmatrix}
\end{equation}
again with respect to the splitting in \eqref{eq:splitlambda}.
Then we see immediately that $I_x$ is in fact orthogonal with respect to $g$, and it satisfies the equation
$
\Omega_x = g(I_x -,-)$.

Next consider the map
\begin{equation}\label{pre-phi}
\phi_{x} := \operatorname{ev}_x\times\operatorname{ev}_{\tau_{\C P^1}(x)} \colon \
\mathcal S \longrightarrow Z_x \times Z_{\tau_{\C P^1}(x)} = Z_x\times \overline{Z}_x.
\end{equation}

\begin{Proposition}\label{p:local biholo}
The map $\phi_x$ in \eqref{pre-phi} restricts to a local biholomorphism
\[\phi_x \colon\ \mathcal{S}'
 \longrightarrow Z_x\times \overline{Z}_x.\]
\end{Proposition}

\begin{proof}
Clearly, the spaces have the common complex dimension $4d$. The kernel of the differential
${\rm d}\phi_x$ is given by \smash{$T_{\operatorname{ev}_{x}}\mathcal S'\cap
T_{\operatorname{ev}_{\tau_{\C P^1}(x)}}\mathcal S'$}, and we have
\[T_{\operatorname{ev}_{x}}\mathcal S'\cap T_{\operatorname{ev}_{\tau_{\C P^1}(x)}}\mathcal S' = \{0\}\]
 by
Lemma \ref{lem:sl2web}. The proposition follows.
\end{proof}

Now we discuss how the above data interact with the real structure $\tau$ on $\mathcal S'$ defined in \eqref{tau}.
Let
\begin{equation}\label{mp}
M' := (\mathcal{S}')^\tau
\end{equation}
be the space of real sections so that we have an embedding $M \hookrightarrow M'$ induced by \eqref{te}.
Hence~$\mathcal{S}'$ is a natural complexification of the real analytic smooth manifold $M'$. Note that
in some examples $M$ is all of $M'$ \big(e.g., for the standard flat hyperK\"ahler
manifolds $\C^{2d}$\big) but not always (see Example \ref{ex:grafting}).

For any $s \in M'$ in \eqref{mp}, the differential ${\rm d}\tau \colon T_s{\mathcal S}'
 \longrightarrow T_s{\mathcal S}'$ is $\mathbb C$-antilinear, involutive
and satisfies the equation
\[
{\rm d}\tau(T_{\operatorname{ev}_{x}}\mathcal S') = \overline{T_{\operatorname{ev}_{\tau_{\C P^1}(x)}}\mathcal S'}.
\]
Indeed, we have for $V \in T_s\mathcal S' = H^0(s^*T_\varpi Z)$ the formula
\[
{\rm d}\tau(V)(x) = {\rm d}\tau_Z(V(\tau_{\C P^1}(x))).
\]
Thus, if $V(x) = 0$, then ${\rm d}\tau(V)(\tau_{\C P^1}(x)) = 0$. This implies in
particular that if $s \in M'$, and~${V \in (T_s\mathcal S')^\tau}$ is \emph{real}, then the
section $V \in H^0\bigl(\C P^1, s^*T_\varpi Z\bigr)$ is either identically zero or it is nowhere
vanishing (a nonzero holomorphic section of ${\mathcal O}_{\C P^1}(1)$ cannot vanish at two
distinct points of $\C P^1$). As a consequence, the map $\operatorname{ev}_x$ in \eqref{t5} gives a
local diffeomorphism from~$M'$ to $Z_x = (M, I_x)$. Moreover, $I_x$ is real in the sense that
${\rm d}\tau \circ I_x = I_x\circ {\rm d}\tau$, and therefore it preserves~${TM' \cong (T\mathcal S')^\tau}$.
Hence it defines an almost complex structure on $M'$. Consequently, we obtain a hypercomplex
structure $\bigl\{(\operatorname{ev}_x)^* I_x \mid x \in \C P^1\bigr\}$ on $M'$.

\begin{Remark}\label{rem:T10M}
The differential of the inclusion map $M \hookrightarrow M'$ is a $\mathbb R$-linear isomorphism
$T_{m}M \stackrel{\sim}{\longrightarrow} T_{s_m}M'$ whose complexification yields an identification
$T_s\mathcal S' \cong T_{s(0)}M\otimes \C$. Under this identification, the
decomposition in \eqref{eq:splitlambda} is actually mapped to the natural decomposition~${T_{m}M\otimes\C \cong T_m^{0,1}M\oplus T_m^{1,0}M}$ with respect to the complex structure $I_x$.
\end{Remark}

The real tangent vectors at $s \in M'$ can be described as follows: Let $V \in
T_{\operatorname{ev}_x}\mathcal S'\vert_s$. Then the tangent vector
\[
V + \tau(V) \in T_{\operatorname{ev}_x}\mathcal S'\vert_s \oplus T_{\operatorname{ev}_{\tau_{\C P^1}(x)}}\mathcal S'\vert_s
\]
is obviously real and we have
\[
g(V+\tau(V), V+\tau(V)) = 2g(V, \tau(V)).
\]
Since the twisted relative symplectic form $\omega$ on $Z$ satisfies
the condition $\tau_Z^*\omega = \overline\omega$, it follows,
by working through the definition of $g$, that $g$ is real in the sense that $\tau^*g = \overline g$. Hence
$g$ induces a~real-valued pseudo-Riemannian metric on $M'$. Note that this immediately forces the restriction of
$\Omega_x$ to $M'$ to be real as well. Pulled back to $(M, I_x)$, the form $\Omega$ is just the K\"ahler form
associated with the pseudo-Riemannian metric and the hermitian almost complex structure~$I_x$.

In summary, we have obtained the following result.

\begin{Proposition}
For every $x \in \C P^1$, the form $\Omega_x$ constructed in \eqref{eq:Ilambda}
defines a holomorphic symplectic form on each component of $\mathcal S'$ that intersects $M'
 = (\mathcal S')^\tau$. On $M \subset M'$ it induces the K\"ahler form $\omega_x$.
\end{Proposition}

\begin{Remark}
So far we have not shown that $\Omega_x$ is actually closed.
One way to show that is to use the Atiyah--Ward transform which will be discussed in detail in Section \ref{ss:AW}.
The bundle~$\mathcal V$ can be seen to arise from this transform applied to the bundle $T_\varpi Z(-1)$.
As such it carries a~natural connection. Its tensor product with the trivial connection on the trivial
bundle $\mathcal H$ can be shown to give the Levi-Civita connection of the holomorphic Riemannian manifold
$(T\mathcal S', g)$. The form $\Omega_x$ can then be shown to be parallel with respect to
the Levi-Civita connection, from which it follows that the form $\Omega_x$ is closed.
In the next subsection we will give another, more direct proof that $\Omega_x$ is closed.
\end{Remark}

\subsection[Alternative description of the holomorphic symplectic form Omega\_x on S']{Alternative description of the holomorphic symplectic form $\boldsymbol{\Omega_x}$ on $\boldsymbol{\mathcal{S}'}$}\label{ss:Omega}

Fix a point $x \in \C P^1$. We now give an alternative description of $\Omega_x$ which is
better suited for performing computations. This alternative description also shows that $\Omega_x$ can be
extended to a~holomorphic two-form on $\mathcal S$, which is typically strictly larger than
$\mathcal S'$.

Consider the diagram
\[
\begin{tikzcd}
Z \arrow[r, leftarrow, "\operatorname{ev}"] \arrow[d, "\varpi"'] & \mathcal F = \C P^1 \times \mathcal{S} \arrow[dl, "\pi_1"'] \arrow[d, "\pi_2"] \\
\C P^1 & \mathcal{S}.
\end{tikzcd}
\]
 We observe that for any $k \in \mathbb{Z}$, the direct image $(\pi_2)_*\pi_1^*\mathcal{O}_{\C P^1}(k)$
on $\mathcal S$ is a trivial vector bundle with fiber $H^0\bigl(\C P^1, \mathcal{O}_{\C P^1}(k)\bigr)$.

Starting with the twisted relative symplectic form $\omega \in H^0\bigl(Z, \bigl(\Lambda^2T^*_\varpi Z\bigr)(2)\bigr)$, its pullback~$\operatorname{ev}^*\omega$ defines a holomorphic section of \smash{$\bigwedge^2\pi_2^*(T^*\mathcal S)(2)$}.
Invoking pushforward to $\mathcal S$, we obtain a~vector-valued holomorphic two-form
\begin{equation}\label{eq:Omegavector}
\Omega \in H^0\bigl(\mathcal S, \Lambda^2 T^*\mathcal S\bigr)\otimes H^0\bigl(\C P^1, \mathcal{O}_{\C P^1}(2)\bigr).
\end{equation}
Now note that $H^0\bigl(\C P^1, \mathcal{O}_{\C P^1}(2)\bigr)$ is the space of all holomorphic vector fields on $\C P^1$ and therefore
it has the structure of a Lie algebra isomorphic to $\sln_2(\C)$; the Lie algebra structure is given by the Lie bracket
operation of vector fields.
Fix an affine coordinate $\lambda$ on $\C P^1$ such that~${\lambda(x) = 0}$ and $\tau(\lambda)
 = -(\overline{\lambda})^{-1}$. Then we obtain the following basis of $H^0\bigl(\C P^1, \mathcal{O}_{\C P^1}(2)\bigr)$
\[
e = \frac{\partial}{\partial\lambda},\qquad h = -2\lambda\frac{\partial}{\partial\lambda}, \qquad f
 = -\lambda^2\frac{\partial}{\partial\lambda},
\]
which satisfies the standard relations
$
[h, e] = 2e$, $ [h, f] = -2f$, $ [e, f] = h$.
We can now use the Killing form $\kappa$ on $\sln_2(\C) = H^0\bigl(\C P^1, \mathcal{O}_{\C P^1}(2)\bigr)$ to define
\begin{equation}\label{def:Omegax}
\widetilde{\Omega}_x := \frac{1}{8 \i}\kappa(\Omega, h) \in
H^0\left(\mathcal S, \bigwedge\nolimits^2 T^*\mathcal S\right).
\end{equation}
Note that $\widetilde{\Omega}_x$ in \eqref{def:Omegax} is independent of the affine coordinate $\lambda$ at $x$.
Using $\kappa(h, h) = \operatorname{tr}(\mathrm{ad}_h\circ\mathrm{ad}_h) = 8$ and $\kappa(h, e)
 = 0 = \kappa(h, f)$ we may rewrite this as follows. A general element $A
 \in H^0\bigl(\C P^1, \mathcal{O}_{\C P^1}(2)\bigr)$ is of the form
\[A = A_ee + A_hh + A_ff = \bigl(A_e -2\lambda A_h-\lambda^2A_f\bigr)\frac{\partial}{\partial\lambda}
=: A(\lambda)\frac{\partial}{\partial\lambda}\]
with $A_e, A_h, A_f \in \C$. Then
\[
\frac{1}{8 \i}\kappa(A, h) = -\i A_h = \frac{\i}{2}\frac{\partial}{\partial\lambda}_{|\lambda =0}A(\lambda).
\]
With this setup in place, we may therefore write for $s \in \mathcal S'$, and tangent vectors
\[V_s, W_s \in T_s\mathcal S' = H^0\bigl(\C P^1, s^*T_{\varpi}Z\bigr),\]
the following
\begin{align}\label{eq:OmegaGeneral}
\widetilde\Omega_x\vert_s(V, W) & = \frac{1}{8\i}\kappa(\operatorname{ev}^*\omega(\pi_2^*V, \pi_2^*W), h) = \i \pder{t}{\lambda}_{|\lambda = 0}\omega_{s(\lambda)} \left(V_s(\lambda), W_s(\lambda)\right).
\end{align}
Recall the non-degenerate two-form $\Omega_x$ defined in \eqref{eq:defOmegalambda}.

\begin{Theorem}\label{thm:Omega}
Let $x \in \C P^1$.
\begin{enumerate}[$(a)$]\itemsep=0pt
\item The two-form $\widetilde \Omega_x \in H^0\bigl(\mathcal S, \bigwedge^2T^*\mathcal S\bigr)$ defined in \eqref{eq:OmegaGeneral} restricts to a holomorphic symplectic form on $\mathcal{S}'$ which is real with respect to $\tau$.

\item Over the open subset $\mathcal S' \subset \mathcal S$,
$\Omega_x\vert_{\mathcal S'} = \widetilde\Omega_x\vert_{\mathcal S'} $.
In particular, $(\mathcal{S}', \Omega_0)$ is a complexification of the $($real analytic$)$ K\"ahler manifold
$(M', \omega_I)$, where $\omega_I$ is the K\"ahler form associated to~${0 \in \C P^1}$.

\item The distributions $T_{\operatorname{ev}_x}\mathcal S'$ and \smash{$T_{\operatorname{ev}_{\tau_{\C P^1}(x)}}\mathcal S'$} are
Lagrangian with respect to $\widetilde\Omega_x$.
\end{enumerate}
\end{Theorem}

\begin{proof}
We choose and fix throughout the proof an affine coordinate $\lambda$ on $\C P^1$ such that $\lambda(x) = 0$ and \smash{$\tau(\lambda)
 = -(\overline{\lambda})^{-1}$}.
To prove the first statement, recall that $\tau_Z^*\omega = \overline{\omega}$.
This implies that $\widetilde{\Omega}_x$ is real with respect to $\tau$.

Next we show that $\widetilde{\Omega}_x$ is a holomorphic symplectic form on $\mathcal{S}'$. Recall the diagram \eqref{eq: doublefibration}. We know from Lemma \ref{Lem:evsubmersion} that $\operatorname{ev}\colon \mathcal F' \longrightarrow Z$ is a holomorphic submersion, which for any fixed value of $\lambda$ maps the fibre $\pi_1^{-1}(\lambda)\subset \mathcal F'$ to the fibre $\varpi^{-1}(\lambda)\subset Z$. Moreover, we know that $\omega$ restricts to a symplectic, hence closed, form on each fibre $\varpi^{-1}(\lambda)$. From these facts, we deduce that the vector-valued form $\Omega$ defined in \eqref{eq:Omegavector} is closed and hence also the restriction of $\widetilde{\Omega}_x$ to the open subset $\mathcal{S}'$.

We will now show part (b), i.e., that $\widetilde{\Omega}_x$ agrees pointwise on $\mathcal S'$ with the form $\Omega_x = g(I_x - , -)$ defined in \eqref{eq:defOmegalambda}. Fix $s\in\mathcal S'$ and take two tangent vectors $V,W\in T_s\mathcal S' = H^0\bigl(\C P^1, N_s\bigr) = \mathcal V_s\otimes \mathcal H_s$, where we use Lemma \ref{lem:TSEH}. We write
$V(\lambda)= a_1 + a_2\lambda$, $ W(\lambda) = b_1+b_2\lambda$
and compute, using the notation of Lemma \ref{lem:TSEH}
\begin{align*}
\widetilde{\Omega}_x(V,W) &= \frac{\i}{2}\frac{\partial}{\partial\lambda}|_{\lambda=0}\omega(V(\lambda),W(\lambda))= \frac{\i}{2}\frac{\partial}{\partial\lambda}|_{\lambda=0} \omega(a_1+a_2\lambda,b_1+b_2\lambda)\\
 &= \frac{\i}{2}\frac{\partial}{\partial\lambda}|_{\lambda=0}\bigl(\omega_{\mathcal V}(a_1,b_1) + (\omega_{\mathcal V}(a_1,b_2)+\omega_{\mathcal V}(a_2,b_1))\lambda + \omega_{\mathcal V}(a_2,b_2)\lambda^2\bigr)\\
 &= \frac{\i}{2}(\omega_{\mathcal V}(a_1,b_2)+\omega_{\mathcal V}(a_2,b_1)).
\end{align*}
On the other hand, viewing $1$ and $\lambda$ as sections of $\mathcal O_{\C P^1}(1)$ which vanish at $\infty\in\C P^1$ and $0\in \C P^1$, respectively, we have by the definition of $\Omega_x$ given in \eqref{eq:defOmegalambda}
\begin{align*}
\Omega_x(V,W) &= -\i(g(a_2\lambda,b_1) - g(a_1,b_2\lambda))= -\i(\omega_{\mathcal V}(a_2,b_1) + \omega_{\mathcal V}(a_1,b_2))\omega_{\mathcal H}(\lambda,1)\\
&= \frac{\i}{2}(\omega_{\mathcal V}(a_2,b_1) + \omega_{\mathcal V}(a_1,b_2))= \widetilde{\Omega}_x(V,W).
\end{align*}

Since $\Omega_x =
\widetilde{\Omega}_x$ the remaining claims follow because they have been established for $\Omega_x$.
\end{proof}

\begin{Remark}
Both definitions of $\Omega_x$ will be used. For example, \eqref{eq:OmegaGeneral} makes sense on all of $\mathcal S$ and thus
can be evaluated on any holomorphic section $s$, even if we do not have the knowledge of the normal bundle of $s$.
\end{Remark}

\section{The hyperholomorphic line bundle and the energy functional}
\label{sectionhyperholo}
\subsection{Rotating circle actions and the hyperholomorphic line bundle}\label{ss:rotating}

Now assume that $M$ is equipped with an action of $S^1$ that preserves the Riemannian metric $g_0$ and also preserves
the family of complex structures $\{I_x\}_{x\in S^2}$. We also assume that the resulting action of $S^1$
on $S^2$ is nontrivial. This implies that the $S^1$-action on $M$ preserves the associative algebra structure of
$({\mathbb R}\cdot {Id}_{TM})\oplus ({\mathbb R}\cdot I) \oplus ({\mathbb R}\cdot J)\oplus ({\mathbb R}\cdot K)$.
Consequently, the action of~$S^1$ on~$S^2$ is a nontrivial rotation.
This means that without loss of generality we may
assume that the $S^1$-action on $M$ preserves $I$ (thus it also preserves $-I$) and rotates the plane spanned by~$J$,~$K$ in the standard way. We therefore call this a \emph{rotating circle action}. The K\"ahler forms~$\omega_I$,~$\omega_J$,~$\omega_K$ and the Killing vector field $X$ on $M$ associated with this circle
satisfy the following:
\begin{equation}\label{eq:RotS1LieDer}
\mathcal L_X\omega_I = 0, \qquad \mathcal L_X\omega_J = \omega_K, \qquad \mathcal L_X\omega_K = -\omega_J.
\end{equation}
Note that the second equation in \eqref{eq:RotS1LieDer} implies that the K\"ahler form $\omega_K$ is exact,
and the third equation in \eqref{eq:RotS1LieDer} implies that the K\"ahler form $\omega_J$ is exact.
Therefore, the manifold $M$ must necessarily be non-compact.

The above $S^1$-action on $M$ evidently induces a holomorphic $S^1$-action on $Z$.
In terms of the identification $Z = \C P^1 \times M$ of the underlying smooth manifold, this
action of $S^1$ is given by
\begin{equation}\label{S1actionTrivialization}
\zeta.(\lambda, m) = (\zeta\lambda, \zeta.m), \qquad \ \zeta \in S^1 .
\end{equation}
Here $\lambda$ is an affine coordinate on $\C P^1$ such that $I = I_0$. The $C^\infty$-vector field
on $Z$ associated to the $S^1$-action in \eqref{S1actionTrivialization} will be denoted by $Y$. Note that
$Y\vert_{Z_0\cup Z_\infty}$ (see \eqref{zl}) is actually a~holomorphic vector field on the divisor
$Z_0 \cup Z_\infty = (M, I)\cup (M, -I)$,
because the $S^1$-action on~$M$ preserves both $I$ and $-I$.

We normalize the affine coordinate $\lambda$ on $\C P^1$ such that the antipodal map $S^2 \longrightarrow S^2$,
$ x \longmapsto -x$, corresponds to the map $\C P^1 \longrightarrow \C P^1$, $ \lambda \longmapsto
-(\overline{\lambda})^{-1}$. The coordinate $\lambda$ is then uniquely determined up to multiplication by a constant
phase ${\rm e}^{\i \theta_0}$, $\theta_0 \in [0, 2\pi)$.

Clearly, the $S^1$-action on $Z$ is compatible with the real structure $\tau_Z$ in \eqref{etz} in the sense that
\begin{gather}\label{eq:S1tau}
\tau_Z(\zeta.z) = (\overline{\zeta})^{-1}.\tau_Z(z) = \zeta .\tau_Z(z)
\end{gather}
for all $z \in Z$ and $\zeta \in S^1$. It is straightforward to check that for any $\zeta \in S^1$ the
twisted relative
symplectic form $\omega$ in \eqref{t6} satisfies the equation
$
\zeta^*\omega = \omega$,
in other words, $\omega$ is $S^1$-invariant.

Let $\mu \colon M \longrightarrow \i \R = \mathfrak{u}(1)$ be a moment map with respect to $\omega_I$
for the $S^1$-action on $M$. Note that since our moment map is complex-valued, the moment map equation takes the form
\begin{equation}\label{eq:MomentMapEq}
{\rm d}\mu(-) = \i \omega_I(X, -).
\end{equation}
Haydys has shown in \cite{Haydys} that the $2$-form $\omega_I + \i\, {\rm dd}_I^c \mu$ is of type $(1, 1)$ with respect to
every complex structure $I_\lambda$, $\lambda \in {\mathbb C}P^1$.

For the remainder of Section~\ref{sectionhyperholo}, we shall make the hypothesis that $\omega_I$ (and hence $\omega_I + \i \, {\rm dd}_I^c \mu$) is \emph{integral}, i.e.,
\begin{equation}\label{hyp:integral}
[\omega_I/2\pi] \in H^2(M, \mathbb{Z}).
\end{equation}
Then there exists a $C^\infty$ hermitian line bundle
$(L_M, h_M)
 \longrightarrow M$
with a compatible hermitian connection $\nabla_M$ whose curvature is $\omega_I + \i\, {\rm dd}_1^c \mu$.
Consequently, if~\smash{$\nabla_M^{(0,1)_\lambda}$} is the $(0,1)$-part of $\nabla_M$ with respect to the complex structure
$I_\lambda$, we obtain a holomorphic line bundle~${\bigl(L_M, \nabla^{(0,1)_\lambda}\bigr) \longrightarrow (M,
I_\lambda)}$. Note that the Chern connection of the hermitian holomorphic line bundle $(L, {\rm e}^{-\mu}h_M)
 \longrightarrow (M, I)$ has curvature $\omega_I$, and hence it is a prequantum line bundle on the K\"ahler
manifold $(M, \omega_I)$.

Denote by
$q \colon Z \longrightarrow M$
the $C^\infty$ submersion
given by the natural projection
${S^2 \times M \longrightarrow M}$.
Note that for each $m\in M$ the fibre of $q^{-1}(m)$ is exactly the image of the twistor line ${s_m\colon \C P^1 \longrightarrow Z}$, where we use the notation of Section~\ref{ss:twistor}.
It follows that
\[L_Z = \bigl(q^*L_M, (q^*\nabla_M)^{0,1}\bigr) \longrightarrow Z\]
is a holomorphic line bundle
over the twistor space $Z$ of $M$, which is trivial along each twistor line.

In \cite{Hitchin-HKQK}, Hitchin provided a twistorial description of the line bundle $L_Z$ exhibiting a natural meromorphic connection $\nabla$ on $L_Z$.
To recall this, observe that the $S^1$-action on $Z$ covers the standard $S^1$-action on $\C P^1$, and therefore the associated holomorphic vector field $Y$ on $Z$ is $\varpi$-related to
$
\sigma := \i\lambda\frac{\partial}{\partial\lambda}
$
on $\C P^1$. Viewing $\sigma$ as a section of $\varpi^*{\mathcal O}_{\C P^1}(2)$ which vanishes on the divisor
$D := Z_0\cup Z_\infty$, we have the short exact sequence
\begin{equation}\label{eq:ses}
\begin{tikzcd}
0 \ar[r] & T^*Z \ar[r, "\cdot \sigma"] & T^*Z(2) \ar[r] & T^*Z(2)_{|D} \ar[r] & 0
\end{tikzcd}
\end{equation}
of coherent analytic sheaves on $Z$.
Hitchin then constructs from the $S^1$-action a certain element~$\varphi \in H^0(D, T^*Z(2)|_D)$.
Explicitly, using the $C^\infty$-splitting $T^*Z = T^*M \oplus T^*\C P^1$ we can write $\varphi$ in terms
of the data on $M$ as
\[
\varphi = \left(\frac{1}{2}\big({\rm d}_2^c\mu+{\rm i}{\rm d}_3^c\mu+2{\rm i}\lambda {\rm d}_1^c\mu+\lambda^2\big({\rm d}_2^c\mu- {\rm i} {\rm d}_3^c\mu\big)\big)\otimes \frac{\partial}{\partial \lambda},
 (\mu {\rm d}\lambda) \otimes \frac{\partial}{\partial \lambda} \right) .
\]
It is shown in \cite[Lemma 3.7]{BeHeRo} that $\varphi$ satisfies the equation
\begin{equation}\label{eq:varphiomega}
\varphi\big\vert_{T_\varpi Z|_D} = -\frac{1}{2}\iota_Y\omega,
\end{equation}
where $\omega$ is the relative symplectic form on $Z$. Note that since the vector field $Y|_D$ is vertical, the formula
in \eqref{eq:varphiomega} makes sense. From the long exact sequence of cohomologies associated to~\eqref{eq:ses}
\[
\begin{tikzcd}
0 \ar[r] & H^0(Z,T^*Z) \ar[r, "\cdot \sigma"] \ar[draw=none]{d}[name=X, anchor=center]{} & H^0(Z,T^*Z(2))
\ar[rounded corners, to path={ -- ([xshift=2ex]\tikztostart.east)
 |- (X.center) \tikztonodes
 -| ([xshift=-2ex]\tikztotarget.west)
 -- (\tikztotarget)}]{dl}[at end]{} \\
& H^0(D,T^*Z(2)_{|D}) \ar[r, "\delta"] & H^1(Z,T^*Z)\ar[r]&\cdots,
\end{tikzcd}
\]
we then obtain that
$\alpha_L := \delta(\varphi) \in H^1(Z, T^*Z)$.
In fact, this element $\alpha_L $ lies in the image of~$H^1\bigl(Z, \Omega_{Z,{\rm cl}}^1\bigr)$, where $\Omega^1_{Z,{\rm cl}}$ denotes the sheaf
of closed $1$-forms on $Z$ (it is not a coherent analytic sheaf). The above class $\alpha_L$ therefore defines an extension
\[
\begin{tikzcd}
0\arrow[r] & \mathcal O \arrow[r] & E \arrow[r] & TZ \arrow[r] & 0
\end{tikzcd}
\]
such that $E \longrightarrow TZ$ is a holomorphic Lie algebroid. The Lie algebroid $E \longrightarrow TZ$ exists even if \eqref{hyp:integral} does not hold.
Under the assumption \eqref{hyp:integral}, the form $\alpha_L$ is integral, so $E \longrightarrow TZ$ is the Atiyah algebroid of a line bundle $L_Z$ and $\alpha_L$ is its Atiyah class.
Explicitly, relative to some open cover $\mathcal U = \{U_i\}$ of $Z$, we have
$
(\alpha_L)_{ij} = g_{ij}^{-1}dg_{ij}
$
and $\bigl\{g_{ij}\in H^0(U_{i}\cap U_j, \mathcal O_Z^*)\bigr\}$ is a~cocycle representing $L_Z$.
Since on the other hand $\alpha_L = \delta(\varphi)$, we may rewrite this using the definition of the connecting homomorphism $\delta$.
Let $\varphi$ be given by $\bigl\{\varphi_i \in H^0(U_i\cap D, T^*Z(2)|_D)\bigr\}$ and take extensions $\bigl\{\widetilde\varphi_i
 \in H^0(U_i, T^*Z(2))\bigr\}$.
Then a representative for $\alpha_L = \delta(\varphi)$ is given by~\smash{$
(\alpha_L)_{ij} = \frac{\widetilde\varphi_i-\widetilde\varphi_j}{\sigma}$}.
Thus, \smash{$A_i = \frac{\widetilde\varphi_i}{\sigma}$} are the local $1$-forms of a meromorphic connection $\nabla$
on $L_Z$, which have a simple pole along the divisor $D$. Its residue along $D$ is given by $\varphi$.

Moreover, Hitchin shows that the curvature $F$ of the meromorphic connection $\nabla$ has the following properties:
\begin{itemize}\itemsep=0pt
\item $\iota_Y F = 0$ and $Y$ spans the annihilator of $ F$.

\item $ F = \frac{\omega}{\sigma}$ on $T_\varpi^*Z|_{\varpi^{-1}(\C^*)}$, where $\omega
 \in H^0\bigl(Z, \bigl(\bigwedge^2T^*_\varpi Z\bigr)(2)\bigr)$ is the relative symplectic form on~${Z \longrightarrow \C P^1}$.
\end{itemize}
We may thus think of $(L_Z, \nabla)$ as a ``meromorphic relative prequantum data'' for the meromorphic
relative symplectic form $\frac{\omega}{\sigma}$.

\subsection[The S\^{}1-action on S and the energy functional]{The $\boldsymbol{S^1}$-action on $\boldsymbol{\mathcal S}$ and the energy functional}

The $S^1$-action on $Z$ obtained from the $S^1$-action on $M$ goes on to
produce an $S^1$-action on $\mathcal S$, which is constructed as follows
$
(\zeta.s)(\lambda) = \zeta.\bigl(s\bigl(\zeta^{-1}\lambda\bigr)\bigr)
$
for all $s \in \mathcal S$, $\zeta \in S^1$ and $\lambda \in \C P^1$.
The evaluation map $\operatorname{ev}\colon \C P^1 \times \mathcal S \longrightarrow Z$ is evidently
equivariant with respect to the diagonal action of $S^1$ on $\C P^1 \times \mathcal S$ and the $S^1$-action on $Z$.
Moreover, the $S^1$-action on $\mathcal S$ is compatible with $\tau$ in the sense
that
$
\tau(\zeta.s) = \zeta.\tau(s)
$
for all $s \in \mathcal S$ and $\zeta \in S^1$.
In particular, the $S^1$-action on $\mathcal S$ preserves the subset ${\mathcal S}^\tau \subset
{\mathcal S}$ fixed pointwise by $\tau$.

\begin{Proposition}
The holomorphic two-form $\Omega_0$ on ${\mathcal S}$ constructed in \eqref{eq:defOmegalambda} is $S^1$-invariant.
\end{Proposition}

\begin{proof}
We use the description of $\Omega_0$ given in Section~\ref{ss:Omega}. Recall that $\omega$ is $S^1$-invariant and
$\operatorname{ev}\colon \mathcal F \longrightarrow Z$ is $S^1$-equivariant. Moreover, the Killing form on
$\sln_2 = H^0\bigl(\C P^1, \mathcal{O}_{\C P^1}(2)\bigr)$ and the element \[h = 2\i\sigma \in H^0\bigl(\C P^1, \mathcal{O}_{\C P^1}(2)\bigr)\]
are also $S^1$-invariant. Hence $\Omega_0 = \frac{1}{8 \i}\kappa( \Omega,h)$ too must be
$S^1$-invariant.
\end{proof}

Let $Y$ be the vector field on $Z$ which is induced by the
$S^1$-action. Since the $S^1$-action commutes with $\tau_Z$ (see \eqref{eq:S1tau}), we conclude that
$Y$ is $\tau_Z$-invariant, meaning ${\rm d}\tau_Z(Y) = Y\circ \tau_Z$.

\begin{Lemma}\label{RelationXandY}
Let $s \in \mathcal S'$, and let $X$ be the vector field associated to the $S^1$-action on ${\mathcal S}'$.
Then $X_s \in T_s\mathcal S' \cong H^0\bigl(\C P^1, N_s\bigr)$ is given by
\begin{equation}\label{Eq: XsYsdots}
X_s(\lambda) = Y_{s(\lambda)} - \i \lambda \dot{s}(\lambda).
\end{equation}
\end{Lemma}

Note that $X_s$ is indeed vertical because
\[{\rm d}\varpi(Y_{s(\lambda)} -\i \lambda\dot{s}(\lambda)) = \i \lambda\frac{\partial}{\partial\lambda} -\i
\lambda\frac{\partial}{\partial\lambda} = 0 .
\]

\begin{proof}
This is a direct computation: take $s \in \mathcal S'$ and $\lambda \in \C P^1$. Then
\[
X_s(\lambda) = \frac{\rm d}{{\rm d}t}\Big\vert_{t=0} \bigl({\rm e}^{\i t}.s\bigr)(\lambda)
 = \frac{\rm d}{{\rm d}t}\Big\vert_{t=0} {\rm e}^{\i t}.\bigl(s\bigl({\rm e}^{-\i t}\lambda\bigr)\bigr)
 = Y_{s(\lambda)} -\i \lambda\dot{s}(\lambda).
\]
More invariantly, we can use the vector field $\sigma$ for the standard $S^1$-action on $\C P^1$, that is,
$\sigma(\lambda)(\lambda) = \i\lambda\frac{\partial}{\partial \lambda}$ to rewrite the relation of Lemma \ref{RelationXandY} as
\begin{equation}\label{Eq: XsYsev}
{\rm d}\operatorname{ev}\circ X = Y\circ\operatorname{ev} - {\rm d}\operatorname{ev}\circ \sigma.
\end{equation}
The fundamental vector field $Y_{\mathcal F}$ for the diagonal $S^1$-action on $\mathcal F
 = \C P^1 \times \mathcal S$ is given by
\[
Y_{\mathcal F}(\lambda,s) = \pi_1^*\sigma(\lambda) + \pi_2^*X(s).
\]
Thus, the formulae \eqref{Eq: XsYsdots} and \eqref{Eq: XsYsev} just become
${\rm d}\operatorname{ev}\circ Y_{\mathcal{F}} = Y\circ \operatorname{ev}$.
\end{proof}

The recent paper \cite{BeHeRo} discusses the holomorphic
energy functional
\begin{equation}\label{rf}
\mathcal E\colon\ \mathcal S \longrightarrow \C,\ \qquad s \longmapsto
 \mathrm{res}_{\lambda=0}(s^*\nabla).
\end{equation}
Note that $\mathcal E(s) = s^*\varphi \in T^*\C P^1(2)\big\vert_{\{0,\infty\}}
 \cong \mathcal{O}_{\C P^1}|_{\{0,\infty\}}$, so ${\mathcal E}(s)$ is indeed a complex number.

The functional $\mathcal E$ in \eqref{rf} has the property that $\iota^*\mathcal E = \mu$ for the natural inclusion
map~${\iota \colon M
 \longrightarrow \mathcal{S}}$ in \eqref{te}. The residue formula in \eqref{rf} implies that
\begin{equation}\label{EnergyResidue}
\overline{\tau^*\mathcal{E}(s)} = -\mathcal{E}(s) + \mathrm{deg}(s^*L_Z).
\end{equation}

To give an explicit formula for the energy $\mathcal{E}$, we recall the section
\[\varphi_0 = \varphi_{|Z_0} \in
\Gamma(Z_0, T_Z^*(2)_{|Z_0})\]
 defined above. Contracting ${\rm d}\lambda \otimes \frac{\partial}{\partial
\lambda}$, we therefore obtain that
\[
\mathcal{E}(s) = s^* \varphi_0 = -\frac{1}{2} \iota_Y\omega \left(\dot{s}(0)- \dot{s}_{s(0)}(0)\right) + \mu(s(0)),
\]
where $\dot{s}(0) = {\rm d}s_0\bigl(\frac{\partial}{\partial \lambda}\bigr)$ etc.
The difference $\dot{s}(0) - \dot{s}_{s(0)}(0)$ accounts for the fact that we are working with
the $C^\infty$-splitting induced by twistor lines; see \cite{BeHeRo} for the details.

The following theorem is a holomorphic extension of \eqref{eq:MomentMapEq} because we have $\Omega_0|_M = \omega_I$.

\begin{Theorem}\label{EnergyMomentMap}
The energy functional $\mathcal{E} \colon \mathcal{S}' \longrightarrow \C$ {\rm(}see \eqref{rf}$)$ is a moment map for the natural
$S^1$-action and holomorphic symplectic structure $\Omega_0$ on $\mathcal{S}'$, in other words,
${\rm d}\mathcal{E}(-) = \i \Omega_0(X, -)$.
In particular, the $S^1$-fixed points in $\mathcal{S}'$ are precisely the critical points of
the functional $\mathcal{E}$.
\end{Theorem}

\begin{proof}
Take any $s \in \mathcal{S}'$, and set $m := s(0)$.
We first compute ${\rm d}_s \mathcal{E}(V)$ for
\[V \in T_{\operatorname{ev}_0}\mathcal S'|_s \subset T_s \mathcal{S}'
 = H^0\bigl(\C P^1, N_s\bigr) .\]
By definition, $V(0) = 0$, so $V(\lambda) = v\otimes\lambda$ for some $v \in H^0(\C P^1, s^*T_{\varpi}Z(-1))$
(see Lemma \ref{lem:TSEH}).
Further, $V$ is representable by a family $s_t$ of sections with $s_t(0) = m$ and
$\partial_{t=0} s_t(\lambda) = V(\lambda)$.
Note that $\partial_{\lambda = 0} V(\lambda)$ is well-defined because $V(0) = 0$.
Then a local computation shows that
\[\partial_{t = 0} \dot{s}_t(0) = \partial_{\lambda = 0} V(\lambda) = v \in T_m M.\]
Hence we conclude that
\[
{\rm d}_s \mathcal{E}(V) = \frac{\partial}{\partial t}\Big\vert_{t = 0} \left( -\frac{1}{2} \iota_Y
\omega\left(\dot{s}_t(0) - \dot{s}_m(0)\right) + \mu(m) \right) = -\frac{1}{2} \iota_Y \omega(v).
\]

Next we consider $X_s \in H^0\bigl(\C P^1, N_s\bigr)$, the fundamental vector field $X$ evaluated at $s$.
Since $Y(m) = X_s(0)$ (see Lemma \ref{RelationXandY}), it follows that
\[
\Omega_0|_s(X_s, V) = \frac{\i}{2}\frac{\partial}{\partial\lambda}\Big|_{\lambda=0}\omega_{s(\lambda)}(X_s(\lambda), V(\lambda))
 = \frac{\i}{2}\omega(Y(0), v) = \frac{\i}{2}\iota_Y\omega(v).
\]

It remains to prove the claim for $V \in T_{\operatorname{ev}_\infty}{\mathcal S}'\big\vert_s$.

To this end, observe that ${\rm d}(\overline{\tau^*\mathcal{E}}) = {\rm d}\mathcal{E}$ by \eqref{EnergyResidue} because
the function $s \longmapsto \mathrm{deg}(s^*L_Z)$ on ${\mathcal S}'$ is locally constant.
Further, we know that $\overline{{\rm d}\tau(T_{\operatorname{ev}_\infty}\mathcal S')} = \tau^*T_{\operatorname{ev}_0}\mathcal S'$,
so that every \smash{$V \in T_{\operatorname{ev}_\infty}\mathcal S'\big\vert_s$} is of the form \smash{${\rm d}\tau_{\tau(s)}(W)$} for \smash{$W
 \in \overline{T_{\operatorname{ev}_0}\mathcal S'\big\vert_{\tau(s)}}$}.
Then we compute the following using the previous result and $\tau^*\Omega = \overline{\Omega}$:
\begin{align*}
{\rm d}\mathcal{E}_s(V) & = \overline{\tau^*{\rm d}\mathcal{E}}_{s}(V) = {\rm d}\overline{\mathcal{E}}_{\tau(s)} (W) = \overline{\Omega}_{\tau(s)}(X_{\tau(s)}, {\rm d}\tau_{\tau(s)}(V)) \\
	&= \bigl(\overline{\tau^*\Omega}\bigr)_s(X_s, V) = \Omega_s(X_s, V).
\end{align*}
This completes the proof.
\end{proof}

\subsection[Critical Points of E]{Critical Points of $\boldsymbol{\mathcal E}$}\label{ss:C*}

In this subsection, we assume that the $S^1$-action on $Z$ extends to a holomorphic action of~$\C^*$
on~$Z$. This means that the vector field $I_ZY$ on $Z$ is complete. Since $\mathcal E$
in \eqref{rf} is holomorphic and $S^1$-invariant,
it follows that $\mathcal E$ is in fact $\C^*$-invariant, under this assumption.
Thus the critical points of $\mathcal E$ are the $\C^*$-fixed points in $\mathcal S'$.

We first examine the $\C^*$-fixed points in $\mathcal{S}$.
Any $\C^*$-fixed point $s\in \mathcal S$ is characterized by~${
s(\zeta\lambda) = \zeta.s(\lambda)}
$
for all $\zeta \in \C^*$ and $\lambda \in \mathbb CP^1$.
In particular, $s \in \mathcal{S}^{\C^*}$ is determined by its value at $\lambda = 1$.
Indeed, $s(\lambda) = \lambda.s(1)$ for $\lambda \in \C^*$, and by continuity we have
\begin{equation}\label{eq:fixedpointlimits}
s(0) = \lim_{\lambda\to 0} \lambda.s(1),\qquad s(\infty) = \lim_{\lambda \to \infty} \lambda.s(1).
\end{equation}
Hence the closures of the $\C^*$-orbits in $Z$ lying over $\C^* \subset \mathbb CP^1$ correspond precisely to the
$\C^*$-fixed points in $\mathcal{S}$.

Conversely, any point $z \in Z_{1}$ potentially determines a $\C^*$-invariant section $s^z \colon \C^*
 \longrightarrow Z$ of~$\varpi$ as follows.
For $\lambda \in \C^*$, set
$
s^z(\lambda) = \lambda.z
$
which is clearly a section over $\C^*$.
If the limits
$
\lim_{\lambda\to 0}\lambda.z$, $ \lim_{\lambda\to\infty}\lambda.z$,
exist in $Z$ (see \eqref{eq:fixedpointlimits}), then the section extends to a $\C^*$-invariant
section
\[s^z\colon\ \C P^1 \longrightarrow Z\]
 of $\varpi$. The existence of these limits has been investigated
in detail in \cite{Simpson-iterated} for $\mathcal{M}_{\mathrm{DH}}$; see
also Section \ref{ss:explicitC*} below.

Clearly, for any fixed point $s \in \mathcal S^{\C^*}$, we have $s(0) \in Z_0 = M$ and $s(\infty)
 \in Z_\infty \cong \overline{M}$ are fixed points of the $\C^*$-actions on $M$ and $\overline{M}$ respectively.
The following gives the converse on $\mathcal{S}' \subset \mathcal{S}$ (also see \cite{Fei01,Fei02,Fei03}).

\begin{Proposition}\label{prp1}
Let $s \in \mathcal S'$ be such that $s(0) \in Z_0^{\C^*}$ and $s(\infty) \in Z_\infty^{\C^*}$. Then
$s \in (\mathcal S')^{\C^*}$.
\end{Proposition}

\begin{proof}
We have $X_s(\lambda) = Y_{s(\lambda)} - \i \lambda\dot{s}(\lambda)$, which implies that $X_s(0) = 0 = X_s(\infty)$.
Since $N_s \cong \mathcal{O}_{\C P^1}(1)^{\oplus 2d}$, this implies that $X_s = 0$, and thus $s$ is a fixed point.
\end{proof}

Proposition \ref{prp1} has the following consequence.

\begin{Corollary}
Let $s \in \mathcal S$ be a section such that $s(0) \in Z_0^{\C^*}$ and $s(\infty) \in Z_\infty^{\C^*}$.
If $s$ is not a fixed point of the $\C^*$-action on $\mathcal S$, then $s$ cannot be contained in $\mathcal S'$, meaning
the normal bundle of $s$ is not isomorphic to $\mathcal{O}_{\C P^1}(1)^{\oplus 2d}$.
\end{Corollary}

We end this subsection with the observation that $\C^*$-fixed points $s \in \mathcal S$ are also
the fixed points under the twisting procedure that was introduced in \cite{BeHeRo}; see also \cite{BHR}.
Recall that a~section~$s\in \mathcal{S}$ is twistable if the section $\widetilde{s}(\lambda) := \lambda^{-1}.s\bigl(\lambda^2\bigr)$ over $\C^*$ extends to a section on all of~$\mathbb CP^1$.

\begin{Proposition}\label{Prop: FixedPointTwist}
Let $s$ be a fixed point for the $\C^*$-action on $\mathcal S$. Then $s$ is twistable, and the twist $\widetilde{s}$ satisfies
the equation $\widetilde{s} = s$.
\end{Proposition}

\begin{proof}
If $s$ is fixed, then $s(\lambda) = \lambda.s(1)$ for all $\lambda \in \C^*$. If follows that
\[\widetilde s(\lambda) = \lambda^{-1}.s\bigl(\lambda^2\bigr) = \lambda^{-1}.\lambda^2.s(1) = \lambda.s(1) = s(\lambda),\]
which completes the proof.
\end{proof}

\subsection[The Atiyah--Ward transform of L\_Z]{The Atiyah--Ward transform of $\boldsymbol{L_Z}$}\label{ss:AW}

Let $\mathcal S^0 \subset \mathcal S'$ stand for the space of all sections $s \in \mathcal S$ such that
\begin{itemize}\itemsep=0pt
\item the normal bundle is isomorphic to $\mathcal{O}_{\C P^1}(1)^{\oplus 2d}$, and

\item the pullback $s^*L_Z$ trivial.
\end{itemize}
Since $s^*L_Z$ trivial if and only if $\deg s^*L_Z = 0$, we conclude that
$\mathcal S^0$ is an open subset of $\mathcal S$, and it is in fact a union of some
connected components of $\mathcal S'$.
Consider the space $\mathcal F^0 = \C P^1 \times \mathcal S^0$ and restrict the diagram\footnote{Here we slightly abuse notation,
however the restriction of the maps will be evident in the following.} \eqref{eq: doublefibration} to $\mathcal F^0$
\begin{equation}\label{td}
\begin{tikzcd}
Z \arrow[r, leftarrow, "\operatorname{ev}"] \arrow[d, "\varpi"'] & {\mathcal F}^0 = \C P^1 \times \mathcal{S}^0 \arrow[dl, "\pi_1"'] \arrow[d, "\pi_2"] \\
\C P^1 & \mathcal{S}^0.
\end{tikzcd}
\end{equation}
Clearly, we have the identification
$
T\mathcal F^0 = \pi_1^*T\C P^1 \oplus \pi_2^*T\mathcal S^0
$
(see \eqref{td}). Since ${\rm d}\operatorname{ev}\colon \pi_2^*T\mathcal S' \longrightarrow \operatorname{ev}^*T_\varpi Z$
is surjective by Proposition \ref{p:local biholo}, we obtain the following commutative diagram with exact rows:
\begin{equation}\label{TFTS}
\begin{tikzcd}
0 \arrow[r] & T_{\operatorname{ev}}\mathcal F^0 \arrow[r] & T\mathcal F^0 \arrow[r,"{\rm d}\operatorname{ev}"] & \operatorname{ev}^*TZ \arrow[r] & 0\\
0 \arrow[r] & T_{\operatorname{ev}}\mathcal F^0 \arrow[u, "\mathrm{id}"]\arrow[r] & \pi_2^*T\mathcal S^0 \arrow[u, hookrightarrow] \arrow[r,"{\rm d}\operatorname{ev}"] &\arrow[u,hookrightarrow] \operatorname{ev}^*T_\varpi Z \arrow[r]& 0.
\end{tikzcd}
\end{equation}
We next describe the Atiyah--Ward transform of $L_Z$ \cite{AW,BaEa,Merkulov} and study how it interacts with the
meromorphic connection on $L_Z$. The Atiyah--Ward construction is of course valid in a more general
context, but for the convenience of the reader we spell out the details relevant for our discussion for the
special case of a line bundle.

The Atiyah--Ward transform of $L_Z$ is a holomorphic line bundle $\mathcal{L}$ over $\mathcal{S}^0$ with a
holomorphic connection $\nabla^{\rm AW}$. To construct this line bundle, recall that $L_Z$ is trivial along $s\bigl(\C P^1\bigr)$
for any~${s \in \mathcal S^0}$. Hence $\operatorname{ev}^*L_Z$ is trivial along $\C P^1\times\{s\} = \pi_2^{-1}(s)$, and the
($0$-th) direct image construction yields the line bundle
$\mathcal L := (\pi_2)_*\operatorname{ev}^*L_Z$.
We observe that there is a natural isomorphism $\pi_2^*\mathcal{L} \cong \operatorname{ev}^*L_Z$. Over any point $(\lambda, s) \in
 \mathcal{F}^0$, the isomorphism is given by evaluating an element of $\pi_2^*\mathcal{L}_{(\lambda, s)} =
 H^0\bigl(\C P^1, s^*L_Z\bigr)$ at $\lambda$.

We next equip $\mathcal{L}$ with a holomorphic connection.
First consider the relative exterior differential ${\rm d}_{\operatorname{ev}}\colon \mathcal O_{\mathcal F^0}
\longrightarrow T^*_{\operatorname{ev}}\bigl(\mathcal F^0\bigr)$ \big(see \eqref{td} for ${\mathcal F}^0$\big) which is defined as the
following composition of maps:
\[
\begin{tikzcd}
{\rm d}_{\operatorname{ev}} \colon \mathcal O_{\mathcal F^0} \arrow[r, "d"] & T^*\mathcal F^0 \arrow[r] & T^*_{\operatorname{ev}}\bigl(\mathcal F^0\bigr).
\end{tikzcd}
\]
By construction, we have $\operatorname{ev}^*\mathcal O_Z \subset \ker {\rm d}_{\operatorname{ev}}$, and for any $f
 \in \mathcal O_{\mathcal F^0}$ we have
$
{\rm d}_{\operatorname{ev}}f = {\rm d}f|_{T_{\operatorname{ev}}\mathcal F^0}$.
Next we tensor the diagram defining ${\rm d}_{\operatorname{ev}}$ with $\operatorname{ev}^*L_Z$ and observe that this gives a well-defined ``relative connection''
\[\nabla_{\operatorname{ev}} \colon\ \operatorname{ev}^*L_Z \longrightarrow T^*_{\operatorname{ev}}\mathcal F^0 \otimes
\operatorname{ev}^*L_Z,\qquad
\nabla_{\operatorname{ev}}(f\cdot\operatorname{ev}^*\nu) := ({\rm d}_{\operatorname{ev}}f)\otimes\operatorname{ev}^*\nu,
\]
for a locally defined trivializing holomorphic section $\nu$ of $L_Z$ and $f \in \mathcal O_{\mathcal F^0}$.
This operator is well-defined because the transition functions for $\operatorname{ev}^*L_Z$ lie in $\operatorname{ev}^*\mathcal O_Z$
and hence they are in the kernel of ${\rm d}_{\operatorname{ev}}$.
The next lemma will enable us to construct $\nabla^{\rm AW}$ from $\nabla_{\operatorname{ev}}$.

\begin{Lemma}\label{isoS0}
There is a natural isomorphism of vector bundles over $\mathcal S^0$
\[
T^*\mathcal S^0 \cong (\pi_2)_* T^*_{\operatorname{ev}}\mathcal F^0, \qquad A \longmapsto
\pi_2^*A_{|T_{\operatorname{ev}}\mathcal{F}^0}
\]
$($see \eqref{td}$)$. In particular, there is an isomorphism
\[
T^*\mathcal{S}^0 \otimes \mathcal{L} \cong (\pi_2)_* \bigl(T^*_{\operatorname{ev}}\mathcal F^0 \otimes \operatorname{ev}^*L_Z\bigr).
\]
\end{Lemma}

\begin{proof}
Dualizing the second row of \eqref{TFTS} yields

\begin{tikzcd}
0\arrow[r] & \operatorname{ev}^*T^*_\varpi Z\arrow[r] & \pi_2^*T^*\mathcal S^0 \arrow[r, "r"] & T^*_{\operatorname{ev}}\mathcal F^0\arrow[r] & 0.
\end{tikzcd}
Here $r$ denotes the restriction of $1$-forms to $T_{\operatorname{ev}}^*\mathcal{F}^0$. For this
short exact sequence, consider the corresponding
long exact sequence of direct image sheaves with respect to $\pi_2$
\[
\begin{tikzcd}
0\arrow[r] & (\pi_2)_*\operatorname{ev}^*(T^*_\varpi Z) \arrow[r] \ar[draw=none]{d}[name=X, anchor=center]{} &
(\pi_2)_*(\pi_2^*T^*\mathcal S^0)
\ar[rounded corners, to path={ -- ([xshift=2ex]\tikztostart.east)
 |- (X.center) \tikztonodes
 -| ([xshift=-2ex]\tikztotarget.west)
 -- (\tikztotarget)}]{dl}[near start]{\pi_{2*}r} \\
& (\pi_2)_*(T^*_{\operatorname{ev}}
\mathcal F^0) \arrow[r] & (\pi_2)^1_*\operatorname{ev}^*(T^*_\varpi Z) \arrow[r] & \cdots.
\end{tikzcd}
\]
The fibers of the vector bundles $ (\pi_2)^q_*\operatorname{ev}^*(T_\varpi^*Z)$, $q = 0, 1$, over
any $s \in \mathcal{S}^0$ satisfy
\[
 H^q\bigl(\C P^1, s^*T_\varpi^*Z\bigr) \cong H^q\bigl(\C P^1, \mathcal{O}(-1)^{\oplus 2d}\bigr) = 0.
\]
Thus the above long exact sequence of cohomologies gives the isomorphism
\[
\pi_{2*}r \colon (\pi_2)_*\bigl(\pi_2^*T^*\mathcal S^0\bigr) \stackrel{\cong}{\longrightarrow} (\pi_2)_*\bigl(T^*_{\operatorname{ev}}\mathcal F^0\bigr)
\]
of vector bundles on $\mathcal S^0$.

By the projection formula (see, e.g., \cite[Chapter 3]{Hart}) or by a direct computation, we~have
\[
(\pi_{2})_*\bigl(\pi_2^*T^*\mathcal{S}^0 \otimes \mathcal{O}_{\mathcal{F}^0}\bigr)
 \cong T^*\mathcal{S}^0 \otimes (\pi_2)_*\mathcal{O}_{\mathcal{F}^0} \cong T^*\mathcal{S}^0.
\]
The last isomorphism follows from the fact that
$(\pi_2)_*\mathcal{O}_{\mathcal{F}^0} \cong \mathcal{O}_{\mathcal{S}^0}$, which in turn follows
from the fact that the fibers of $\pi_2$ are connected. Since we have
$\operatorname{ev}^*L_Z \cong \pi_2^*\mathcal{L}$ over $\mathcal{S}^0$, the second
statement in the lemma is again derived from the projection formula.
\end{proof}

\begin{Proposition}\label{Prop: AtiyahWard}
The operator
\[
\nabla^{\rm AW} = (\pi_2)_*(\nabla_{\operatorname{ev}})\colon\ \mathcal L
 \longrightarrow T^*\mathcal{S}^0\otimes\mathcal L
 \cong (\pi_2)_*\left(T^*_{\operatorname{ev}}\mathcal F^0\otimes\operatorname{ev}^*L_Z\right)
\]
induces a natural holomorphic connection on $\mathcal L$.
This connection is trivial on the submanifolds~${\mathcal S_z^0 = \pi_2\bigl(\operatorname{ev}^{-1}(z)\bigr)}$ for all $z \in Z$.
This property determines $\nabla^{\rm AW}$ completely as long as $\mathcal S_z^0$ is connected.
\end{Proposition}

We call the above connection $\nabla^{\rm AW}$ on $\mathcal L$ the \emph{Atiyah--Ward connection}.

\begin{proof}
The $\mathcal O_{\mathcal S^0}$-module structure on the sheaf of sections of the
vector bundle $\mathcal L = (\pi_2)_*\operatorname{ev}^*L_Z$ is as follows. Let $S \subset
\mathcal S^0$ be an open subset, and take any $f \in \mathcal O_{\mathcal S^0}(S)$.
Then $f$ acts on~${\mathcal L(S) = H^0\bigl(\C P^1\times S, \operatorname{ev}^*L_Z\bigr)}$ by multiplication
with $\pi_2^*f$. Hence for any $\psi \in \mathcal L(S)$ we have
\[
\nabla^{\rm AW}(f\psi) =
\nabla_{\operatorname{ev}}(\pi_2^*f\psi) = {\rm d}_{\operatorname{ev}}\pi_2^*f\otimes\psi +
f\nabla_{\operatorname{ev}}\psi.
\]
We see that ${\rm d}_{\operatorname{ev}}\pi_2^*f\otimes\psi + f\nabla_{\operatorname{ev}}\psi$ is obtained from
$\pi_2^*\bigl({\rm d}f\otimes\psi + f\nabla^{\rm AW}\psi\bigr)$ by restricting to the subbundle
$T_{\operatorname{ev}}\mathcal F^0$. Note that this is uniquely determined by
Lemma \ref{isoS0}, so we conclude that{\samepage
\[
\nabla^{\rm AW}(f\psi) = {\rm d}f\otimes\psi + f\nabla\psi.
\]
Therefore, $\nabla^{\rm AW}$ is indeed a connection.}

To show that this connection $\nabla^{\rm AW}$ is trivial on $\mathcal{S}_z^0$ for each $z \in Z$, observe
that any section of $\operatorname{ev}^* L_Z \longrightarrow {\mathcal F}^0$ (see \eqref{td}) of the form
$\operatorname{ev}^*\psi$ is actually parallel
(covariant constant) for the relative connection
$\nabla_{\operatorname{ev}}$ by the formula
$\nabla_{\operatorname{ev}}(f\cdot\operatorname{ev}^*\psi) =
({\rm d}_{\operatorname{ev}}f)\otimes\operatorname{ev}^*\psi $.

Now the fibers of $\operatorname{ev}$ are of the form
$
\operatorname{ev}^{-1}(z) = \{\varpi(z)\}\times\mathcal S_z^0$.
Then we get a~frame~$\psi$ of~$\operatorname{ev}^*L_Z|_{\operatorname{ev}^{-1}(z)}$ by putting
$
\psi(\varpi(z), s) = l
$
for any fixed $0 \neq l \in (L_Z)_z$.
Note that this construction makes sense because we have
$\operatorname{ev}(\varpi(z),s) = s(\varpi(z)) = z$ for all $s \in
\mathcal S_z^0$. The frame~$\psi$ induces a~natural trivialization $\psi^{\mathcal L}$ of
$\mathcal L$ along $\mathcal S_z^0$. Since $\nabla_{\operatorname{ev}}\psi = 0$, we have $\nabla^{\rm AW}
\psi^{\mathcal L} = 0$ along~$\mathcal S_z^0$. Thus~$\nabla^{\rm AW}$ is trivial on $\mathcal S_z^0$.

Now let $\widetilde\nabla = \nabla^{\rm AW} +A$ be another holomorphic connection on $\mathcal{L}$
which is trivial along each~$\mathcal S_z^0$. Then the holomorphic $1$-form $A$ restricted to
$\mathcal{S}_z^0$ is exact for all $z \in Z$. Then a result of Buchdahl in \cite{Buchdahl} implies that
we can find $f \in \mathcal O\bigl(\mathcal F^0\bigr)$ such that \smash{${\rm d}_{\operatorname{ev}}f =
(\pi_2^*A)|_{T_{\operatorname{ev}}\mathcal F^0}$}. Note that $f$ is constant along the compact fibers
of $\pi_2$. Therefore, there exists \smash{$\widetilde f \in \mathcal O\bigl(\mathcal S^0\bigr)$} such that we have
\smash{$
f = \widetilde{f}\circ \pi_2
$}
if $\mathcal S_z^0$ is connected. It now follows that
\[
\pi_2^*\bigl(A-{\rm d}\widetilde f\bigr)|_{T_{\operatorname{ev}}\mathcal F^0} = 0,
\]
and this implies that $A = {\rm d}\widetilde{f}$. Thus, $\nabla^{\rm AW}$ and $\widetilde\nabla =
\nabla^{\rm AW} +{\rm d}\widetilde{f}$ are actually gauge-equivalent.
\end{proof}

\begin{Corollary}
On the real submanifold $M \subset \mathcal S^0$ the connection $\nabla^{\rm AW}$
coincides with the hyperholomorphic connection on $L_M$.
Thus the curvature $F^{\rm AW}$ of $\nabla^{\rm AW}$ satisfies the equation
\[
F^{\rm AW} = \Omega_0 + \i\, {\rm d}I_0{\rm d}\mathcal E
\]
on every component of $\mathcal S^0$ that meets $M$ $($see \eqref{rf} for $\mathcal E)$.
\end{Corollary}

\begin{proof}
The evaluation map $\operatorname{ev}\colon \mathcal F^0 \longrightarrow
Z$ restricts to a diffeomorphism $\epsilon\colon M\times \C P^1 \longrightarrow
Z$. Thus we have a natural map
$q = \pi_2\circ\epsilon^{-1}\colon Z
 \longrightarrow M$.
Note that
\[q^*\mathcal L = q^*(\pi_2)_*\operatorname{ev}^*L_Z
 = \bigl(\epsilon^{-1}\bigr)^*\pi_2^*(\pi_2)_*\operatorname{ev}^*L_Z
 = \bigl(\epsilon^{-1}\bigr)^*\operatorname{ev}^*L_Z = L_Z.\]
On the other hand, we have, by definition, that $q^*L_M = L_Z$.

Now the hyperholomorphic connection $\nabla_M$ is real analytic. Thus we can complexify and
extend it to a holomorphic connection $\widetilde\nabla$ on $\mathcal
L \longrightarrow \mathcal S^0$ at least locally over a neighbourhood of~${M \subset
\mathcal S^0}$. Note that its curvature is $\Omega_0 + \i \, {\rm d}I_0 {\rm d} \mathcal E$, which is in
fact the complexification of the curvature of $\nabla_M$.

Since there is a unique twistor line through each point in $z$, we see that ${\mathcal S}_z^0$ intersects $M$
in a unique point, namely $q(z)$. Moreover, as the twistor line $q(z)$ through $z$ also passes through~$\tau_Z(z)$, we see that \smash{$q(z) = \mathcal S_z^0\cap \mathcal S_{\tau(z)}^0\cap M$}. We obtain a splitting
\[
T_{q(z)}\mathcal S^0 = T_{q(z)}\mathcal S_z^0 \oplus T_{q(z)}\mathcal S_{\tau(z)}^0,
\]
which can be identified with the splitting \smash{$T_{q(z)}M\otimes \C = T^{1,0,\lambda}M \oplus T^{0,1,\lambda}M$}.
Since $\nabla_M$ has curvature of type $(1,1)$ on $(M, I_\lambda)$ (here~$\lambda
 = \varpi(z)$), its complexification $\widetilde\nabla$ is flat along $\mathcal S_z^0$
\big(and~\smash{$\mathcal S_{\tau(z)}^0$}\big) for all $z \in Z$.
By Proposition \ref{Prop: AtiyahWard}, the connections $\nabla$ and $\widetilde \nabla$ must be
gauge equivalent over a neighbourhood\footnote{Note that there exists a neighbourhood $M
 \subset U \subset \mathcal{S}$ such that $\mathcal{S}_z\cap U$ is connected for
every $z \in Z$, because~${M\cap \mathcal{S}_z = \{ q(z)\}}$.} of $M$ in $\mathcal S^0$.

The curvature $F^{\rm AW}$ is a holomorphic two-form on $\mathcal S^0$ which is
compatible with the real structure. On the other hand, the form $\Omega_0 + {\rm d}I_0{\rm d}\mathcal E$ is also
holomorphic and real and coincides with $F^{\rm AW}$ on the real submanifold $M$. Thus, the two 2-forms
must coincide on every component of $\mathcal S^0$ that meets $M$.
\end{proof}

\begin{Corollary}\label{om0curaw}
The curvature $F^{\nabla^0}$ of the holomorphic connection
$\nabla^0 = \nabla^{\rm AW} -\i I_0{\rm d}\mathcal E$ satisfies the equation
$F^{\nabla^0} = \Omega_0$ on every component of $\mathcal{S}^0$ that intersects $M$, where
$\mathcal E$ is defined in~\eqref{rf}.
\end{Corollary}

\subsection{The Atiyah--Ward transform and the meromorphic connection}

We want to describe the relationship between the meromorphic connection $\nabla$ on $L_Z$ and the
Atiyah--Ward transform $\bigl(\mathcal{L}, \nabla^{\rm AW}\bigr)$ of $L_Z$. As the first step, we observe
that the Atiyah class~$\operatorname{ev}^*\alpha_L$ of $\operatorname{ev}^*L_Z
 \longrightarrow \mathcal F^0$ vanishes because
$\operatorname{ev}^*L_Z = \pi_2^*\mathcal{L}$ admits the holomorphic connection~$\pi_2^*\nabla^{\rm AW}$.
This has the following implications: On $\mathcal F^0$ consider
\[\widehat\sigma
 = \operatorname{ev}^*\sigma \in H^0\bigl({\mathcal F}^0, \pi_1^*{\mathcal O}_{{\mathbb C}P^1}(2)\bigr),\]
which vanishes on the divisor
\[
\widehat{D} = \operatorname{ev}^{-1}(D) = \bigl(\{0\}\times\mathcal S^0\bigr)\cup \bigl(\{\infty\}\times\mathcal S^0\bigr).
\]
We have the short exact sequence of sheaves
\[
\begin{tikzcd}
0 \ar[r] &T^*\mathcal{F}^0 \ar[r, "\cdot \widehat\sigma"] & T^*\mathcal{F}^0(2) \ar[r] &
T^*\mathcal{F}^0(2)_{|\widehat{D}} \ar[r] & 0
\end{tikzcd}
\]
with the corresponding long exact sequence of cohomologies
\begin{equation}\label{eq:les}
\begin{tikzcd}
0 \ar[r] & H^0\bigl(\mathcal F^0, T^*\mathcal F^0\bigr) \ar[r, "\cdot \widehat\sigma"] \ar[draw=none]{d}[name=X, anchor=center]{} & H^0\bigl(\mathcal F^0,
T^*\mathcal F^0(2)\bigr)
\ar[rounded corners, to path={ -- ([xshift=2ex]\tikztostart.east)
 |- (X.center) \tikztonodes
 -| ([xshift=-2ex]\tikztotarget.west)
 -- (\tikztotarget)}]{dl}[near start]{} \\
 & H^0\bigl(\widehat{D}, T^*\mathcal F^0(2)_{|\widehat{D}}\bigr)
\ar[r, "\widehat{\delta}"] & H^1\bigl(\mathcal F^0, T^*\mathcal F^0\bigr)\ar[r]&\cdots.
\end{tikzcd}
\end{equation}
The construction of $L_Z$ implies that $\operatorname{ev}^*\alpha_L
 = \widehat{\delta}(\operatorname{ev}^*\varphi)$, which vanishes by the previous observation.
Hence there exists $\phi \in H^0\bigl(\mathcal F^0, T^*\mathcal F^0(2)\bigr)$ such that
$\phi|_{\widehat{D}} = \operatorname{ev}^*\varphi$ by the long exact sequence in \eqref{eq:les}.

\begin{Proposition}\label{p:holconn}
The connection $\operatorname{ev}^*\nabla - \frac{\phi}{\widehat\sigma}$ on
$\operatorname{ev}^*L_Z = \pi_2^*\mathcal{L}$ is holomorphic.
\end{Proposition}

\begin{proof}
From the construction of $\operatorname{ev}^*\nabla$, we know that it has connection $1$-forms
\[
\widehat{A}_i = \frac{\operatorname{ev}^*\widetilde{\varphi}_i}{\widehat\sigma}
\]
with respect to the open covering $\operatorname{ev}^*\mathcal U = \bigl\{\operatorname{ev}^{-1}(U_i)\bigr\}$.
Here $\mathcal U = \{U_i\}$ and \smash{$\widetilde{\varphi}_i$} are as in the construction of
$\nabla$ in Section \ref{ss:rotating}.
The connection \smash{$\operatorname{ev}^*\nabla - \frac{\phi}{\widehat\sigma}$} has the connection $1$-forms~\smash{$
\widehat{A}_i - \frac{\phi_i}{\widehat\sigma}$}.
But \smash{$\phi_i|_{\widehat D}$} and \smash{$\operatorname{ev}^*\varphi_i$} coincide on \smash{$\operatorname{ev}^{-1}(U_i) \cap
\widehat{D}$} and the forms $\widehat{A}_i$ have no other poles.
Hence the proposition follows.
\end{proof}

We next explicitly construct such a section $\phi \in H^0\bigl(\mathcal F^0,
T^*\mathcal F^0(2)\bigr)$ and describe the difference form $B$ between the holomorphic
connections \smash{$\operatorname{ev}^*\nabla -\frac{\phi}{\widehat\sigma}$} and \smash{$\pi_2^*\nabla^{\rm AW}$}.
To do so, we first observe that
\begin{equation}\label{eq:TW}
T^*\mathcal F^0(2) \cong \bigl(\pi_1^*T^*\C P^1 \oplus \pi_2^*T^*\mathcal{S}^0\bigr)(2)
 \cong \mathcal{O}_{\mathcal F^0} \oplus \pi_2^*T^*\mathcal{S}^0(2).
\end{equation}
We describe $\phi$ component-wise, that is, $\phi = (\beta, \gamma)$ with respect to the
splitting \eqref{eq:TW}.
Given any $V \in T_{s} \mathcal{S}^0$,
construct a section $\gamma \in H^0\bigl(\mathcal F^0, \pi_2^*T^*\mathcal{S}^0\bigr)(2))$ in the following
way:
\[
\gamma_{(s,\lambda)} (V) := -\frac{1}{2}\omega(X_s(\lambda), V(\lambda))
 = -\frac{1}{2}(\iota_X\operatorname{ev}^*\omega)_{(\lambda,s)}(V).
\]
The first component $\beta$ of $\phi$ is
\smash{$
\beta := \mathcal E \mathrm{id}_{\pi_1^*\mathcal O(2)}$},
where $\mathcal E$ is defined in \eqref{rf}, which is in
\[H^0\bigl(\mathcal F^0, \mathrm{End}(\pi_1^*\mathcal{O}(2))\bigr)
 = H^0\bigl(\mathcal F^0, \pi_1^*T^*\C P^1(2)\bigr)
 \cong H^0\bigl(\mathcal F^0, \mathcal O_{\mathcal F^0}\bigr).\]

\begin{Lemma}\label{Lemma:holoconn}
The element $\phi = (\beta, \gamma) \in H^0\bigl(\mathcal F^0, T^*\mathcal F^0(2)\bigr)$
constructed above has the desired properties. Moreover, $\phi|_{T_{\operatorname{ev}}\mathcal F^0}
 = 0$.
\end{Lemma}

\begin{proof}
We may write $\beta = \mathcal E {\rm d}\lambda\otimes\frac{\partial}{\partial\lambda}$
on $\{\lambda\neq\infty\}$, and hence at $\lambda = 0$
\[
\operatorname{ev}^*\varphi_{(s,0)}\left(\frac{\partial}{\partial \lambda}, 0\right) =
\varphi_{s(0)} (\dot{s}(0)) = \mathcal{E}(s) \frac{\partial}{\partial \lambda}
 = \beta_{(s,0)}\left(\frac{\partial}{\partial\lambda}\right).
\]
Similarly, at $\lambda = \infty$,
\[
\operatorname{ev}^*\varphi_{(s,\infty)}\left(\frac{\partial}{\partial\widetilde \lambda}, 0\right) =
\varphi_{s(\infty)} (\dot{s}(\infty)) = \mathcal{E}(s) \frac{\partial}{\partial
\widetilde\lambda} = \beta_{(s,\infty)}\left(\frac{\partial}{\partial\lambda}\right).
\]
Now we check the component $\gamma$. Let $V \in T_s\mathcal S^0$. Then, by \eqref{eq:varphiomega},
\[
\operatorname{ev}^*\varphi_{(s,0)}(V) = \varphi_{s(0)}(V(0)) = -\frac{1}{2} (\iota_Y\omega)_{s(0)}(V(0)).
\]
On the other hand, by \eqref{Eq: XsYsdots},
\begin{align*}
-\frac{1}{2}(\iota_X\operatorname{ev}^*\omega)_{(s,0)}(V)& = -\frac{1}{2}\omega(X(0), V(0)) = -\frac{1}{2}\omega(Y(0), V(0)) \\
& = -\frac{1}{2}\iota_Y\omega(V(0)) = \operatorname{ev}^*\varphi_{(s,0)}(V).
\end{align*}
Similar considerations apply at $\lambda = \infty$. Since $\gamma
 = -\frac{1}{2}\iota_X(\operatorname{ev}^*\omega)|_{T\mathcal S^0}$ and $T_{\operatorname{ev}}\mathcal F^0
 \subset \pi_2^*T\mathcal S^0$, it follows that $\phi|_{T_{\operatorname{ev}}\mathcal F^0} = 0$.
\end{proof}

To understand the relationship of the holomorphic connection
\smash{$
\widehat\nabla := \operatorname{ev}^*\nabla - \frac{\phi}{\widehat\sigma}
$}
 with the Atiyah--Ward connection, we need the following lemma.

\begin{Lemma}\label{Lemma:holconnpullback}
There is a unique holomorphic connection $\nabla^{\mathcal{L}}$ on $\mathcal{L}$ such that
$\pi_2^*\nabla^{\mathcal L} = \widehat{\nabla}$ on~${\pi_2^*\mathcal L = \operatorname{ev}^*L_Z}$.
\end{Lemma}

\begin{proof}
We will prove the more general statement that any holomorphic connection on $\pi_2^*\mathcal L$
is pulled back from a unique holomorphic connection on $\mathcal L$. Let $F$ be the curvature
of such a~connection which is a holomorphic two-form on $\mathcal{F}^0$. Further, let $V$ be a
vector field on $\mathcal S^0$ which we view as a vector field on $\mathcal F^0 = \C P^1\times
\mathcal S^0$ by assigning zero vector field along $\C P^1$.
Then~$\iota_VF^\mathsf{D}$ pulls back to a holomorphic $1$-form on $\C
P^1\times\{s\}$, which must vanish, since there are no non-trivial holomorphic $1$-forms on $\C
P^1$. This shows that $F$ is purely horizontal with respect to the projection $\pi_2\colon \C P^1\times
\mathcal S^0 \longrightarrow \mathcal S^0$.

We can thus find local frames for $\pi_2^*\mathcal L$ which are parallel along the fibers of
$\pi_2$ \big(which are~$\mathbb CP^1$\big). More precisely, there exists an open cover of $\C P^1\times
\mathcal S^0$ by open subsets of the form $\C P^1\times U$, where~${U \subset \mathcal S^0}$ is open,
together with holomorphic frames for $\operatorname{ev}^*L_Z$ over $\C P^1\times U$ that are parallel
(covariant constant) in the $\C P^1$-direction. The transition functions with respect to
these local frames are
holomorphic functions $\C P^1\times (U\cap V) \longrightarrow \C^*$, which are actually
constant along the $\C
P^1$-fibers. Hence the transition functions are of the form $g_{UV}\circ \pi_2$, where the
collection~$\{g_{UV}\}$ defines the line
bundle $\mathcal L$ on $\mathcal S$. Moreover, the holomorphic connection on~$\pi_2^*\mathcal L$ uniquely descends to $\mathcal L$, concluding the proof.
\end{proof}

\begin{Theorem}\label{AWH}
The holomorphic connections $\nabla^{\mathcal{L}}$ and $\nabla^{\rm AW}$ on $\mathcal{L}$ coincide.
\end{Theorem}

\begin{proof}
Clearly, $\nabla^{\mathcal{L}} = \nabla^{\rm AW} + B$ for a unique holomorphic $1$-form $B$
on $\mathcal{S}^0$. Consider the pullback connection
\begin{equation}\label{eq:difference}
\widehat\nabla = \pi_2^*\nabla^{\mathcal{L}} = \pi_2^*\nabla^{\rm AW} + \pi_2^*B.
\end{equation}
Let $U \subset Z$ be an open subset such that $L_Z$ has a local frame $\nu$ on $U$. Hence
$\psi := \operatorname{ev}^*\nu$ is a~local frame of $\operatorname{ev}^*L_Z$ on $\widehat{U} :=
\operatorname{ev}^{-1}(U)$. Observe that along $\mathcal{S}^0_z$, for $z \in U$, the frame~$\psi(\varpi(z), -)$ coincides with the parallel frame \big(with respect to $\nabla^{\rm AW}$\big)
constructed in the proof of Proposition~\ref{Prop: AtiyahWard}. Since we have $\operatorname{ev}^{-1}(z) =
\{\varpi(z)\} \times \mathcal{S}^0_z$, it follows that $\pi_2^*\nabla^{\rm AW}\psi$ restricted to
$T_{\operatorname{ev}}\mathcal{F}^0 = \ker {\rm d}\operatorname{ev}$ actually vanishes on all of \smash{$\widehat{U}$}.

Next let $(\lambda, s) \in \widehat{U} - \widehat{D}$, and \smash{$X \in \bigl(T_{\operatorname{ev}}\mathcal{F}^0\bigr)_{(\lambda,s)}$}.
Then we compute
\[
\widehat{\nabla}_X \psi = (\operatorname{ev}^*\nabla)_X \psi - \frac{\phi}{\widehat{\sigma}}(X)\psi
 = 0.
\]
Here we made use of the definition of $T_{\operatorname{ev}}\mathcal{F}^0$ and Lemma \ref{Lemma:holoconn}.
By continuity, the restriction of~$\widehat{\nabla}\psi$ to
$T_{\operatorname{ev}}\mathcal{F}^0$ vanishes on all of $\widehat{U}$ as well.
Since $\mathcal{F}^0$ can be covered by open subsets $\widehat{U}$ as above, from
\eqref{eq:difference} it follows that $\pi_2^*B_{|T_{\operatorname{ev}}\mathcal{F}^0} = 0$.
But the isomorphism $T^*\mathcal{S}^0 \cong \pi_{2*}T_{\operatorname{ev}}\mathcal{F}^0$ is given
by $A \longmapsto \pi_2^*A_{|T_{\operatorname{ev}}\mathcal{F}^0}$ (cf.\ Lemma \ref{isoS0}),
so we have $B = 0$, which implies that $\nabla^{\mathcal{L}} = \nabla^{\rm AW}$.
\end{proof}

The existence of a holomorphic connection on $\mathcal S^0$ with curvature $\Omega_0$ combine together
with Theorem~\ref{thm:Omega} to allow us to obtain the following results on the global
structure of $\mathcal S^0$.

Recall
\[\mathcal{S}^0_z = \operatorname{ev}_{\varpi(p)}^{-1}(z) \cap
\mathcal{S}^0\]
for $z \in Z$.

\begin{Theorem}\label{Thm:S0neqZ0Zinfty}
Suppose $0 \neq [\omega_1] \in H^2(M, {\mathbb C})$, and let $z \in Z_\infty$.
Then $\mathcal{S}^0_z$ cannot be isomorphic to $Z_0$.
In particular, the local diffeomorphism
\[
\operatorname{ev}_0\times\operatorname{ev}_\infty
 \colon\ \mathcal S^0 \longrightarrow Z_0\times Z_\infty
\]
cannot extend to a global diffeomorphism.
\end{Theorem}

\begin{proof}
Since $T\mathcal S_z^0$ is Lagrangian for the symplectic form $\Omega_0$, the Atiyah--Ward
connection $\nabla^{\rm AW}$ pulls back to a flat connection on \smash{$\mathcal L^{(z)} :=
\mathcal L|_{\operatorname{ev}_\infty^{-1}(z)}$}. Consequently, we have $c_1\bigl(\mathcal{L}^{(z)}, {\mathbb C}\bigr)
 = 0 \in H^2\bigl(\operatorname{ev}_{\infty}^{-1}(z)\bigr)$.
But the first Chern class of $L_M$ is $[\omega_1] \neq 0 \in H^2(Z_0, {\mathbb C})$.
\end{proof}

\begin{Remark}
If $Z_{0} = M$ is not simply-connected,
then $\operatorname{ev}_0 \colon \mathcal S^0_p \longrightarrow Z_0$ might be a covering map.
In fact, this happens for the rank $1$ Deligne--Hitchin moduli space; see \cite{BH}.
It would be interesting to understand what happens if
 $Z_{0}$ is simply-connected, for example in the case of $\operatorname{SL}(2,\C)$-Deligne--Hitchin moduli spaces.
\end{Remark}

\subsection{An alternative construction of Hitchin's meromorphic connection}\label{subsec:alternativeconn}

We can reverse the above constructions to obtain the meromorphic Hitchin connection on a~hyperholomorphic bundle over the twistor space $Z$
of a hyperK\"ahler manifold with a rotating circle action.
We only sketch the steps, as the technical details are already explained above.
The necessary twistor data are
\begin{itemize}\itemsep=0pt
\item the twistor space $Z$ of a hyperK\"ahler manifold with a rotating circle action;
\item the twisted relative holomorphic symplectic form
$\omega \in H^0\bigl(Z, \bigl(\bigwedge^2 T_{\varpi}^*\bigr)(2)\bigr).$
\end{itemize}
Adding the additional structure of a ``prequantum'' hyperholomorphic line bundle $\mathcal L$,
we obtain from the Atiyah--Ward construction that $\text{ev}^*\mathcal L
\longrightarrow \mathcal S^0$ is equipped with the holomorphic connection $\pi_2^*\nabla^{\rm AW}$.
From Theorem~\ref{AWH}, we have
\[\operatorname{ev}^*\nabla =
\pi_2^*\nabla^{\rm AW}+ \frac{\phi}{\widehat\sigma}.\]
Note that not only $\operatorname{ev}^*\nabla$ is the pullback of a meromorphic connection on the twistor space but also
$\pi_2^*\nabla^{\rm AW}$ together with the meromorphic 1-form $\frac{\phi}{\widehat\sigma}$ are pullbacks from $Z$
as well.
Clearly, the corresponding objects on $Z$ are neither a holomorphic connection nor a meromorphic 1-form but only real analytic objects.

We shall give a local construction of the holomorphic connection $\nabla^{\rm AW}$ on $\mathcal S^0$
which only involves the above listed twistor data: As explained in Section \ref{sectionaboutsections}, we obtain from the twistor data the holomorphic symplectic form
$\Omega_0.$ Locally, the space
$\mathcal S^0$ is identified with $Z_0\times Z_\infty$ by evaluation, and there is the moment map $\mathcal E$
of the circle action on $\mathcal S^0$; see Theorem~\ref{EnergyMomentMap}.
Define the vector bundle automorphism
\[\eta \colon\ T^*(Z_0\times Z_\infty) \longrightarrow
T^*(Z_0\times Z_\infty)\]
that acts on $T^*Z_0$ (resp.\ $T^* Z_\infty$) as multiplication
by $+1$ (resp.\ $-1$). It is evident that there are locally defined holomorphic 1-forms
$\alpha_U$ with $\eta(\alpha_U) = \alpha_U$ and
\begin{equation}\label{connonS}
{\rm d}\alpha_U =
\Omega_0+\i \, {\rm d}I_0{\rm d}\mathcal E.\end{equation}
We may think of the forms $\alpha_U$ as holomorphic connection 1-forms. Note that the corresponding cocycle is the pullback of a cocycle on $Z_0.$

By Corollary \ref{om0curaw}, the connection defined by \eqref{connonS} is flat on $\mathcal
S_z$ for every $z \in Z$, and therefore it is locally trivial. Hence, we recover the Atiyah--Ward
connection at least locally. By adding the 1-form $\frac{\phi}{\widehat\sigma}$, which is
constructed from the circle action and from the twisted symplectic form via Lemma
\ref{Lemma:holoconn}, we obtain (the pullback of) the meromorphic connection on $Z$.

\subsubsection[tau-sesquilinear forms and a generalized Chern connection]{$\boldsymbol{\tau}$-sesquilinear forms and a generalized Chern connection}

It seems appropriate to put the above local construction into a global context.
We need an additional structure on $\mathcal L \longrightarrow \mathcal S$, which can
be regarded as the complexification of the hermitian metric $h_M$ on $L_M$.

\begin{Definition}
Let $\mathcal L \longrightarrow \mathcal S^0$ be a holomorphic line bundle, and
let $\tau \colon \mathcal S^0 \longrightarrow \mathcal S^0$ be an anti-holomorphic involution.
A non-degenerate pairing
\[ \langle\cdot, \cdot\rangle \colon \ \mathcal L\times \tau^*\overline{\mathcal L} \longrightarrow \C\]
is called a \emph{holomorphic $\tau$-sesquilinear form on $\mathcal L$} if for all
local holomorphic sections $v, w \in \Gamma(\mathcal U, \mathcal L)$ defined on some
$\tau$-invariant open subset
$\mathcal U \subset \mathcal S$, the function $\langle\langle v, w\rangle\rangle$ on
$\mathcal U$ defined by~${\langle\langle v, w\rangle\rangle (s) := \langle v(s), \tau^*w(s)\rangle}$
is holomorphic and it satisfies the identity
\smash{$\tau^*\overline{\langle\langle v, v \rangle\rangle} = \langle\langle v, v \rangle\rangle$}.
\end{Definition}

Let $M' = \bigl({\mathcal S}^0\bigr)^\tau \subset {\mathcal S}^0$ be the set of real points.
It follows that $\langle v_s, v_s\rangle \in {\mathbb R}$ for all $s \in M'$ and~$v_s \in {\mathcal L}_s$.
Note that for functions $f, g \colon \mathcal U \longrightarrow \C$, we have
\[\langle\langle f v, g v\rangle\rangle = f (\tau^*\overline{g}) \langle\langle v, v\rangle\rangle.\]
Thus, up to a sign, a holomorphic $\tau$-sesquilinear form is a complexification of a
hermitian metric over the locus of real points.

Assume that there is a holomorphic involution
$\eta \in H^0\bigl(\mathcal S^0, \operatorname{End}\bigl(T^*\mathcal S^0\bigr)\bigr)$
such that
\begin{equation}\label{Eq: etareal}
\eta(\tau^*\overline\alpha) = -\tau^*\overline{\eta(\alpha)}\qquad \text{for any}\ \alpha \in \Omega^{1,0}\bigl(\mathcal S^0\bigr) .
\end{equation}
Further, assume that the holomorphic line bundle $\mathcal L
 \longrightarrow \mathcal S^0$ is
given by a cocycle
of the form
\begin{equation}\label{taucocylce}\{f_{i,j}\colon \mathcal U_{i,j} \longrightarrow \C^*\},\end{equation}where the derivatives ${\rm d} f_{i,j}$ are in
the $+1$ eigenspaces of $\eta$, that is,
\[{\rm d}f_{i,j} = \frac{1}{2}( {\rm d}f_{i,j}+\eta ({\rm d}f_{i,j}) ).\]
Under these assumptions, we have the following.

\begin{Lemma}\label{Lemma: ccc}
There is a unique holomorphic {\em Chern}-connection
$\nabla$ on $\mathcal L \longrightarrow \mathcal S^0$ with
\begin{equation}\label{Eq: CCC1}
d (\langle\langle \sigma, \sigma\rangle\rangle) =
\langle\langle\nabla \sigma, \sigma\rangle \rangle+ \langle\langle \sigma, \nabla \sigma\rangle\rangle
\end{equation}
for any holomorphic section $\sigma \in \Gamma(\mathcal U, \mathcal L)$ and any
$\tau$-invariant open subset $\mathcal U \subset \mathcal S^0$.
\end{Lemma}

\begin{proof}
Take a collection of local trivializations by nowhere-vanishing holomorphic sections~$
\sigma_i \in \Gamma(\mathcal U_i, \mathcal L)$ satisfying the
condition that locally we have $\sigma_i
 = \operatorname{ev}_0^*v_i$, where each open subset $\mathcal U_i \subset {\mathcal S}^0$
is $\tau$-invariant, such that the
transition functions satisfy the condition~${{\rm d}f_{i,j} = \eta({\rm d}f_{i,j})}$, where~$f_{i,j} = \frac{\sigma_i}{\sigma_j}$. Then, writing
$
\nabla\sigma_i = A_i\sigma_i$,
we have $\eta(A_i) = A_i$. The condition in \eqref{Eq: CCC1} then implies that
\[
{\rm d}\log\langle\langle \sigma_i, \sigma_i\rangle\rangle =
 \bigl(A_i+ \tau^*\overline{A_i}\bigr)\langle\langle \sigma_i, \sigma_i\rangle\rangle.
\]
By \eqref{Eq: etareal}, we have $\eta\bigl(\tau^*\overline{A_i}\bigr) = -\tau^*\overline{A_i}$ and thus
\[
A_i = \frac{1}{2}({\rm d}\log\langle \langle \sigma_i, \sigma_i\rangle\rangle+
\eta {\rm d}\log \langle\langle \sigma_i, \sigma_i\rangle\rangle).
\]
This shows the uniqueness of the connection.

For the existence of connection, set
\[\nabla \sigma_i = \frac{1}{2}({\rm d}\log\langle \langle \sigma_i, \sigma_i\rangle\rangle+
\eta {\rm d}\log \langle\langle \sigma_i, \sigma_i\rangle\rangle)\otimes \sigma_i.\]
Arguing analogously just as for the usual Chern connection, we see that this defines a connection
with the desired properties.
\end{proof}

The main example we have in mind arises as follows. Consider a hyperK\"ahler manifold~$M$ equipped with a rotating $S^1$-action and associated hyperholomorphic hermitian line bundle~${L_M \longrightarrow M}$ as described in Section~\ref{ss:rotating}. Let $\mathcal S^0_M$ be a connected component of $\mathcal S^0$ that meets~$M$ (for the definition of $\mathcal S^0$, see the beginning of Section~\ref{ss:AW}). Assume that there exists a~global $\tau$-sesquilinear form on the bundle
$\mathcal L \longrightarrow S^0_M$ obtained from the bundle $L_Z$ via the Atiyah--Ward transform as in Section~\ref{ss:AW}.
Consider the decomposition in \eqref{eq:splitlambda} for~${x = 0 \in \C P^1}$,
and the involution which is $-$Id on $T_{\text{ev}_0}\mathcal S^0_M$ and $+$Id on
$T_{\text{ev}_\infty}\mathcal S^0_M$.
Denote its dual endomorphism by~$\eta$. It satisfies \eqref{Eq: etareal}.
Now, $\mathcal L$ agrees with the hyperholomorphic line bundle~$L_M$ and $T_{\text{ev}_\infty}\mathcal S^0_M$ agrees with $T^{1,0}M$ (see Remark \ref{rem:T10M}) on $M\subset \mathcal S^0_M$. It follows that $\mathcal L$ admits a cocycle as in~\eqref{taucocylce}.
Applying Lemma \ref{Lemma: ccc} we obtain a natural connection, which
can be shown to be the Atiyah--Ward connection.
It would be interesting to find natural expressions for the holomorphic
$\tau$-sesquilinear form in concrete examples such as the Deligne--Hitchin moduli spaces;
see also Section \ref{ss:thelinebundleonMDH} below and the following Remark \ref{re6}.

\begin{Remark}\label{re6}
The existence of a $\tau$-sesquilinear form on $\mathcal L \longrightarrow \mathcal S^0$
follows almost
automatically from the existence of an Atiyah--Ward type connection $\nabla$ on
$\mathcal L \longrightarrow \mathcal S^0$: assume that there exists a hermitian metric $h$ on
$\mathcal L_{\mid
M} \longrightarrow M$ over the real points $M \subset \mathcal S^0.$ Consider a local
holomorphic frame
$\sigma \in \Gamma(\mathcal U, \mathcal L)$ on a simply-connected $\tau$-invariant open subset
$\mathcal U \subset \mathcal S^0$ such that $M \cap \mathcal U \not= \emptyset$. Write $\nabla
\sigma = \alpha\otimes\sigma.$ Then
$\alpha+\tau^*\overline{\alpha}$ is closed; its exterior derivative is a $(2,0)$-form
which vanishes on $M$, and hence vanishes on its complexification $\mathcal S^0$. Integrating this closed
form on the simply connected set $\mathcal U$ produces the $\tau$-sesquilinear form
via
\[\langle\langle(\mu_1\sigma(p), \mu_2\sigma(\tau(p))\rangle\rangle =
\mu_1\overline\mu_2 h(\sigma(p), \sigma(p))\exp\left(\int_s^p\alpha+\tau^*\overline{\alpha}\right)\]
for all $\mu_1, \mu_2 \in \C$ and $p \in \mathcal U$.
\end{Remark}

\section[Space of holomorphic sections of the Deligne--Hitchin moduli
space]{Space of holomorphic sections of the Deligne--Hitchin\\ moduli
space} \label{s:DeligneHitchin}

In this section, we illustrate the general theory described in Sections~\ref{sectionaboutsections} and~\ref{sectionhyperholo} for the
Deligne--Hitchin twistor space of the moduli space of Higgs bundles on a compact Riemann
surface.

\subsection{Hitchin's self-duality equations}\label{ss: MSD}

Let $\Sigma$ be a compact Riemann surface of genus at least two; denote its holomorphic
cotangent bundle by $K_\Sigma$. Consider
a smooth complex rank $n$ vector bundle $E$ of degree zero with structure group in ${\rm SU}(n)$.

Denote by $\mathcal A(E)$ the space of ${\rm SU}(n)$-connections on $E$. Sending any \smash{$\nabla
 \in \mathcal A(E)$} to its $(0,1)$-part \smash{$\overline{\partial}^\nabla$} we may
identify $\mathcal A$ with the space of holomorphic structures on $E$ inducing the
trivial holomorphic structure $\overline{\partial}$ on the determinant bundle
$\det E = \bigwedge^n E$ of $E$, that is, $\smash{\bigl(\det E, \overline{\partial}^\nabla\bigr)}
 = \bigl(\mathcal O_\Sigma, \overline{\partial}\bigr)$. Thus, $\mathcal A$ is an affine space modelled
on $\Omega^{0,1}(\Sigma, \sln(E))$, where $\sln(E)$ is the bundle of trace-free endomorphisms of $E$. The product
\[
T^*\mathcal A = \mathcal A\times \Omega^{1,0}(\Sigma, \sln(E))
\]
can be thought of as the cotangent bundle of $\mathcal A$. We will denote its elements by
$\bigl(\overline{\partial}^\nabla, \Phi\bigr)$, and we will call $\Phi$ the \emph{Higgs field} of the pair
$\bigl(\overline{\partial}, \Phi\bigr)$.
Formally, $T^*\mathcal A$ carries a flat hyperK\"ahler structure that can be described as
follows. The Riemannian metric is given by the $L^2$-inner product on $\Omega^{1}(\Sigma, \sln(E))$.
The almost complex structures $I, J, K = IJ$ act on a tangent vector~$\smash{(\alpha, \phi) \in T_{(\overline{\partial}^\nabla,\Phi)}T^*\mathcal A = \Omega^{0,1}(\Sigma,\sln(E)) \oplus \Omega^{1,0}(\Sigma,\sln(E))}$ as
$
I(\alpha, \phi) = (\i\alpha, \i\phi)$, $
J(\alpha, \phi) = (-\phi^*, \alpha^*)$.
The K\"ahler forms are given by
\begin{gather}
\omega_I((\alpha, \phi), (\beta, \psi)) = \int_{\Sigma}\operatorname{tr}(\alpha^*\wedge\beta-\beta^*\wedge\alpha + \phi\wedge\psi^*-\psi\wedge\phi^*),\nonumber\\
(\omega_J + \i \omega_K)((\alpha, \phi), (\beta, \psi)) = 2\i\int_\Sigma \operatorname{tr}(\psi\wedge\alpha-\phi\wedge\beta).\label{eq:Msdkahlerforms}
\end{gather}
The group $\mathcal G := \Gamma(\Sigma, {\rm SU}(E))$ of unitary gauge transformations
of $E$ acts (on the right) on~$T^*\mathcal A$~as
\begin{gather}\label{eq:gauge}
\bigl(\overline{\partial}^\nabla, \Phi\bigr).g = \bigl(g^{-1} \circ \overline{\partial}^\nabla \circ g,
 g^{-1}\Phi g\bigr).
\end{gather}
This action preserves the flat hyperK\"ahler structure, and, formally, the vanishing condition for the associated hyperK\"ahler moment map yields \emph{Hitchin's self-duality equations}
\begin{align}\label{eq:SD}
F^\nabla + [\Phi\wedge\Phi^*] = 0,\qquad
\overline{\partial}^\nabla\Phi = 0.
\end{align}
The second equation implies that $\Phi \in H^0(\Sigma, \sln(E)\otimes K_\Sigma)$ with respect to the holomorphic structure \smash{$\overline{\partial}^\nabla$}. Note that these equations imply that the connection $\nabla + \Phi +\Phi^*$ is flat.

Let $\mathcal H \subset T^*\mathcal A$ be the set of solutions of \eqref{eq:SD}. The moduli space of solutions to the self-duality equations is the quotient
\begin{equation}\label{eq:MSD}
\mathcal{M}_{\mathrm{SD}} := \mathcal{M}_{\mathrm{SD}}(\Sigma, \operatorname{SL}(n,\C)) := \mathcal H/\mathcal G .
\end{equation}
Formally, it is the hyperK\"ahler quotient of $T^*\mathcal A$ by the action of $\mathcal G$.

A solution $\bigl(\overline{\partial}^\nabla, \Phi\bigr)$ to the self-duality equations \eqref{eq:SD}
is called \emph{irreducible} if its stabilizer in~$\mathcal G$ is the center of
${\rm SU}(n)$. We write $\mathcal H^{\rm irr} \subset \mathcal{H}$ for the set of irreducible
solutions. It is known that~$
\mathcal{M}_{\mathrm{SD}}^{\rm irr} = \mathcal H^{\rm irr}/\mathcal G
$
is the smooth locus of $\mathcal{M}_{\mathrm{SD}}$ in \eqref{eq:MSD}. It is equipped with a hyperK\"ahler metric induced
by the flat hyperK\"ahler structure on $T^*\mathcal A$ described above. Moreover,
there is a~rotating circle action on $\mathcal{M}_{\mathrm{SD}}$ given by
\begin{equation}\label{eq:S1Msd}
\zeta.\bigl(\overline{\partial}^\nabla, \Phi\bigr) = \bigl(\overline{\partial}^\nabla, \zeta\Phi\bigr),
\end{equation}
$\zeta \in S^1$. It is Hamiltonian with respect to $\omega_I$, and the map
\begin{equation}\label{eq:muI}
\mu_I \colon\ \mathcal{M}_{\mathrm{SD}} \longrightarrow
\i \R, \qquad \mu_I\bigl(\overline{\partial}^\nabla, \Phi\bigr) = - \int_\Sigma\operatorname{tr}(\Phi\wedge\Phi^*)
\end{equation}
restricts to a natural moment map on $\mathcal{M}_{\mathrm{SD}}^{\rm irr}$.

The complex manifold \smash{$\bigl(\mathcal{M}_{\mathrm{SD}}^{\rm irr}, I\bigr)$} can be described as follows. An \emph{$\operatorname{SL}(n,\C)$-Higgs
bundle} is a~pair consisting of a holomorphic vector bundle $\bigl(E, \overline{\partial}_E\bigr)$, such
that $\det\bigl(E, \overline{\partial}_E\bigr) = \mathcal O_\Sigma$, together with a Higgs field
$\Phi \in H^0(\Sigma, \sln(E)\otimes K_\Sigma)$. The group $\mathcal G_{\C} =
\Gamma(\Sigma, \operatorname{SL}(E))$ acts on the set of Higgs bundles by the same rule as
in \eqref{eq:gauge}. A Higgs bundle $\bigl(\overline{\partial}_E, \Phi\bigr)$ is called \emph{stable}
if every $\Phi$-invariant holomorphic subbundle $E'$ has negative degree, and it
is called \emph{polystable} if it is a direct sum of stable Higgs bundles
of degree zero. We denote by
\[\mathcal{M}_{\mathrm{Higgs}}^{\rm st} := \mathcal{M}_{\mathrm{Higgs}}^{\rm st}(\Sigma,\operatorname{SL}(n,\C))\qquad \text{and}\qquad \mathcal{M}_{\mathrm{Higgs}}^{\rm ps} :=
\mathcal{M}_{\mathrm{Higgs}}^{\rm ps}(\Sigma,\operatorname{SL}(n,\C))\]
the moduli spaces of stable and polystable Higgs bundles
respectively. It is known that the map~\smash{$\bigl(\overline{\partial}^\nabla, \Phi\bigr) \longmapsto
 ((E, \overline{\partial}^\nabla),\Phi)$} induces a biholomorphism
~\smash{$
\bigl(\mathcal{M}_{\mathrm{SD}}^{\rm irr}, I\bigr) = \mathcal{M}_{\mathrm{Higgs}}^{\rm st}
$}
that extends to a~homeomorphism~\smash{$
\mathcal{M}_{\mathrm{SD}} = \mathcal{M}_{\mathrm{Higgs}}^{\rm ps}$}.
The circle action described in \eqref{eq:S1Msd} extends to a holomorphic $\C^*$-action on $\mathcal{M}_{\mathrm{Higgs}}$
\begin{equation}\label{eq:C*Mh}
\zeta\cdot\bigl(\overline{\partial}, \Phi\bigr) = \bigl(\overline{\partial}, \zeta\Phi\bigr).
\end{equation}
While the complex structure $I$ is induced from the complex structure on $\Sigma$, the
complex structure~$J$ has a more topological origin. Consider the space $\mathcal A_\C$ of $\operatorname{SL}(n,\C)$ connections on $E$. The group~$\mathcal G_\C$ acts on $\mathcal A_\C$ as
$
\nabla.g = g^{-1}\circ\nabla\circ g$.
A connection $\nabla \in \mathcal A_\C$ is called \emph{irreducible} if its stabilizer in
$\mathcal G_\C$ is the center of $\operatorname{SL}(n,\C)$, and $\nabla$ is called \emph{reductive} if it is
isomorphic to
a direct sum of irreducible connections, in other words, if any $\nabla$-invariant subbundle $E' \subset E$
admits a~$\nabla$-invariant complement. Let
$\mathcal{M}_{\mathrm{dR}} := \mathcal{M}_{\mathrm{dR}}(\Sigma, \operatorname{SL}(n,\C))$
be the moduli space of reductive flat $\operatorname{SL}(n,\C)$-connections on $E$. Its smooth locus is given
by~\smash{$\mathcal{M}_{\mathrm{dR}}^{\rm irr}$}, the moduli space of irreducible flat $\operatorname{SL}(n,\C)$-connections. Then the
map $\smash{\bigl(\overline{\partial}^\nabla, \Phi\bigr) }\longmapsto \nabla+\Phi+\Phi^*$ is a biholomorphism~\smash{$
\bigl(\mathcal{M}_{\mathrm{SD}}^{\rm irr}, J\bigr) \cong \mathcal{M}_{\mathrm{dR}}^{\rm irr}
$}
that extends to a homeomorphism~\smash{$
\mathcal{M}_{\mathrm{SD}} \cong \mathcal{M}_{\mathrm{dR}}$}.
The Riemann--Hilbert correspondence gives a biholomorphism
\[
\mathcal{M}_{\mathrm{dR}} \cong {\rm Hom}(\pi_1(\Sigma), \operatorname{SL}(n,\C))^{\mathrm{red}}/\operatorname{SL}(n,\C) =:
 \mathcal{M}_{\mathrm{B}}(\Sigma,\operatorname{SL}(n,\C)) =: \mathcal{M}_{\mathrm{B}}.
\]
The space $\mathcal{M}_{\mathrm{B}}$ is known as the Betti moduli space.

\subsection{The Deligne--Hitchin moduli space}

We now recall Deligne's construction of the twistor space $Z(\mathcal{M}_{\mathrm{SD}})$ via $\lambda$-connections;
see \cite{Si-Hodge}.

Take any $\lambda \in \mathbb C$. A holomorphic $\operatorname{SL}(n,\C)$ $\lambda$-connection on the
${\rm SU}(n)$-vector bundle $E \longrightarrow
\Sigma$ is a pair $\bigl(\overline{\partial}, D\bigr)$, where $\overline{\partial}$ is a~holomorphic
structure on $E$ and $D\colon \Gamma(\Sigma,E) \longrightarrow \Omega^{1,0}(\Sigma,E)$ is a~differential
operator such that
\begin{itemize}\itemsep=0pt
\item for any $f\in \C^\infty(\Sigma)$ we have $D(fs) = \lambda s\otimes \partial f + fDs$,

\item $D$ is holomorphic, that is, $\overline{\partial}\circ D + D\circ\overline{\partial} = 0$, and

\item as a holomorphic vector bundle we have $\bigl(\bigwedge^n E, \overline{\partial}\bigr) =
\bigl({\mathcal O}_\Sigma, \overline{\partial}_\Sigma\bigr)$, and the holomorphic differential
operator on $\Lambda^n E$ induced by $D$ coincides with $\lambda\partial_\Sigma$.
\end{itemize}
Therefore, an $\operatorname{SL}(n,\C)$ $0$-connection is an $\operatorname{SL}(n,\C)$-Higgs bundle, and an
$\operatorname{SL}(n,\C)$ $1$-connection is a usual holomorphic $\operatorname{SL}(n,\C)$-connection. More
generally, the operator $\overline{\partial}_E + \lambda^{-1}D(\lambda)$ is a~holomorphic
connection for every $\lambda \neq 0$.

The group $\mathcal G_\C$ of complex gauge transformations of $E$ acts on the set of holomorphic $\lambda$-connections in the usual way
\[
\bigl(\overline{\partial}, D\bigr)\cdot g =
\bigl(g^{-1}\circ \overline{\partial}\circ g , g^{-1} \circ D\circ g\bigr).
\]
We again call $\bigl(\overline{\partial}, D\bigr)$ stable if any $D$-invariant holomorphic subbundle
$E' \subset E$ has negative degree, and \smash{$\bigl(\overline{\partial}, D\bigr)$} is
called polystable if it is isomorphic to a direct sum of stable $\lambda$-connections
of degree zero. We note that if $\lambda \not= 0$, then the degree of any $\lambda$-connection
is automatically zero. Equivalently, $\bigl(\overline{\partial}, D\bigr)$ is stable, with $\lambda \not= 0$, if it is
irreducible, that is, its stabilizer in $\mathcal G_\C$ is given by the center of $\operatorname{SL}(n, \C)$.
The Hodge-moduli space $\mathcal{M}_{\mathrm{Hod}} := \mathcal{M}_{\mathrm{Hod}}(\Sigma, \operatorname{SL}(n,\C))$ is the moduli space of
polystable holomorphic $\lambda$-connections, where $\lambda$ varies over $\C$
\begin{align*}
\mathcal{M}_{\mathrm{Hod}} :={}& \mathcal{M}_{\mathrm{Hod}}(\Sigma, \operatorname{SL}(n, \C)) \\={}&
\bigl\{\bigl(\lambda, \overline{\partial}, D\bigr) \mid \lambda \in \C,
\text{$\bigl(\overline{\partial}, D\bigr)$ polystable holomorphic $\lambda$-connection}\}/\mathcal G_\C.
\end{align*}
It has a natural holomorphic projection to $\C$ given by
\begin{equation}\label{evpn}
\varpi\colon\ \mathcal{M}_{\mathrm{Hod}} \longrightarrow \C,\qquad \big[\lambda , {\overline{\partial}}, D \big] \longmapsto \lambda.
\end{equation}
The smooth locus $\mathcal{M}_{\mathrm{Hod}}^{\rm irr}$ of the Hodge moduli space $\mathcal{M}_{\mathrm{Hod}}$
coincides with the locus of irreducible $\lambda$-connections, when $\lambda \not= 0$.

The $\C^*$-action in \eqref{eq:C*Mh} extends to a natural $\C^*$-action on $\mathcal{M}_{\mathrm{Hod}}$ covering the standard $\C^*$-action on $\C$ by
$
\zeta\cdot\bigl(\lambda, \overline{\partial}, D\bigr) = \bigl(\zeta\lambda, \overline{\partial}, \zeta D\bigr)$.

The map $\bigl(\lambda, \overline{\partial}, D\bigr) \longmapsto \bigl(\overline{\partial} + \lambda^{-1}D,
 \lambda\bigr)$ induces a biholomorphism
\[
\varpi^{-1}(\C^*) \cong \mathcal{M}_{\mathrm{dR}} \times\C^* \cong \mathcal{M}_{\mathrm{B}} \times\C^*,
\]
where $\varpi$ is the map in \eqref{evpn}.

The Deligne--Hitchin moduli space $\mathcal{M}_{\mathrm{DH}} := \mathcal{M}_{\mathrm{DH}}(\Sigma, \operatorname{SL}(n,\C))$ is obtained by gluing $\mathcal{M}_{\mathrm{Hod}}$
and $\overline{\mathcal{M}_{\mathrm{Hod}}} \cong \mathcal{M}_{\mathrm{Hod}}\bigl(\overline{\Sigma}, \operatorname{SL}(n,\C)\bigr)$ over $\C^*$ via the
Riemann--Hilbert correspondence
\begin{align*}
\mathcal{M}_{\mathrm{DH}} := & \mathcal{M}_{\mathrm{DH}}(\Sigma, \operatorname{SL}(n,\C))
 = \bigl(\mathcal{M}_{\mathrm{Hod}}(\Sigma,\operatorname{SL}(n,\C)) \dot\cup
\mathcal{M}_{\mathrm{Hod}}\bigl(\overline{\Sigma}, \operatorname{SL}(n,\C)\bigr)\bigr)/{\sim},
\end{align*}
where
$
\big[\bigl(\lambda, \overline{\partial}, D\bigr)\big] \sim \big[\bigl(\lambda, \lambda^{-1} D,
\lambda^{-1}\overline{\partial}\bigr)\big]
$
for any $\big[\bigl(\lambda, \overline{\partial}, D\bigr)\big] \in \mathcal{M}_{\mathrm{Hod}}$ with $\lambda \neq 0$.

The projections from the respective Hodge moduli spaces to $\C$ glue to give a holomorphic
projection
$\varpi\colon \mathcal{M}_{\mathrm{DH}} \longrightarrow \mathbb CP^1 $.
The smooth locus of $\mathcal{M}_{\mathrm{DH}}$ coincides with the locus
$\mathcal{M}_{\mathrm{DH}}^{\rm irr}$ of irreducible $\lambda$-connections, when $\lambda \not= \{0, \infty\}$, which in turn coincides with the twistor
space~$Z\bigl(\mathcal{M}_{\mathrm{SD}}^{\rm irr}\bigr)$ of the hyperK\"ahler manifold $\mathcal{M}_{\mathrm{SD}}^{\rm irr}$.

The space $\mathcal{M}_{\mathrm{B}}\times\C^*$ admits an anti-holomorphic involution $\tau_{\mathcal{M}_{\mathrm{B}}}$ covering the
antipodal involution \smash{$\lambda \longmapsto -(\overline{\lambda})^{-1}$} of $\C^*$ which is
constructed as follows. For any $\rho \in {\rm Hom}(\pi_1(\Sigma), \operatorname{SL}(n,\C))$, define
\[
\tau_{\mathcal{M}_{\mathrm{B}}}(\rho, \lambda) = \bigl(\overline{\rho^{-1}}^T, -(\overline{\lambda})^{-1}\bigr).
\]
This produces the following antiholomorphic involution $\tau_{\mathcal{M}_{\mathrm{DH}}}$ of $\mathcal{M}_{\mathrm{DH}}$:
\begin{gather*}
\tau_{\mathcal{M}_{\mathrm{DH}}} \colon\ \mathcal{M}_{\mathrm{DH}} \longrightarrow \mathcal{M}_{\mathrm{DH}},\\
 \big[\bigl(\lambda, \overline{\partial}, D\bigr)\big] \longmapsto \big[\bigl(-(\overline{\lambda})^{-1}, (\overline{\lambda})^{-1}\overline{D}^*,
 -(\overline{\lambda})^{-1}\overline{\overline\partial}^*\bigr)\big]
 = \big[\bigl((-\overline\lambda), \overline{\overline\partial}^*, -\overline{D}^*\bigr)\big].
\end{gather*}
This involution $\tau_{\mathcal{M}_{\mathrm{DH}}}$ is compatible with the $\C^*$-action in the following sense.
For $\zeta \in \C^*$,
\[
\tau_{\mathcal{M}_{\mathrm{DH}}}\bigl(\zeta.\bigl(\lambda, \overline\partial, D\bigr)\bigr) =
\overline \zeta^{-1}.\tau_{\mathcal{M}_{\mathrm{DH}}}\bigl(\lambda \overline\partial, D\bigr).
\]

The $\C^*$-action on $\mathcal{M}_{\mathrm{Hod}}$ extends to a $\C^*$-action on $\mathcal{M}_{\mathrm{DH}}$. For any $\zeta \in \C^*$,
we have
\[
\zeta \cdot \big[\lambda, {\overline{\partial}}, D(\lambda) \big] = \big[\zeta \lambda, {\overline{\partial}}, \zeta D(\lambda)\big].
\]
Clearly, this action covers the natural $\C^*$-action on $\mathbb CP^1$.
Therefore, the $\C^*$-fixed points of $\mathcal{M}_{\mathrm{DH}}$ are given by
\[
\mathcal{M}_{\mathrm{DH}}^{\C^*} \cong \mathcal{M}_{\mathrm{Higgs}}(\Sigma, \operatorname{SL}(n,\C))^{\C^*} \amalg \mathcal{M}_{\mathrm{Higgs}}\bigl(\overline{\Sigma}, \operatorname{SL}(n,\C)\bigr)^{\C^*},
\]
that is, by the locus of complex variations of Hodge structures on $\Sigma$
and $\overline{\Sigma}$ (see \cite{SimpsonIHES1992}).

The twisted relative symplectic form on \smash{$\mathcal{M}_{\mathrm{DH}}^{\rm irr}$} can be described as follows. Let
\smash{$V_i = \bigl(\dot{\overline{\partial}_i}, \dot D_i\bigr)$}, $i = 1, 2$, be a pair of tangent
vectors to the fiber $\varpi^{-1}(\lambda)$. Then
\begin{equation}\label{eq:omegalambdaMDH}
\omega_\lambda(V_1, V_2) = 2\i\int_\Sigma \operatorname{tr}\bigl(\dot D_2\wedge \dot{\overline{\partial}}_1-\dot D_1\wedge\dot{\overline{\partial}}_2\bigr).
\end{equation}
Note that at $\lambda = 0$ this exactly resembles $\omega_J +\i \omega_K$ as defined in \eqref{eq:Msdkahlerforms}.

\subsection[The line bundle on M\_DH\^{}irr]{The line bundle on $\boldsymbol{\mathcal{M}_{\mathrm{DH}}^{\rm irr}}$}\label{ss:thelinebundleonMDH}

Let us describe the holomorphic line bundle $L_Z$ for $Z = \mathcal{M}_{\mathrm{DH}}^{\rm irr}$ along each fiber of
$\varpi$. First consider the moduli space of holomorphic $\operatorname{SL}(n,\C)$-connections $\mathcal{M}_{\mathrm{dR}}$. Let $\mathcal{M}_{\mathrm{dR}}'$ be the Zariski open subset of $\mathcal{M}_{\mathrm{dR}}^{\rm irr}$ such that the underlying holomorphic vector
bundle is stable. So we have a holomorphic map \[ f\colon \ \mathcal{M}_{\mathrm{dR}}' \longrightarrow
\mathcal{N} , \] where $\mathcal{N}$ is the moduli space of stable
$\operatorname{SL}(n,\C)$-bundles, that sends any pair $(\bigl(E, {\overline{\partial}}_E\bigr), \nabla^E)$ to~$\bigl(E, {\overline{\partial}}_E\bigr)$.

The pullback map $f^*\colon \mathrm{Pic}(\mathcal{N}) \longrightarrow \mathrm{Pic}(\mathcal{M}_{\mathrm{dR}}')$ is an isomorphism~\cite{BR}. On the other hand, $\mathrm{Pic}(\mathcal{M}_{\mathrm{dR}}') = {\mathbb Z}$, and holomorphic line bundles
on $\mathcal{M}_{\mathrm{dR}}'$ are uniquely determined by their first class~\cite{DN}.
Also the restriction map $\mathrm{Pic}\bigl(\mathcal{M}_{\mathrm{dR}}^{\rm irr}\bigr) \longrightarrow \mathrm{Pic}(\mathcal{M}_{\mathrm{dR}}')$ is an isomorphism,
because the codimension of the complement of $\mathcal{M}_{\mathrm{dR}}'$ inside $\mathcal{M}_{\mathrm{dR}}^{\rm irr}$ is at least two
(see \cite[p.~6, Lemma~3.1]{BMunoz2}, \cite[p.~202--203]{BMunoz1}).
Therefore, we have $\mathrm{Pic}\bigl(\mathcal{M}_{\mathrm{dR}}^{\rm irr}\bigr) = \mathrm{Pic}(\mathcal{N}) \cong \mathbb{Z}$, and
the holomorphic line bundles on $\mathcal{M}_{\mathrm{dR}}^{\rm irr}$ are uniquely determined by their first Chern class.

Let $L_1$ denote the restriction of $L_Z$ to $Z_1 = \mathcal{M}_{\mathrm{dR}}^{\rm irr}$. On $L_1$, the meromorphic connection
$\nabla$ on $L_Z$ induces a holomorphic connection and its curvature on $Z_1$ is the holomorphic symplectic
form $2\omega_J+2\i\omega_I$, which is cohomologous to $2\i\omega_I$, since $\omega_J$ is exact. So by
Chern--Weil theory, we have $c_1(L_1) = [\omega_I]$, where $[\omega_I]$ is the cohomology class of $\omega_I$.

It follows that $L_1$ is the holomorphic line bundle on $Z_1$ determined by $[\omega_I]$. To
calculate~$[\omega_I]$, note that the above projection $f$ has a $C^\infty$-section that sends
any stable vector bundle $\bigl(E, {\overline{\partial}}_E\bigr)$ to the unique unitary flat connection on $E$ \cite{NS}.
The restriction
of $\omega_I$ to the image of that section coincides with the standard K\"ahler form on $\mathcal{N}$.
On the other hand, the first Chern class of the determinant line bundle $\xi$ on $\mathcal{N}$ is the
K\"ahler form \cite{Qu}. In fact, the curvature for the Quillen metric on $\xi$ is the K\"ahler form.

We conclude that $L_1$ is holomorphically isomorphic to the determinant bundle on $Z_1$, which
coincides with the pullback of the determinant line bundle on $\mathcal{N}$ because the determinant
line bundle is functorial. Similarly, the line bundles on the moduli space of stable Higgs
bundles, namely $Z_0 = \mathcal{M}_{\mathrm{Higgs}}^{\rm st}$, are uniquely determined by their first Chern class. Let $L_0$
denote the restriction of $L_Z$ to $Z_0$. Since $L_1$ is isomorphic to the determinant line
bundle, and since the first Chern class of the family of line bundles $L_t \longrightarrow
Z_t$, $t \in \C
P^1$, is independent of $t$, it follows that~$L_0$ is isomorphic to the determinant line
bundle. For a description of $L_Z$ on all of $Z = \mathcal{M}_{\mathrm{DH}}^{\rm irr}$, see \cite[Section~3.7]{Hitchin-HKQK}.

\subsection{Irreducible and admissible sections}

In this subsection, we recall some concepts and definitions from \cite{BHR} and \cite{BeHeRo} on
sections of the twistor projection $\varpi\colon \mathcal{M}_{\mathrm{DH}} \longrightarrow
\C P^1$. We write $\mathcal S_{\mathcal{M}_{\mathrm{DH}}}$ for
the space of holomorphic sections. Since $\mathcal{M}_{\mathrm{DH}}$ is a complex space, a similar argument as in
Proposition \ref{prop:SpaceOfSections} shows that $\mathcal{S}_{\mathcal{M}_{\mathrm{DH}}}$ is a complex space as
well. It is equipped with an antiholomorphic involution $\tau$ defined in \eqref{tau}, with
$\tau_{\mathcal{M}_{\mathrm{DH}}}$ playing the role of $\tau_Z$.

\begin{Definition} \label{Def: irred_section}
A holomorphic section $s \in \mathcal S_{\mathcal{M}_{\mathrm{DH}}}$ is \emph{irreducible} if the image of
$s$ is contained in $\mathcal{M}_{\mathrm{DH}}^{\rm irr}$.
These are precisely the sections of the twistor space of $\mathcal{M}_{\mathrm{SD}}^{\rm irr}$.
\end{Definition}

If $s\colon \C P^1 \longrightarrow \mathcal{M}_{\mathrm{DH}}$ is an irreducible section, then by
 \cite[Lemma 2.2]{BHR} (see also \cite[Remark~1.11]{BeHeRo}) it admits a holomorphic lift
\begin{equation} \label{Eq: LiftPowerseries}
\widehat s(\lambda) = \bigl(\lambda, \overline{\partial}(\lambda), D(\lambda)\bigr)
 = \left(\lambda, \overline\partial + \lambda\Psi+ \sum_{j=2}^\infty\lambda^j\Psi_j, \lambda\partial + \Phi
+ \sum_{j=2}^\infty\lambda^j\Phi_j\right), \qquad\lambda \in \C,
\end{equation}
to the space of holomorphic $\lambda$-connections of class $C^k$. Here
$\Psi, \Psi_j \in \Omega^{0,1}(\mathfrak{sl}(E))$ and $\Phi, \Phi_j \in \Omega^{1,0}(\mathfrak{sl}(E))$. There is
also a lift $^-\widehat s$ on $\mathbb CP^1\setminus\{0\}$.

The lifts $\widehat{s}$, $^-\widehat{s}$ over $\C$ and $\mathbb CP^1\setminus\{0\}$, respectively, allow us to interpret $s|_{\C^*}$ as a $\C^*$-family of flat connections
\[
^+\nabla^\lambda = \overline\partial(\lambda) + \lambda^{-1}D(\lambda) =
\lambda^{-1}\Phi +\nabla + \cdots ,
\]
and $^-\nabla^\lambda$ defined similarly over $\C P^1\setminus\{0, \infty\}$. Here we write
$\nabla = \overline\partial + \partial$ using the notation of equation
\eqref{Eq: LiftPowerseries}.
There exists a holomorphic $\C^*$-family $g(\lambda)$ of $\mathrm{GL}(n,\C)$-valued
gauge transformations, unique up to multiplication by a holomorphic scalar function,
such that $^+\nabla^\lambda.g(\lambda) = {}^-\nabla^\lambda$ (see \cite{BHR}).

\begin{Definition}\label{admissible}
We call a holomorphic section $s \in \mathcal S_{\mathcal{M}_{\mathrm{DH}}}$ {\it admissible} if it admits a
lift $\widehat s$ on $\mathbb C$ of the form
${\widehat s}(\lambda) = \bigl(\lambda, \overline{\partial}+\lambda \Psi, \lambda\partial+\Phi\bigr)$
for a Dolbeault operator $\overline\partial$ of type $(0,1)$, a
Dolbeault operator $\partial$ of type $(1,0)$, a~$(1,0)$-form~$\Phi$
and a $(0,1)$-form $\Psi$, such that $\bigl(\overline{\partial}, \Phi\bigr)$ and $(\partial, \Psi)$
are semi-stable Higgs pairs on $\Sigma$ and $\overline\Sigma$ respectively.
\end{Definition}

Now suppose that $s \in \mathcal S_{\mathcal{M}_{\mathrm{DH}}}^\tau$ is a real section, that is,
$s = \tau_{\mathcal{M}_{\mathrm{DH}}}\circ s\circ \tau_{\C P^1}$.
If we have a~lift~$\nabla^\lambda$ on $\mathbb C \subset \mathbb CP^1$ of $s$, then
for every $\lambda \in \mathbb C^*$ there is a gauge
transformation $g(\lambda) \in \mathcal G_\C$ such that
\[
\nabla^\lambda.g(\lambda) = \overline{\nabla^{-(\overline{\lambda})^{-1}}}^*.
\]

\begin{Remark}\label{Rem: twistorlines}
If $\bigl(\overline\partial^\nabla, \Phi\bigr)$ is a solution to the self-duality equations, then
the associated twistor line is given by the $\C^*$-family of flat $\operatorname{SL}(n,\C)$-connections
$
\nabla^{\lambda} = \lambda^{-1}\Phi + \nabla + \lambda\Phi^{*_h}$.
By the non-abelian Hodge correspondence, this family gives rise to an equivariant harmonic
map~${f\colon \widetilde\Sigma \longrightarrow \operatorname{SL}(n,\C)/{\rm SU}(n)}$ from the universal cover.
If the solution \smash{$\bigl(\overline{\partial}^\nabla, \Phi\bigr)$} is irreducible, then so is the associated section in $\mathcal S^\tau_{\mathcal{M}_{\mathrm{DH}}}$.
\end{Remark}

In the next example, we describe a large class of irreducible sections for $n = 2$.

\begin{Example}\label{exm1}
Let $\nabla$ be an irreducible flat $\operatorname{SL}(2,\C)$-connection on the rank two vector bundle~${V
\longrightarrow \Sigma}$.
We assume that \smash{$\overline\partial^\nabla$} is not strictly semi-stable, but it is either stable or
it is strictly
unstable. We also assume the analogous condition for the holomorphic structure
$\partial^\nabla$ on~$\overline\Sigma$. Due to \cite{BDH}, there exist irreducible flat $\operatorname{SL}(n,\C)$-connections for which this assumption does not hold. Under the assumption, we obtain a section
$s = s^\nabla$ as follows. Over $\C^*$ consider the (constant) family of flat connections
$^+\nabla^\lambda = \nabla$.
If \smash{$\overline\partial^\nabla$} is stable
\smash{$\widehat s(\lambda) = \bigl(\lambda, \overline\partial^\nabla, \lambda \partial^\nabla\bigr)$}
is the lift of an irreducible section of $\mathcal{M}_{\mathrm{DH}}$ over $\C \subset \C P^1$.
If \smash{$\overline\partial^\nabla$} is unstable, we consider its destabilizing subbundle~${L \subset
V}$ of positive degree. The connection
induces a nilpotent Higgs field
$\Phi$ on the holomorphic vector bundle $L\oplus (V/L) = L\oplus L^*$ via
\smash{$\Phi = \pi^{V/L}\circ\nabla_{\mid L}$}.
This is a special case of \cite{Simpson-iterated} and it can be interpreted from a gauge
theoretic point of view (see also \cite[Section~4]{BeHeRo} for details):
Consider a complementary subbundle $\widetilde{L} \subset V$ of $L$, together with the family of
gauge-transformations
\[g(\lambda) = \begin{pmatrix} 1&0\\0&\lambda \end{pmatrix}\]
on $V = L\oplus \widetilde{L}$.
The gauge equivalent family of flat connections
$^+\nabla^\lambda = \nabla.g(\lambda)$
then corresponds to the lift
\[\widehat s(\lambda) = \bigl( \lambda, \overline\partial^{\nabla.g(\lambda)}, \partial^{\nabla.g(\lambda)}\bigr)\]
which is well-defined on all of $\C$ and $^+\widehat s(0)$ may be identified with the irreducible (stable) Higgs pair
$(L\oplus L^*, \Phi)$. Similarly, if the holomorphic structure $\partial^\nabla$ on $\overline{\Sigma}$ is stable, then
\[^-\widehat{s}(\lambda) = \bigl(\partial^\nabla, \lambda^{-1}\overline{\partial}^\nabla, \lambda^{-1}\bigr)\]
is the lift of an irreducible section over $\C P^1\setminus \{0\}$. If the holomorphic structure $\partial^\nabla$ on $\overline{\Sigma}$ is unstable, then a family $^-g\bigl(\lambda^{-1}\bigr)$ constructed in an analogous manner as $g(\lambda)$ above gives a lift of the form
\[^-\widehat{s}(\lambda) = \bigl(\partial^{\nabla.^-g(\lambda)}, \lambda^{-1}\overline{\partial}^{\nabla.^-g(\lambda)}, \lambda^{-1}\bigr)\]
defined on all of $\C P^1\setminus\{0\}$ such that $^-s(\infty)$ is an irreducible Higgs pair on $\overline{\Sigma}$.
\end{Example}

\section[Energy functional on sections of the Deligne--Hitchin moduli space]{Energy functional on sections of the Deligne--Hitchin\\ moduli space}\label{ss:GeometrySectionDH}

\subsection{The energy as a moment map}

It was proven in \cite[Corollary 3.11]{BeHeRo} that the energy of an irreducible section $s$ with lift $\widehat{s}$ as in~\eqref{Eq: LiftPowerseries} is given by
\begin{equation}\label{eq:EnergyOfLift}
\mathcal{E}(s) = \frac{1}{2\pi \i} \int_{\Sigma} \operatorname{tr}( \Phi \wedge \Psi).
\end{equation}
In particular, this integral is independent of the lift $\widehat{s}$.
The reader should be aware of the different pre-factors in \eqref{eq:EnergyOfLift} and in \eqref{eq:muI}.
In particular, if we think of $\mathcal E$ as the energy of a~harmonic map, it should be real-valued, while we want a~moment map
for the $S^1$-action to be~$\i\R$-valued. Working with the pre-factor~$\frac{1}{2\pi \i}$ also has the advantage that we get fewer factors of~$2\pi \i$ in the statements of the results below.

\begin{Remark}
As pointed out in \cite[Remark 2.3]{BHR}, the energy is defined
for all local sections around $\lambda = 0$ which admit a lift as in \eqref{Eq: LiftPowerseries}.
\end{Remark}

Let us write again \smash{$\mathcal{S}' = \mathcal S_{\mathcal{M}_{\mathrm{DH}}}'$} for the space of \emph{irreducible}
sections whose normal bundle is isomorphic to $\mathcal{O}_{\C P^1}(1)^{\oplus 2d}$. Take any $s \in
\mathcal S'$. In terms of lifts of sections, a tangent vector~${V\in T_s \mathcal{S}'}$ is expressed as
follows. Let $\widehat{s}$ be a lift of $s$ as in \eqref{Eq: LiftPowerseries}, and denote the
curvature of the connection~${\partial+{\overline{\partial}}}$ by \smash{$F^{\partial+{\overline{\partial}}} =
\overline{\partial}\partial + \partial\overline{\partial}$}. Expanding the integrability
condition
\begin{equation}\label{Eq: IntegrabilityCondition}
\overline{\partial}(\lambda)D(\lambda) + D(\lambda)\overline{\partial}(\lambda) = 0
\end{equation}
in powers of $\lambda$, the zeroth and first-order coefficients yield
\begin{align}\label{Eq: IntegrabilityConditionfirstorder}
\overline{\partial}\Phi = 0,\qquad
F^{\partial + {\overline{\partial}}} + [\Phi\wedge \Psi] = 0 .
\end{align}
Consider a family of sections \smash{$(s_t \in {\mathcal S}_{\mathcal{M}_{\mathrm{DH}}})_t$} with $s = s_0$ which
represents $V \in T_s \mathcal{S}'$.
The corresponding (lifted) infinitesimal variation \smash{$\dot{\widehat{s}} =
 \bigl(\lambda, \dot{{\overline{\partial}}}(\lambda), \dot{D}(\lambda)\bigr)$} satisfies the linearisation
of~\eqref{Eq: IntegrabilityCondition}, that is,
\begin{equation}\label{Eq: Intsdot}
{\overline{\partial}}(\lambda)\bigl(\dot D(\lambda)\bigr) + D(\lambda)\bigl(\dot{{\overline{\partial}}}(\lambda)\bigr) = 0.
\end{equation}
Expanding $\dot{\widehat{s}}$ into a power series
 \[
\dot{\widehat{s}}(\lambda) = \left(\sum_{k=0}^\infty\psi_k\lambda^k,
\sum_{k=0}^\infty \varphi_k\lambda^k, \lambda \right) ,
\]
for $\varphi_k \in \Omega^{1,0}(\mathfrak{sl}(E))$, $\psi_k \in \Omega^{0,1}(\mathfrak{sl}(E))$, the linearisation
of \eqref{Eq: IntegrabilityConditionfirstorder} becomes
\[
\overline{\partial}\varphi_0 + [\psi_0\wedge\Phi] = 0,\qquad
\overline{\partial}\varphi_1 + \partial\psi_0 + [\varphi_0\wedge \Psi] + [\Phi\wedge \psi_1] = 0 .
\]
Variations along the gauge orbit of $\widehat{s}$ are determined by infinitesimal gauge
transformations~${\C \ni \lambda \longmapsto \xi(\lambda) \in
\Gamma(\Sigma, \sln(E))}$ and are of the form
\begin{equation}\label{Eq:InfGauge}
\bigl(\lambda, {\overline{\partial}}(\lambda) \xi(\lambda), D(\lambda) \xi(\lambda)\bigr).
\end{equation}
By expanding $\xi(\lambda) = \sum_{k=0}^\infty \xi_k \lambda^k$, we get with \eqref{Eq:InfGauge} and \eqref{Eq: LiftPowerseries}
\begin{gather*}
{\overline{\partial}}(\lambda)\xi(\lambda) = {\overline{\partial}}\xi_0 + \bigl({\overline{\partial}}\xi_1 + [\Psi, \xi_0]\bigr)\lambda + O\bigl(\lambda^2\bigr), \\
D(\lambda)\xi(\lambda) = [\Phi, \xi_0] +(\partial\xi_0 + [\Phi, \xi_1])\lambda + O\bigl(\lambda^2\bigr) .
\end{gather*}
Now let $s \in \mathcal{S}'$ with lift $\widehat{s}$ over $\C$, and consider
$V_j \in T_s \mathcal{S}'$, $j = 1, 2$, represented by
\[
\dot{\widehat{s}}_j = \bigl(\lambda, \dot{{\overline{\partial}}}_j(\lambda), \dot{D}_j(\lambda)\bigr)
 = \bigl(\lambda, \psi^{(j)}_0 + \psi^{(j)}_1\lambda , \varphi_0^{(j)} + \varphi_1^{(j)}\lambda \bigr) + O\bigl(\lambda^2\bigr).
\]
Then we define, recalling the definition of $\omega_\lambda$ given in \eqref{eq:omegalambdaMDH},
\begin{align}
\widehat{\Omega}_{\widehat{s}}(V_1, V_2) & =
-\frac{\i}{2}\frac{\partial}{\partial\lambda}_{|\lambda=0}\omega_\lambda(V_1(\lambda), V_2(\lambda))\nonumber\\
& =-\frac{\i}{2} \frac{\partial}{\partial \lambda}_{|\lambda = 0} 2\i \int_\Sigma \operatorname{tr}\bigl(-\dot{D}_1(\lambda)
\wedge \dot{{\overline{\partial}}}_2(\lambda) + \dot{D}_2(\lambda) \wedge \dot{{\overline{\partial}}}_1(\lambda)\bigr)\nonumber \\
& = \int_\Sigma \operatorname{tr}\bigl(-\varphi_0^{(1)} \wedge \psi_1^{(2)} + \varphi_0^{(2)} \wedge \psi_1^{(1)} -
\varphi_1^{(1)} \wedge \psi_{0}^{(2)} + \varphi_1^{(2)} \wedge \psi_0^{(1)}\bigr).\label{o1}
\end{align}
We view $\widehat{\Omega}$ in \eqref{o1} as a two-form on the infinite-dimensional space of germs of
sections of $\varpi$ at~${\lambda = 0}$. Note that the formula for $\widehat{\Omega}$ is
exactly \eqref{eq:OmegaGeneral} in the present context. It is straightforward to check that the form
$\widehat{\Omega}$ is closed.

\begin{Proposition}\label{p:OmegaHol}
The two-form $\widehat{\Omega}$ in \eqref{o1} descends to a holomorphic two-form on the space of
irreducible sections, which on $\mathcal S_{\mathcal{M}_{\mathrm{DH}}}'$ coincides with the holomorphic
symplectic form $\Omega_0$ defined in \eqref{eq:OmegaGeneral}.
\end{Proposition}

\begin{proof}
Since the form $\widehat{\Omega}$ is closed, it suffices to show that the one-form at $\widehat s$ obtained by
contracting $\widehat{\Omega}_{\widehat{s}}$ with any vertical tangent vector at $\widehat{s}$ vanishes.

To this end, let $\widehat{s}$ be a germ of a section near $\lambda = 0$, and
let $\xi(\lambda) = \sum_{k=0}^\infty\xi_k\lambda^k$ be an infinitesimal gauge transformation.
The corresponding tangent vector $V_1$ is represented by
\[\dot{\widehat{s}}_1 = \bigl(\lambda, {\overline{\partial}}(\lambda)\xi(\lambda), D(\lambda)\xi(\lambda)\bigr).\]
Then for an arbitrary tangent vector $V_2$ represented by $\dot{\widehat{s}}_2 =
\bigl(\lambda, \dot{\overline{\partial}}(\lambda), \dot D(\lambda)\bigr)$, we find the following:
\begin{align*}
\widehat{\Omega}_{\widehat{s}}(V_1, V_2) & = \frac{\partial}{\partial \lambda}_{|\lambda = 0} \int_\Sigma\operatorname{tr}\bigl(-\dot D(\lambda)\wedge {\overline{\partial}}(\lambda)\xi(\lambda) + D(\lambda)\xi(\lambda)\wedge\dot{\overline{\partial}}(\lambda)\bigr) \\
\text{(Stokes)\qquad} & = \frac{\partial}{\partial \lambda}_{|\lambda = 0} \int_\Sigma\operatorname{tr}\bigl(\bigl({\overline{\partial}}(\lambda)\bigl(\dot D(\lambda)\bigr) + D(\lambda)\bigl(\dot{\overline{\partial}}(\lambda)\bigr)\bigr)\xi(\lambda)\bigr) = 0,
\end{align*}
where we used \eqref{Eq: Intsdot}.
This shows that $\widehat{\Omega}$ descends to $\mathcal S'$.
\end{proof}

Theorem~\ref{EnergyMomentMap} thus allows us to make the following conclusion.

\begin{Corollary}\label{cor:EnergyMDHC*}
The restriction of $2\pi \i\mathcal E\colon \mathcal S'_{\mathcal{M}_{\mathrm{DH}}}
\longrightarrow \C$ is a holomorphic moment map for the natural $\C^*$-action on $\mathcal S'_{\mathcal{M}_{\mathrm{DH}}}$ with
respect to the holomorphic symplectic form $\Omega_0$. In particular, the
$\C^*$-fixed points in $\mathcal S'_{\mathcal{M}_{\mathrm{DH}}}$ are exactly the critical points of \smash{$\mathcal E\big\vert_{\mathcal S'_{\mathcal{M}_{\mathrm{DH}}}}$}.
\end{Corollary}

\subsection[Explicit description of some C\^{}*-fixed sections]{Explicit description of some $\boldsymbol{\C^*}$-fixed sections}\label{ss:explicitC*}

Corollary \ref{cor:EnergyMDHC*} shows a close relationship between $\C^*$-fixed points in $\mathcal
S_{\mathcal{M}_{\mathrm{DH}}}$ and the energy functional. We therefore examine the $\C^*$-fixed points more closely in this
section. Before explicitly determining the $\C^*$-fixed \emph{irreducible} sections, we first
observe the following.

\begin{Lemma}\label{lem:C*fixedMdh}
The locus $\mathcal{S}_{\mathcal{M}_{\mathrm{DH}}}^{\C^*} \subset \mathcal{S}_{\mathcal{M}_{\mathrm{DH}}}$ of all $\C^*$-fixed sections is in a
natural bijection with~$\mathcal{M}_{\mathrm{dR}}$, the moduli space of flat completely reducible $\operatorname{SL}(n,\C)$-connections.

In particular, the critical points of $\mathcal{E} \colon \mathcal{S}_{\mathcal{M}_{\mathrm{DH}}}'
\longrightarrow \C$ correspond to an open subset of $\mathcal{M}_{\mathrm{dR}}^{\rm irr}$, the moduli space of flat irreducible $\operatorname{SL}(n,\C)$-connections.
\end{Lemma}

\begin{proof}
Take $\nabla \in \mathcal{M}_{\mathrm{dR}}$.
As in Section \ref{ss:C*}, we obtain the following $\C^*$-invariant section~$s_\nabla \colon \allowbreak\C^*
\longrightarrow \mathcal{M}_{\mathrm{DH}}$
\begin{equation}\label{eq:sNablaOverC*}
 s_\nabla(\lambda) = \big[\bigl(\lambda, {\overline{\partial}}^\nabla, \lambda \partial^\nabla\bigr) \big],\qquad
\partial^\nabla = \nabla^{1,0}, \qquad {\overline{\partial}}^\nabla = \nabla^{0,1}.
\end{equation}
By a very crucial result of Simpson (\cite{Si-Hodge} for existence and \cite{Simpson-iterated}
for a more explicit approach), the limits of $s_\nabla(\lambda)$ for
$\lambda \to 0$ and $\lambda\to \infty$
\emph{always} exist in $\mathcal{M}_{\mathrm{Higgs}}(\Sigma, \operatorname{SL}(n,\C))$ and $\mathcal{M}_{\mathrm{Higgs}}\bigl(\overline{\Sigma}, \operatorname{SL}(n,\C)\bigr)$ respectively.
The resulting section, also denoted by $s_{\nabla} \in \mathcal{M}_{\mathrm{DH}}$, is $\C^*$-invariant by continuity.
Evaluation of sections $s \colon \C P^1 \longrightarrow \mathcal{M}_{\mathrm{DH}}$ at $\lambda = 1$
gives the inverse of the map $\nabla \longmapsto s_\nabla$.

The last statement in the lemma is a direct consequence of Theorem~\ref{EnergyMomentMap} and
Corollary \ref{cor:EnergyMDHC*}.
\end{proof}

We next determine explicitly the $\C^*$-fixed sections $s \in \mathcal{S}_{\mathcal{M}_{\mathrm{DH}}}$ such that $s$
is irreducible over~$\C$; this will be done by using some results of \cite{CollierWentworth}. In terms of Lemma
\ref{lem:C*fixedMdh}, these are precisely the sections $s_\nabla$ such that $s_\nabla(0)$ is
stable. Indeed, since irreducibility is an open condition, $s_\nabla(\lambda)$ is an
irreducible $\lambda$-connection for all $\lambda$ close to $0$. Using the $\C^*$-invariance, we
see that $s(\lambda)$ is irreducible for every $\lambda \in \C$.

For any $\C^*$-fixed section $s_\nabla$, its values at $0$ and $\infty$ are $\C^*$-fixed Higgs
bundles on $\Sigma$ and $\overline{\Sigma}$ respectively.
These are called complex variations of Hodge structures (VHS).
Let $\bigl({\overline{\partial}}, \Phi\bigr)$ be any VHS on $\Sigma$.
The fact that $\bigl({\overline{\partial}}, \Phi\bigr)$ is a $\C^*$-fixed point yields a holomorphic splitting
\begin{equation}\label{Eq: VHSsplitting}
E = \bigoplus_{j=1}^lE_j
\end{equation}
into a direct sum of holomorphic vector bundles.
With respect to the splitting in \eqref{Eq: VHSsplitting}, ${\overline{\partial}}$ and~$\Phi$ are given in the
following block form:
\begin{equation}\label{Eq: VHSblockmatrices}
\overline{\partial} = \begin{pmatrix}
\overline{\partial}_{E_1} & 0 &\dots & \dots & 0\\
0 & \overline{\partial}_{E_2} & \ddots & & \vdots\\
\vdots & \ddots & \ddots & \ddots & \vdots\\
\vdots & & \ddots &\ddots & 0\\
0 & \dots & \dots & 0 & \overline{\partial}_{E_l}
\end{pmatrix},\qquad
\Phi = \begin{pmatrix}
0 & \dots & \dots & \dots & 0\\
\Phi^{(1)} & \ddots & & &\vdots\\
0 & \Phi^{(2)} & \ddots & &\vdots\\
\vdots & \ddots & \ddots & \ddots & \vdots\\
0 & \dots & 0 & \Phi^{(l-1)} & 0
\end{pmatrix},
\end{equation}
where $\Phi^{(j)} \in H^0(\Sigma, {\rm Hom}(E_j, E_{j+1})\otimes K_\Sigma)$.
The sheaf $\mathfrak{sl}(E)$ of trace-free holomorphic endomorphisms of $E$ further decomposes into
\[
\mathfrak{sl}(E) = \bigoplus_{k\in\mathbb{Z}} \mathfrak{sl}(E)_k,\qquad
\mathfrak{sl}(E)_k = \{\psi \in \mathfrak{sl}(E) \mid \psi(E_i) \subset E_{i-k} \}.
\]
By construction, we have $\Phi \in H^0(\Sigma, K_\Sigma\otimes \mathfrak{sl}(E)_{-1})$.
To define the next notion, let
\begin{equation}\label{eq:def_N+}
N_+ = \bigoplus_{k>0} \mathfrak{sl}(E)_k ,\qquad N_- = \bigoplus_{k<0 } \mathfrak{sl}(E)_k, \qquad
\mathbb{L} = \mathfrak{sl}(E)_0.
\end{equation}
Note that $N_+$ (resp.\ $N_-$) is the subspace of $\mathfrak{sl}(E)$ consisting of endomorphisms
of $E$ that are strictly upper (resp.\ lower) triangular form with blocks with respect to the
splitting in \eqref{Eq: VHSsplitting}, while $\mathbb{L}$ is the space of
elements of $\mathfrak{sl}(E)$ with diagonal blocks.

Now let $\bigl({\overline{\partial}}, \Phi\bigr) \in \mathcal{M}_{\rm Higgs}(\operatorname{SL}(n,\C))$ be a \emph{stable} complex variation of Hodge structures.
Then the BB-slice \cite[Definition~3.7]{CollierWentworth} through $\bigl({\overline{\partial}}, \Phi\bigr)$ is defined by
\begin{align}\label{eq:BB-slice}
\mathcal{B}^+_{({\overline{\partial}},\Phi)} =
\left\{ (\beta, \phi) \in \Omega^{0,1}(N_+)
\oplus \Omega^{1,0}(\mathbb{L}\oplus N_+)\,\bigg\vert\,
\begin{array}{@{}l@{}}
D''(\beta, \phi)+ [\beta\wedge\phi] = 0,\\
D'(\beta, \phi) = 0
\end{array}
\right\}.
\end{align}
Here we denote by
$
D := {\overline{\partial}}+\partial^h +\Phi+\Phi^{*_h}
$
the non-abelian Hodge connection associated to $\bigl({\overline{\partial}}, \Phi\bigr)$ with harmonic metric $h$, and set
$
D'' := {\overline{\partial}}+\Phi$, $ D' := \partial^h+\Phi^{*_h}$.
Hence the equations in \eqref{eq:BB-slice} are explicitly given by
\begin{gather}
D''(\beta, \phi) + [\beta\wedge\phi] = {\overline{\partial}} \phi + [(\Phi+\phi)\wedge\beta]
 = 0,\nonumber \\
D'(\beta, \phi) = \partial^{h}\beta +[\Phi^{*_h}\wedge\phi] = 0.\label{eq:BB-slice2}
\end{gather}
Note that $\mathcal{B}^+_{({\overline{\partial}},\Phi)}$ is a finite-dimensional affine space. Then, \cite[Theorem 1.4\,(3)]{CollierWentworth} states
that the map
\[
p \colon\ \mathcal{B}^+_{({\overline{\partial}},\Phi)}\times \C \longrightarrow
\mathcal{M}_{\mathrm{Hod}},\qquad ((\beta, \phi), \lambda) \longmapsto \big[\lambda, {\overline{\partial}} + \lambda \Phi^* +\beta,
 \lambda\partial^h+\Phi +\phi\big]
\]
is a holomorphic embedding onto the ``attracting set''
\[
W\bigl({\overline{\partial}}, \Phi\bigr) = \bigl\{ m \in \mathcal{M}_{\mathrm{Hod}}^{\rm irr} \mid \lim_{\zeta \to 0}\zeta\cdot m = \bigl({\overline{\partial}}, \Phi\bigr)\bigr\}
\]
and it is compatible with the obvious projections to $\C$.
In particular, if $W^\lambda\bigl({\overline{\partial}}, \Phi\bigr)$ denotes the intersection of $W\bigl({\overline{\partial}}, \Phi\bigr)$ with
the fiber $\varpi^{-1}(\lambda)$, then $W^\lambda\bigl({\overline{\partial}}, \Phi\bigr)$ is biholomorphic to the
affine space \smash{$\mathcal{B}^+_{({\overline{\partial}},\Phi)}$} via the map $p_{\lambda} := p(\bullet, \lambda)$.
Thus, \smash{$\mathcal{M}_{\mathrm{Hod}}^{\rm irr}$} is stratified by affine spaces.

Given \smash{$(\beta, \phi) \in \mathcal{B}^+_{({\overline{\partial}}, \Phi)}$}, we can use Lemma \ref{lem:C*fixedMdh} and \eqref{eq:sNablaOverC*} to define the $\C^*$-fixed section
\[
s_{(\beta, \phi)} := s_{p_1(\beta, \phi)} \in \mathcal S_{\mathcal{M}_{\mathrm{DH}}}.
\]
As observed earlier, $s_{(\beta, \phi)}$ is an irreducible section over $\C \subset \C P^1$
but not necessarily over all of~$\C P^1$.

\begin{Proposition}\label{Prop:C*_sectionLiftC}
Over $\C$, the $\C^*$-fixed section $s_{(\beta, \phi)}$ may be expressed as
\begin{equation}\label{eq:C*_sectionLiftC}
s_{(\beta,\phi)}(\lambda) = \left[\lambda, {\overline{\partial}} + \lambda(\Phi^{*_h}+\beta_1) + \sum_{j=2}^l\lambda^j\beta_j, \Phi + \lambda\partial^h + \sum_{j=0}^l\lambda^{j+1}\phi_j\right],
\end{equation}
where \smash{$\beta = \sum_{j=1}^l\beta_j$}, with $\beta_j \in \Omega^{0,1}(\mathfrak{sl}(E)_j)$ and
$\phi = \sum_{j=0}\phi_j$ with $\phi_j \in \Omega^{1,0}(\mathfrak{sl}(E)_j)$.
\end{Proposition}

\begin{proof}
Let
$\nabla = p_1(\beta, \phi) = [D+\beta+\phi]$
so that
\smash{$
{\overline{\partial}}^\nabla = {\overline{\partial}} + \Phi^{*_h} + \beta$}, $ \partial^\nabla = \partial^h + \Phi + \phi$.
Hence $s_{(\beta, \phi)} = s_\nabla$ is given by
\[
s_{(\beta,\phi)}(\lambda) = \big[\lambda, {\overline{\partial}}+(\Phi^{*_h}+\beta) ,
 \lambda \partial^h+\lambda\Phi+\lambda \phi \big]
\]
for $\lambda \in \C^*$ (see \eqref{eq:sNablaOverC*}).
This does \emph{not} give a lift of \smash{$s_{(\beta,\phi)}$} over all of $\C$, unless the
holomorphic bundle $\bigl(E, {\overline{\partial}}\bigr)$ is stable, in which case we must have $\beta = 0$ and
$\Phi = 0$.

To construct a lift over all of $\C$, we use the $\C^*$-family of gauge transformations
\begin{equation}\label{Eq: glambdaVHS}
g(\lambda) = \lambda^{m}\begin{pmatrix}
\lambda^{1-l}\mathrm{id}_{E_1} & 0 &\dots & \dots & 0\\
0 & \lambda^{2-l}\mathrm{id}_{E_2} & \ddots & & \vdots\\
\vdots &0& \ddots & \ddots & \vdots\\
\vdots & & \ddots &\ddots & 0\\
0 & \dots & \dots & 0 & \lambda^0\mathrm{id}_{E_l}
\end{pmatrix},
\end{equation}
where \smash{$m = \frac{1}{n}\sum_{j=1}^l(l-j)\mathrm{rk}(E_j)$}, in order to ensure
that $\det g(\lambda) = 1$. Then any $\xi \in \mathfrak{sl}(E)_j$ satisfies
\smash{$
g(\lambda)^{-1}\xi g(\lambda) = \lambda^j\xi$}.
Let $\beta = \sum_{j=1}^l\beta_j$, with $\beta_j \in \Omega^{0,1}(\mathfrak{sl}(E)_j)$, and similarly
$\phi = \sum_{j=0}\phi_j$ with $\phi_j \in \Omega^{1,0}(\mathfrak{sl}(E)_j)$.
Then using $\Phi \in H^0(K\otimes\mathfrak{sl}(E)_{-1})$ and $\Phi^{*_h} \in
\Omega^{0,1}(K\otimes\mathfrak{sl}(E)_1)$, we get that
\begin{gather*}
\bigl({\overline{\partial}}+(\Phi^{*_h}+\beta) , \lambda \partial^h+\lambda\Phi+\lambda \phi\bigr).g(\lambda)\\
 \qquad= \left({\overline{\partial}} + \lambda(\Phi^{*_h}+\beta_1) + \sum_{j=2}^l\lambda^j\beta_j,
\Phi + \lambda\partial^h + \sum_{j=0}^l\lambda^{j+1}\phi_j\right).
\end{gather*}
The result follows.
\end{proof}

Next we discuss the implications for the $\C^*$-fixed leaves of the foliation $\mathcal{F}^+$
on $\mathcal{S}' = \mathcal{S}_{\mathcal{M}_{\mathrm{DH}}}'$. Recall that these leaves consist, in particular, of
\emph{irreducible} sections \big(on all of $\C P^1$\big), by de\-finition. We denote by
\smash{$\mathcal{S}_{({\overline{\partial}}, \Phi)}'$} all sections in $\mathcal{S}'$ which pass through the stable
complex variation~of Hodge~structure \smash{$\bigl({\overline{\partial}}, \Phi\bigr) \in \mathcal{M}_{\mathrm{Higgs}}^{\C^*}$} at $\lambda = 0$.

\begin{Proposition}\label{Prop:classifying_C*-fixed_sections}
The $\C^*$-fixed point locus \smash{$(\mathcal{S}'_{({\overline{\partial}}, \Phi)})^{\C^*}$} is isomorphic to an
open and non-empty subset of the affine space \smash{$\mathcal{B}^+_{({\overline{\partial}}, \Phi)}$}.
\end{Proposition}

\begin{proof}
Consider the section \smash{$s_{(\beta, \phi)} \colon \C P^1 \longrightarrow
 \mathcal{M}_{\mathrm{DH}}$} for \smash{$(\beta, \phi) \in \mathcal{B}^+_{({\overline{\partial}}, \Phi)}$} which is irreducible over~$\C$.
Because the complement of \smash{$\mathcal{M}_{\mathrm{Higgs}}^{\rm irr}\big(\overline{\Sigma}, \operatorname{SL}(n,\C)\big)$} in
$\mathcal{M}_{\mathrm{Higgs}}\big(\overline{\Sigma}, \operatorname{SL}(n,\C)\big)$ is closed and of codimension at least two
(cf.\ \cite{Faltings}), it follows that \smash{$s_{(\beta, \phi)}$} is an irreducible section~for \smash{$(\beta, \phi) \in \mathcal{B}^+_{({\overline{\partial}}, \Phi)}$} lying in an open and dense
subset of \smash{$\mathcal{B}^+_{({\overline{\partial}}, \Phi)}$}.

Note that $(\beta, \phi) = (0, 0)$ corresponds to the twistor line \smash{$s_{({\overline{\partial}}, \Phi)}$}
through $\bigl({\overline{\partial}}, \Phi\bigr)$, which lies in~$\mathcal{S}'$. Since $\mathcal{S}'$ is open and
non-empty in the space of all irreducible sections, we therefore see that the irreducible and
$\C^*$-fixed section $s_{(\beta, \phi)}$ has the desired normal bundle for $(\beta, \phi)$ in an
open and non-empty subset \smash{$U \subset \mathcal{B}^+_{({\overline{\partial}}, \Phi)}$}. Altogether, we obtain the
isomorphism
\[
p_1^{-1} \circ \operatorname{ev}_1 \colon\ (\mathcal{S}_{({\overline{\partial}}, \Phi)}')^{\C^*}
 \overset{\cong}{\longrightarrow} U.
\]
This completes the proof.
\end{proof}

From Theorem~\ref{EnergyMomentMap}, we immediately obtain the following.

\begin{Corollary}
The locus of critical points $s \in \mathcal{S}'$ of $\mathcal{E} \colon \mathcal{S}'
 \longrightarrow \C$ is isomorphic to an open and non-empty subset in $\mathcal{M}_{\mathrm{dR}}^{\rm irr}$.
It is foliated by leaves which are isomorphic to open and non-empty subsets of affine spaces.
\end{Corollary}

\begin{proof}
The first statement follows, invoking a genericity argument, from Lemma \ref{lem:C*fixedMdh}.
The second one is a consequence of Proposition \ref{Prop:classifying_C*-fixed_sections}.
\end{proof}

\begin{Remark}
Let $s\colon \C P^1 \longrightarrow \mathcal \mathcal{M}_{\mathrm{DH}}$ be a $\C^*$-fixed section such that $s(0) =
 \bigl({\overline{\partial}}, \Phi\bigr)$ and $s(\infty) = (\partial, \Psi)$ are \emph{stable} VHS on $\Sigma$ and
$\overline{\Sigma}$ respectively, with respective splittings, of the underlying smooth
vector bundle $E$, of the form
\[
E = \bigoplus_{j=1}^lE_j, \qquad E = \bigoplus_{j=1}^{l'}E'_j.
\]
With respect to these splittings, the respective holomorphic structures are diagonal and the
Higgs fields $\Phi$ and $\Psi$ are lower triangular as in \eqref{Eq: VHSblockmatrices}. Then we
have the BB-slices \smash{$\mathcal B^+_{({\overline{\partial}},\Phi)}(\Sigma)$} and~\smash{$\mathcal
B^+_{(\partial,\Psi)}(\overline{\Sigma})$}. By Proposition
\ref{Prop:classifying_C*-fixed_sections} and its analog on $\overline{\Sigma}$, we see that, on
the one hand, $s$ corresponds to \smash{$(\beta, \phi) \in \mathcal B^+_{({\overline{\partial}},\Phi)}(\Sigma)$}, and on
the other hand it corresponds to \smash{$\big(\widetilde\beta, \widetilde\phi\big) \in \mathcal
B^+_{(\partial,\Psi)}(\overline{\Sigma})$}. Therefore, we obtain two distinguished lifts of $s$
over $\C$ and $\C^*\cup\{\infty\}$ of the form
\begin{gather*}
s(\lambda) = \big[\lambda, \widehat{s}_{(\beta,\phi)}(\lambda)\big]_\Sigma\\
\phantom{s(\lambda)}{}=
\left[\lambda, {\overline{\partial}} + \lambda(\Phi^{*_h}+\beta_1) + \sum_{j=2}^l\lambda^j\beta_j, \Phi + \lambda\partial^h + \sum_{j=0}^l\lambda^{j+1}\phi_j\right]_\Sigma,
\\
s(\lambda) = \big[\lambda^{-1}, \widehat{s}_{(\widetilde\beta,\widetilde\phi)}\bigl(\lambda^{-1}\bigr)\big]_{\overline{\Sigma}}\\
\phantom{s(\lambda)}{} = \left[\lambda^{-1}, \partial + \lambda^{-1}\bigl(\Psi^{*_{\widetilde h}}+\widetilde\beta_1\bigr)
+ \sum_{j=2}^{l'}\lambda^{-j}\widetilde\beta_j, \Psi + \lambda^{-1}\overline{\partial}^{\widetilde h} + \sum_{j=0}^{l'} \lambda^{-(j+1)}\widetilde\phi_j\right]_{\overline{\Sigma}}.
\end{gather*}
Let $g_0$ be a gauge transformation such that
\[\bigl(\partial + \overline{\partial}^{\widetilde h} + \Psi + \Psi^{*_{\widetilde h}}+
\widetilde\beta +\widetilde\phi\bigr).g_0 = {\overline{\partial}} + \partial^h + \beta+\phi.\]
Going through the proof of Proposition \ref{Prop:C*_sectionLiftC} and writing $g(\lambda)$ and
$\widetilde g\bigl(\lambda^{-1}\bigr)$ for the respective $\C^*$-families of gauge transformations, we get that
\[
\widehat{s}_{(\beta,\phi)}(\lambda)
 = \widehat{s}_{(\widetilde\beta,\widetilde\phi)}\bigl(\lambda^{-1}\bigr).\widetilde{g}\bigl(\lambda^{-1}\bigr)^{-1}g_0g(\lambda)
\]
for any $\lambda \in \C^*$.
\end{Remark}

In general, starting only with the lift $\widehat{s}_{(\beta,\phi)}$ over $\C$ obtained above,
it seems hard to determine explicitly the lift
$\widehat{s}_{(\widetilde\beta,\widetilde\phi)}\bigl(\lambda^{-1}\bigr)$ over $\C P^1\setminus\{0\}$ or
even the limiting VHS $s_{(\beta,\phi)}(\infty)$. The next two examples discuss some situations
in which the limit can actually be computed.

\begin{Example}
Suppose the holomorphic structure $\partial^h + \Phi + \phi$ is stable on $\overline{\Sigma}$. Then we can argue as follows. For
$\lambda \in \C^*$ we can write, using the Deligne gluing
\begin{align*}
s_{(\beta,\phi)}(\lambda) &= \big[\lambda,{\overline{\partial}}+\Phi^{*_h}+\beta, \lambda\bigl(\partial^h + \Phi + \phi\bigr)\big]_{\Sigma} = \big[\lambda^{-1},\partial^h + \Phi + \phi,\lambda^{-1}\bigl({\overline{\partial}}+\Phi^{*_h}+\beta\bigr)\big]_{\overline{\Sigma}}.
\end{align*}
Under our assumption that $\partial^h + \Phi + \phi$ is stable, this allows us to conclude $s_{(\beta,\phi)}(\infty)
 = \bigl(\partial^h + \Phi + \phi,0\bigr)$. We will see in the proof of Theorem~\ref{Thm: Hyperholnontrivial} that this situation does in
fact occur, at least for rank $2$ bundles.
\end{Example}

\begin{Example}
Consider the rank two case meaning $n = 2$. If $s$ is the twistor line through a~VHS $\bigl({\overline{\partial}}, \Phi\bigr)$
on $\Sigma$, then we have $E = V \oplus V^*$, where $V$ is a line bundle with $0 < \deg V
 \leq g-1$ and $V^* = \ker\Phi$. Then $s(\infty) = \bigl(\partial^h, \Phi^{*_h}\bigr)$ and the corresponding
splitting is $E = V^*\oplus V$. Note that, since $\overline{\Sigma}$ and $\Sigma$ come with opposite
orientations, we have $\deg V^* > 0$, as a vector bundle on~$\overline{\Sigma}$. Then
$\widetilde g\bigl(\lambda^{-1}\bigr) = g(\lambda)$ in this case, as the order is reversed. The associated lifts
are thus just the lifts of $s$ over $\C$ and $\C^*$ given by the harmonic metric, that is, the associated
solution of the self-duality equations.
\end{Example}

\begin{Example}[{grafting sections}]\label{ex:grafting}
In \cite{Hel} a special class of $\C^*$-invariant sections of $\mathcal{M}_{\mathrm{DH}}(\Sigma, \allowbreak\operatorname{SL}(2,\C))$, called \emph{grafting sections}, have been constructed by using grafting of projective structures on $\Sigma$.
We recover them from the previous proposition as follows.
	
Consider the $\C^*$-fixed stable Higgs bundle $\bigl({\overline{\partial}}, \Phi\bigr)$ with
\[
E = K_\Sigma^{\frac{1}{2}}\oplus K^{-\frac{1}{2}}_\Sigma, \qquad
\Phi =
\begin{pmatrix}
0 & 0 \\
1 & 0
\end{pmatrix},
\]
where \smash{$K_\Sigma^{\frac{1}{2}}$} is a square root of the canonical bundle $K_\Sigma$.
To determine \eqref{eq:def_N+} in this example, we define \smash{$E_1:=K_\Sigma^{\frac{1}{2}}$}, \smash{$E_2:=K_\Sigma^{-\frac{1}{2}}$}.
Then we see that
$
N_+\cong K_\Sigma$, $ N_-=K^{-1}_\Sigma$, $ L \cong \mathcal{O}_\Sigma$. From \eqref{eq:BB-slice2} it follows
that \smash{$(0, \phi) \in \mathcal{B}_{({\overline{\partial}},\Phi)}^+$} if and only if ${\overline{\partial}} \phi = 0$ and
$[\phi\wedge\Phi^{*_h}] = 0$.
Hence $\phi$ is of the form
\[
\phi = \begin{pmatrix}
0 & q \\
0 & 0
\end{pmatrix},\qquad q \in H^0\big(\Sigma, K_\Sigma^{\otimes 2}\big),
\]
with respect to the splitting $E = E_1\oplus E_2$.
For those $q$ such that the monodromy of the corresponding flat connection at $\lambda = 1$ is real, the
sections $s_{(0,\phi)}$ are precisely the grafting sections of~\cite[Section~2.1]{Hel}. Since $\beta = 0$ in
this case, we see that the energy of a grafting section is the same as the energy of the twistor line
associated with the stable Higgs pair~$\bigl({\overline{\partial}}, \Phi\bigr)$. If the monodromy of the corresponding flat
connection is real, then \cite{Hel} shows that the section~$s_{(0,\phi)}$ is real and
it defines an element of \smash{$(\mathcal S'_{\mathcal{M}_{\mathrm{DH}}})^\tau$}, in particular it has the correct normal bundle~\smash{$\mathcal{O}_{\C P^1}(1)^{\oplus 2d}$}. But the section $s_{(0,\phi)}$ is not admissible and thus cannot correspond
to a~solution of the self-duality equations. This shows that we have
\[\mathcal{M}_{\mathrm{SD}}(\Sigma,\operatorname{SL}(2,\C)) \subsetneq (\mathcal S_{\mathcal{M}_{\mathrm{DH}}}')^\tau.\]
\end{Example}

\subsection[The energy of a C\^{}*-fixed section]{The energy of a $\boldsymbol{\C^*}$-fixed section}

Proposition \ref{Prop:C*_sectionLiftC} gives concrete formulas for all $\C^*$-fixed points
$s \in \mathcal{S}_{\mathcal{M}_{\mathrm{DH}}}^{\C^*}$ such that $s(0)$ is a~stable VHS.
We next compute the energy of such sections.

\begin{Proposition}\label{Prop: energyC*-fixed_sections}
Let $\bigl({\overline{\partial}}, \Phi\bigr)$ be a stable $\C^*$-fixed $\operatorname{SL}(n,\C)$-Higgs bundle, and let $s_{(\beta,\phi)}$
be the $\C^*$-fixed section corresponding to \smash{$(\beta, \phi) \in \mathcal{B}^+_{({\overline{\partial}},\Phi)}$}.
Its energy is given by
\[
\mathcal{E}(s_{(\beta,\phi)}) = \mathcal{E}(s_0) = \sum_{k=2}^l(k-1)\deg(E_k),
\]
where $s_0$ is the twistor line through $\bigl({\overline{\partial}}, \Phi\bigr)$.
\end{Proposition}

\begin{proof}
Write $s_{(\beta,\phi)}$ in a form as in \eqref{eq:C*_sectionLiftC}. Then from the definition of
$\mathcal{E}$ it follows immediately that
\[
\mathcal{E}(s_{(\beta,\phi)}) = \mathcal{E}(s_0)+\frac{1}{2\pi \i} \int_\Sigma \operatorname{tr}(\Phi \wedge \beta_1).
\]
Next we will show that $\int_\Sigma \operatorname{tr}(\Phi \wedge \beta_1) = 0$. To this end, let us write
\[
\Phi = \sum_{k=1}^{l-1}\Phi^{(k)},\qquad \beta_1 = \sum_{k=1}^{l-1}\beta^{(k)},
\]
where $\Phi^{(k)} \in \Omega^{1,0}(\operatorname{Hom}(E_k, E_{k+1}))$, $ \beta^{(k)} \in
\Omega^{0,1}(\operatorname{Hom}(E_{k+1}, E_k))$; see the block form in \eqref{Eq: VHSblockmatrices}. It follows that
\[
\operatorname{tr}(\Phi\wedge\beta_1) = \sum_{k=1}^l\operatorname{tr}_{E_k}\bigl(\Phi^{(k-1)}\wedge\beta^{(k-1)}\bigr).
\]
Note that each of the above summands $\Phi^{(k-1)}\wedge\beta^{(k-1}$ belongs to
$\Omega^{1,1}(\mathrm{End}(E_k))$, and we have adopted the convention that $\Phi^{(k)}
= 0 = \beta^{(k)}$ if $k = 0, l$.

Now, equation \eqref{eq:BB-slice2} implies that
$
{\overline{\partial}} \phi_0 + [\Phi\wedge\beta_{1}] = 0
$
and we can write
\[
[\Phi\wedge\beta_1] = \sum_{k=1}^{l-1}\Phi^{(k-1)}\wedge\beta^{(k-1)} + \beta^{(k)}\wedge\Phi^{(k)}.
\]
Thus, for each $k = 1, \dots, l$,
\[
{\overline{\partial}} \phi_0^{(k)} + \Phi^{(k-1)}\wedge\beta^{(k-1)} + \beta^{(k)}\wedge\Phi^{(k)} = 0.
\]
Consider the case of $k = l$
\[
{\overline{\partial}} \phi_0^{(l)} + \Phi^{(l-1)}\wedge\beta^{(l-1)} = 0.
\]
Taking the trace of this equation and integrating over $\Sigma$, we find, using Stokes' theorem,
that
\[
\int_\Sigma\operatorname{tr}_{E_l}\bigl(\Phi^{(l-1)}\wedge\beta^{(l-1)}\bigr) = 0.
\]
Now assume that $\int_\Sigma\operatorname{tr}_{E_{k+1}}\bigl(\Phi^{(k)}\wedge\beta^{(k)}\bigr) = 0$
for all $k \geq k_0$. Then we have
\[
{\overline{\partial}} \phi_0^{(k_0)} + \Phi^{(k_0-1)}\wedge\beta^{(k_0-1)} + \beta^{(k_0)}\wedge\Phi^{(k_0)} = 0.
\]
Taking the trace and integrating yields
\begin{align*}
0 &= \int_\Sigma\operatorname{tr}_{E_{k_0}}\bigl(\Phi^{(k_0-1)}\wedge\beta^{(k_0-1)} + \beta^{(k_0)}\wedge\Phi^{(k_0)}\bigr) \\
&= \int_\Sigma\operatorname{tr}_{E_{k_0}}\bigl(\Phi^{(k_0-1)}\wedge\beta^{(k_0-1)}\bigr) - \int_\Sigma\operatorname{tr}_{E_{k_0+1}}\bigl(\Phi^{(k_0)}\wedge\beta^{(k_0)}\bigr) \\
&= \int_\Sigma\operatorname{tr}_{E_{k_0}}\bigl(\Phi^{(k_0-1)}\wedge\beta^{(k_0-1)}\bigr).
\end{align*}
It follows inductively that $\int_\Sigma\operatorname{tr}(\Phi\wedge\beta_1) = 0$.

It remains to compute the energy of the twistor line $s_0$. To this end, we observe that
\[
\mathcal{E}(s_0) = \frac{1}{2\pi \i}\int_\Sigma \operatorname{tr}(\Phi\wedge \Phi^{*_h}) =
 \frac{1}{2\pi \i}\int_\Sigma\sum_{k=2}^l \operatorname{tr}_{E_k}\bigl(\Phi^{(k-1)}\wedge \bigl(\Phi^{k-1}\bigr)^{*_h}\bigr)
 = \sum_{k=2}^l\mathcal E_k(s_0),
\]
where we put \smash{$\mathcal E_k(s_0) =
\frac{1}{2\pi \i}\int_\Sigma \operatorname{tr}_{E_k}\bigl(\Phi^{(k-1)}\wedge \bigl(\Phi^{k-1}\bigr)^{*_h}\bigr)$}
for $k \geq 2$.
The equation $F^{\nabla^h} + {[\Phi\wedge\Phi^{*_h}] = 0}$ is in the form of diagonal blocks, with respect
to the splitting $E = \bigoplus_{k=1}^lE_k$, with components
\[
F^{\nabla^h_{E_k}} + \Phi^{(k-1)}\wedge \bigl(\Phi^{(k-1)}\bigr)^{*_h} + \bigl(\Phi^{(k)}\bigr)^{*_h}\wedge \Phi^{(k)} = 0.
\]
This gives the following recursive relations:
\begin{align*}
\mathcal E_k(s_0) & = \frac{1}{2\pi \i}\int_\Sigma\operatorname{tr}_{E_k}\bigl(\Phi^{(k-1)}\wedge \bigl(\Phi^{(k-1)}\bigr)^{*_h}\bigr)\\
&= \frac{\i}{2\pi}\int_\Sigma\operatorname{tr}_{E_k}\bigl(F^{\nabla^h_{E_k}}\bigr) +\frac{1}{2\pi \i}\int_\Sigma\operatorname{tr}_{E_{k+1}}\bigl(\Phi^{(k)}\wedge\bigl(\Phi^{(k)}\bigr)^{*_h}\bigr)= \deg(E_k) + \mathcal E_{k+1}(s_0).
\end{align*}
Thus, if $k = l$, we find that
$
\mathcal E_l(s_0) = \deg(E_l) $,
and for general $k$ we get that
\[
\mathcal E_k(s_0) = \sum_{j=k}^{l-1}\deg(E_j) + \mathcal E_l(s_0) = \sum_{j=k}^l\deg(E_j).
\]
Therefore,
\[
\mathcal E(s_0) = \sum_{k=2}^l\mathcal E_k(s_0) = \sum_{k=2}^l\sum_{j=k}^l\deg(E_j)
 = \sum_{k=2}^l(k-1)\deg(E_k),
\]
and this completes the proof.
\end{proof}

\subsection[The second variation of the Energy at a C\^{}*-fixed section]{The second variation of the Energy at a $\boldsymbol{\C^*}$-fixed section}

Next we study the second variation of the energy functional $\mathcal E$ at a $\C^*$-fixed point.

Examining the proof of Proposition \ref{Prop:C*_sectionLiftC}, we can check explicitly that the sections $s_{(\beta,\phi)}$ satisfy,
for any $\zeta \in \C^*$, the relation
\smash{$\zeta.\widehat{s}_{(\beta,\phi)} = \widehat{s}_{(\beta,\phi)}.g(\zeta)^{-1}$}.
Moreover, if we use the notation of equation~\eqref{Eq: glambdaVHS} and put
\[
\xi = \begin{pmatrix}
(m+1-l)\mathrm{id}_{E_1} & 0 &\dots & \dots & 0\\
0 & (m+2-l)\mathrm{id}_{E_2} & \ddots & & \vdots\\
\vdots &0& \ddots & \ddots & \vdots\\
\vdots & & \ddots &\ddots & 0\\
0 & \dots & \dots & 0 & m\mathrm{id}_{E_l}
\end{pmatrix},
\]
then $[\xi, \cdot]$ acts as multiplication by $k$ on $\mathfrak{sl}(E)_k$ and we see that
\begin{equation}\label{eq:Dlambda}
-\i\lambda\frac{\rm d}{{\rm d}\lambda}\overline{\partial}(\lambda) = \overline{\partial}(\lambda)\xi(\lambda),\qquad
\i D(\lambda) - \i\lambda\frac{\rm d}{{\rm d}\lambda}D(\lambda)= D(\lambda)\xi(\lambda)
\end{equation}
with
\[\bigl(\overline{\partial}(\lambda), D(\lambda)\bigr) =
\left({\overline{\partial}} + \lambda(\Phi^{*_h}+\beta_1) + \sum_{j=2}^l\lambda^j\beta_j, \Phi + \lambda\partial^h + \sum_{j=0}^l\lambda^{j+1}\phi_j\right).\]
For $\xi(\lambda) = \sum_{k = 0}^\infty \xi_k \lambda^k$, we deduce from \eqref{eq:Dlambda} the following equations
\begin{gather}\label{eq:PhiPsi}
0 = {\overline{\partial}} \xi_0 ,\qquad \Phi = [\Phi, \xi_0] ,\qquad
-\Psi = {\overline{\partial}} \xi_1 + [\Psi, \xi_0] ,\qquad 0 = [\Phi, \xi_1] + \partial \xi_0.
\end{gather}
We can now compute the second variation of $\mathcal E$ at such fixed points.

\begin{Proposition}\label{prop1}
The second variation of $\mathcal E$ at a $\C^*$-fixed point $s$ with lift $\widehat{s}$ as in \eqref{Eq: LiftPowerseries} is given by
\[
{\rm d}^2\mathcal{E}(\dot s)= \frac{1}{2\pi \i}\int_\Sigma\operatorname{tr}(\psi_0\wedge[\varphi_1,\xi] + \varphi_1\wedge[\psi_0,\xi] + \psi_1\wedge[\varphi_0,\xi]) +\varphi_0\wedge[\psi_1, \xi] + 2\varphi_0\wedge\psi_1) .
\]
\end{Proposition}
\begin{proof}
Let $(s_t)$ be a family of sections with $s_0 = s$. We compute, using the notation for~$\widehat{s}$ and~$\dot s$ as in Section \ref{ss:GeometrySectionDH},
\[
2\pi \i\frac{{\rm d}^2}{{\rm d}t^2}|_{t=0}\mathcal E(s_t) =
\int_\Sigma\operatorname{tr}\bigl(\Phi\wedge\dot\psi_1 + \dot\varphi_0\wedge\Psi + 2\varphi_0\wedge\psi_1\bigr).
\]
Since the section $s_0 = s$ is fixed by the action of $\C^*$, we can use \eqref{eq:PhiPsi} (with $\xi
= \xi_0$, $\xi_1 = 0$) to write
\begin{align*}
2\pi \i\frac{{\rm d}^2}{{\rm d}t^2}|_{t=0}\mathcal E(s_t) ={}&
\int_\Sigma\operatorname{tr}\bigl([\Phi,\xi]\wedge\dot\psi -\dot\phi\wedge[\Psi,\xi] + 2\phi\wedge\psi\bigr)\\
={}& \int_\Sigma\operatorname{tr}\bigl(-\xi\bigl(\big[\Phi\wedge\dot\psi\big] + \big[\dot\phi\wedge\Psi\big]\bigr) + 2\phi\wedge\psi\bigr)\\
={}& \int_\Sigma\operatorname{tr}\bigl(\xi\bigl(\overline{\partial}\dot\varphi_1 +
\partial\dot\psi_0 + 2[\psi_0\wedge\varphi_1] + 2[\varphi_0\wedge \psi_1]\bigr) + 2\varphi_0\wedge\psi_1\bigr)\\
({\rm using}~\overline{\partial}\xi = 0 = \partial\xi) ={}& \int_\Sigma\operatorname{tr}\bigl(\xi\left(2[\psi_0\wedge\varphi_1] + 2[\varphi_0\wedge \psi_1]\right) + 2\varphi_0\wedge\psi_1\bigr)\\
={}& \int_\Sigma\operatorname{tr}\bigl(\psi_0\wedge[\varphi_1,\xi] + \varphi_1\wedge[\psi_0\wedge\xi] + \psi_1\wedge[\varphi_0,\xi] \\
&+\varphi_0\wedge[\psi_1\wedge\xi] + 2\varphi_0\wedge\psi_1\bigr).
\end{align*}
In the third equation from above, we made use of the second linearisation of \eqref{Eq: IntegrabilityConditionfirstorder}.
\end{proof}

Proposition \ref{prop1} shows that
the second variation is closely related to the infinitesimal $\C^*$-action on the tangent space. The following proposition is obtained.

\begin{Proposition}
Let $\dot s(\lambda) = \bigl(\lambda, \dot{\overline{\partial}}(\lambda), \dot D(\lambda)\bigr)
 = \bigl(\lambda, \sum_{k=0}^\infty\psi_k\lambda^k, \sum_{k=0}^\infty\varphi_k\lambda^k\bigr)$ be an infinitesimal deformation
of the critical point $s \in \mathcal S$. Suppose that $\dot s$ satisfies
\[ [\psi_0, \xi] = n_0\psi_0,\qquad [\psi_1, \xi] = n_1\psi_1,\qquad [\varphi_0, \xi]
 = m_0\varphi_0,\qquad [\varphi_1, \xi] = m_1\varphi_1\]
for some $m_i, n_i \in \mathbb{Z}$. Then
\[{\rm d}^2\mathcal{E}(\dot s) = \frac{1}{2\pi \i}\int_\Sigma\operatorname{tr}((m_1+n_0)\psi_0\wedge\varphi_1 + (m_0+n_1+2)\psi_1\wedge\varphi_0).\]
\end{Proposition}

\begin{Remark}
Note that this resembles the discussion surrounding \cite[equation (8.10)]{HitchinLie}. In fact, it does reproduce Hitchin's result in the
case where $s$ is the twistor line corresponding to a~$\C^*$-fixed point in $\mathcal \mathcal{M}_{\mathrm{Higgs}} $
and the deformation $\dot s$ is real, so that $\psi_1 = \varphi_0^*, \psi_0 = -\varphi_1^*$.
\end{Remark}

\subsection{Sections and the degree of the hyperholomorphic line bundle}
Our previous results together with the energy can be used to show that the space of irreducible sections is not connected.
We begin with the following proposition.

\begin{Proposition}\label{Prop: deghyperhol_C*fixedsection}
Let $\bigl({\overline{\partial}}, \Phi\bigr)$ be a stable $\C^*$-fixed Higgs bundle, and let
\smash{$s_{(\beta,\phi)}$} be a $\C^*$-fixed section corresponding to \smash{$(\beta, \phi) \in \mathcal{B}^+_{({\overline{\partial}},\Phi)}$}.
If $s_{\beta,\phi}(\infty)$ is given by a VHS on $\overline{\Sigma}$ with under\-lying~holomorphic bundle \smash{$E
 = \bigoplus_{k=1}^{l'}E'_k$}, then the following holds
\[
\deg(s_{(\beta,\phi)}^*L_Z) = \sum_{k=1}^l(k-1)\deg(E_k) + \sum_{k=1}^{l'}(k-1)\deg(E_k').
\]
\end{Proposition}

\begin{proof}
Proposition \ref{Prop: energyC*-fixed_sections} allows us to compute $\mathcal E_0(s)$ and $\mathcal E_\infty(s)$.
The assertion now follows from the formula
\[\deg(s_{(\beta,\phi)}^*L_Z) = \mathcal E(s_{(\beta,\phi)})+\mathcal E_\infty(s_{(\beta,\phi)}).\]
This completes the proof.
\end{proof}

\begin{Theorem}\label{Thm: Hyperholnontrivial}
There exist irreducible sections $s$ of $\varpi\colon \mathcal{M}_{\mathrm{DH}} (\Sigma,\operatorname{SL}(2,\C)) \longrightarrow \C P^1$ such that the pullback $s^* L_Z$ of the holomorphic line bundle
$L_Z \longrightarrow \mathcal{M}_{\mathrm{DH}}(\Sigma,\operatorname{SL}(2,\C))$ has non-zero degree. In particular, the space of irreducible sections is not connected.
\end{Theorem}

\begin{proof}
Let \smash{$K_\Sigma^{\frac{1}{2}}$} be a square-root of the canonical line bundle $K_\Sigma$.
Consider the uniformization (Fuchsian) flat connection
\[\nabla^{\rm Fuchs}=\begin{pmatrix} \nabla^{K_\Sigma^{\frac{1}{2}}} & 1^*\\ 1 & \nabla^{K_\Sigma^{-\frac{1}{2}}} \end{pmatrix}\]
on the rank two vector bundle \smash{$K_\Sigma^{\frac{1}{2}}\oplus K_\Sigma^{-\frac{1}{2}}$}.
For a generic holomorphic quadratic differential~${q \in H^0\bigl(\Sigma, K_\Sigma^2\bigr)}$, the anti-holomorphic structure
\[\begin{pmatrix} \partial_{K_\Sigma^{\frac{1}{2}}} & q\\ 1 & \partial_{K_\Sigma^{-\frac{1}{2}}} \end{pmatrix}\]
is stable (that is, it defines a stable holomorphic vector bundle on $\overline\Sigma$).
Then,
\[\nabla := \nabla^{\rm Fuchs}+\begin{pmatrix}0&q\\0&0\end{pmatrix}\]
gives a $\C^*$-invariant section $s_\nabla \in \mathcal S_{\mathcal{M}_{\mathrm{DH}}}$ by the construction of Lemma~\ref{lem:C*fixedMdh}.
In view of Proposition~\ref{Prop: energyC*-fixed_sections},
the energy at $\lambda = 0$ is given by
\[\deg{\bigl(K_\Sigma^{-\frac{1}{2}}\bigr)} = 1-g \neq 0 .\]
By assumption, $\partial^\nabla$ is stable, so the anti-Higgs field of $s$ at~${\lambda = \infty}$
vanishes, and the energy at~${\lambda=\infty}$ is given by $\mathcal E_\infty = 0.$
Finally, we have
\[\deg(s^*L_Z) = \mathcal E(s)+\mathcal E_\infty(s) \neq 0\]
by the residue formula for the pullback, under $s$, of the meromorphic connection to $\C P^1$ (see~\cite[Section~3]{BeHeRo}).
\end{proof}

\begin{Remark}
By the general theory explained in Section~\ref{ss:rotating}, the $C^\infty$ line bundle underlying~${L_Z \longrightarrow \mathcal{M}_{\mathrm{DH}}(\Sigma,\operatorname{SL}(n,\C))}$ agrees on $\mathcal{M}_{\mathrm{DH}}^{\rm irr}(\Sigma,\operatorname{SL}(n,\C))$ with the pullback of a line bundle on~$\mathcal{M}_{\mathrm{SD}}^{\rm irr}(\Sigma,\operatorname{SL}(n,\C))$ under the $C^\infty$ projection $\mathcal{M}_{\mathrm{DH}}^{\rm irr}(\Sigma,\operatorname{SL}(n,\C))\cong S^2\times \mathcal{M}_{\mathrm{SD}}^{\rm irr}(\Sigma,\operatorname{SL}(n,\C))\to \mathcal{M}_{\mathrm{SD}}^{\rm irr}(\Sigma,\operatorname{SL}(n,\C))$. Since the twistor lines are just the fibres of this projection, it follows that~$\deg(s^*L_Z) =0$ for any twistor line $s$.
\end{Remark}

Given a(n irreducible) section $s\in \mathcal{S}_{\mathcal{M}_{\mathrm{DH}}}$, it is in general very difficult to compute its normal bundle $N_s$.
However,
by using the methods of \cite{Hel}, it can be shown that the $\C^*$-fixed points considered in the proof of
Theorem~\ref{Thm: Hyperholnontrivial} do not have
normal bundles of generic type, i.e., their normal bundle admits holomorphic sections with double zeros.

\subsection*{Acknowledgements}
 We thank Danu Thung for pointing out a number of typos and errors in an earlier version of this manuscript. We thank the anonymous referees for their careful reading of the manuscript.
F.B.\ is supported by the DFG Emmy Noether grant AL 1407/2-1. I.B.\ is partially supported by a J.C.~Bose
Fellowship (JBR/2023/000003). S.H.\ is supported by the DFG grant HE 6829/3-1 of the DFG priority program SPP 2026 Geometry at
Infinity.

\pdfbookmark[1]{References}{ref}
\LastPageEnding

\end{document}